\documentclass[dvips]{amsart}
\usepackage[dvips]{graphics}
\usepackage{amscd,amssymb,amsxtra,latexsym,pictex,epsfig,psfrag}
\newcommand{\vanish}[1]{}

\newcommand{\R}{\mathbb{R}}
\newcommand{\Q}{\mathbb{Q}}

\newcommand{\Z}{\mathbb{Z}}

\newcommand{\id}{\operatorname{id}}

\newcommand{\pointsat}{\mathrel{\to\mkern-10mu-}}

\newcommand{\degg}{\operatorname{deg}}
\newcommand{\sign}{\operatorname{sign}}
\newcommand{\Aut}{\operatorname{Aut}}
\newcommand{\Der}{\operatorname{Der}}

\newcommand{\into}{\hookrightarrow}
\newcommand{\iso}{\cong}
\newcommand{\yright}{\stackrel{y}{\longrightarrow}}
\newcommand{\yleft}{\stackrel{y}{\longleftarrow}}
\newcommand{\isoto}{\stackrel{\cong}{\longrightarrow}}
\newcommand{\app}{\stackrel{(R2)}{\rightsquigarrow}}
\newcommand{\fsp}{{{\mathfrak{sp}}(2n)}}

\newcommand{\lhor}{\rule[.7ex]{.2in}{.2pt}}%horiz. line for lattice
\theoremstyle{plain}
\newtheorem*{corollary}{Corollary}
\newtheorem{lemma}{Lemma}
\newtheorem{theorem}{Theorem}
\newtheorem*{theorem*}{Theorem}

\newtheorem*{reduction}{Reduction lemma}
\newtheorem{proposition}{Proposition}
\theoremstyle{definition}
\newtheorem{definition}{Definition}
\newtheorem*{intuition}{Intuition}

\theoremstyle{remark}
\newtheorem*{remark}{Remark}

\newcommand{\adams}{adams78:_infin}
\newcommand{\bm}{bar-natan01:_graph}
\newcommand{\bl}{bocklandt:_neckl_lie}
\newcommand{\bll}{bergeron98:_combin}
\newcommand{\boothby}{boothby86:_rieman}
\newcommand{\kbrown}{brown94:_cohom}
\newcommand{\cv}{conant:_infin_operat_graph_graph_homol}
\newcommand{\drinfeld}{drinfelprimed90:_hopf_gal_q_q}

\newcommand{\fh}{fulton91:_repres}
\newcommand{\fuks}{fuks86:_cohom_lie}
\newcommand{\gelfand}{gelprimefand70:_cohom_lie}
\newcommand{\ferenc}{gerlits02:_calcul}
\newcommand{\getzler}{getzler94:_batal_vilkov}
\newcommand{\getzkap}{getzler98:_modul}
\newcommand{\getzlerkap}{getzler95:_cyclic}
\newcommand{\ginz}{ginzburg01:_non}
\newcommand{\gk}{ginzburg94:_koszul}
\newcommand{\jacobson}{jacobson79:_lie}
\newcommand{\joyal}{joyal81:_une}

\newcommand{\karoubi}{karoubi95:_formes}
\newcommand{\kassel}{kassel95:_quant}
\newcommand{\kont}{kontsevich94:_feynm}
\newcommand{\god}{kontsevich93:_formal}
\newcommand{\lazard}{lazard95:_lois}
\newcommand{\lsv}{loday97:_operad}
\newcommand{\loday}{loday98:_cyclic}
\newcommand{\lodayo}{loday97:_overv_leibn}
\newcommand{\lodayetal}{loday01:_dialg}

\newcommand{\markl}{markl99:_cyclic}
\newcommand{\may}{may72}
\newcommand{\mccleary}{mccleary01}
\newcommand{\mcduff}{mcduff98:_introd}
\newcommand{\mm}{milnor65:_hopf}

\newcommand{\spivak}{spivak79}
\newcommand{\stasheff}{stasheff63:_homot_h}
\newcommand{\sweedler}{sweedler69:_hopf}
\newcommand{\thurston}{thurston:_integ_expres_vassil_knot_invar}
\newcommand{\voronov}{voronov:_notes}
\newcommand{\weibel}{weibel94}

\begin{document}
\title{Symplectic operad geometry and graph homology} 
\author{Swapneel Mahajan}

\address{Department of Mathematics\\
Cornell University\\
Ithaca, NY 14853}
\email{swapneel@math.cornell.edu}
%\date{\today}
%\date{September 20, 1999}

%\thanks{Preliminary version; comments and suggestions are welcome.}

\begin{abstract}
A theorem of Kontsevich relates the homology 
of certain infinite dimensional Lie algebras
to graph homology.
We formulate this theorem using the language 
of reversible operads and mated species.
All ideas are explained using a pictorial calculus
of cuttings and matings.
The Lie algebras are constructed as 
Hamiltonian functions on a symplectic operad manifold.
And graph complexes are defined for any mated species.
The general formulation gives us many examples
including a graph homology for groups.
We also speculate on the role of deformation theory
for operads in this setting.
\end{abstract}

\maketitle

\section{Introduction}  

This paper is my humble tribute to the genius of Maxim Kontsevich.
Needless to say, the credit for any new ideas that occur here
goes to him, and not me.
For how I got involved in this wonderful project,
see the historical note (\ref{ss:hn}).

In the papers~\cite{\god,\kont},
Kontsevich defined three Lie algebras and
related their homology with classical invariants,
including the homology of 
the group of outer automorphisms of a free group and
mapping class groups.
He showed this by reducing the homology computation
in each case to three graph complexes.
His main theorem can be informally stated as

\begin{theorem*}
$H_*(\text{Lie algebra}) = H_*(\text{Graph complex}) = 
H_*(\text{Group}).$
\end{theorem*}

In this paper, we consider only the first part of this theorem.
The definitions of the Lie algebras are motivated by
classical symplectic geometry.
Kontsevich considered three worlds - commutative, associative 
and Lie.
He developed formal analogues of classical symplectic geometry
(which is the commutative case)
for the associative and Lie worlds.
Each of the three Lie algebras was then defined as
Hamiltonian functions on a ``symplectic non-commutative manifold''
with the bracket being the analogue of the usual Poisson bracket.
The symplectic Lie algebra $\fsp$ was a Lie subalgebra
of all three Lie algebras.

In the commutative case, Kontsevich defined a chain complex using graphs.
The homology of this chain complex is known as
graph homology.
He also gave an analogous definition in the associative case
using ribbon graphs
and in the Lie case using more complicated graphs,
which one may call Lie graphs.

The connection between an algebraic or geometric object
like a Lie algebra on one side and 
a combinatorial object like a graph on the other side
is the content of Kontsevich's theorem.
The main tool in proving this connection is to use
classical invariant theory of $\fsp$.

\subsection{The goal of this paper}
From the unified way in which Kontsevich treated the three cases,
one expected them to be a part of a more general theory.
In fact, in \cite{\god}, Kontsevich made the following statement.
\begin{quotation}
Our formalism could be extended to the case of
``Koszul dual pairs of quadratic operads''
(see \cite{\gk}) including Poisson algebras and,
probably, operator algebras etc.
\end{quotation}
This paper arose in an attempt to understand this statement.
In mathematics, just as we have the notion of a 
theorem, proposition, corollary, etc, 
so should we have the notion of a ``Kontsevich sentence'',
the above being a classic example.

Earlier we used the term ``world'' to stand for
one of the three words - commutative, associative 
and Lie.
A possible mathematical substitute for this term is
the word ``operad''.
It is not possible to develop an analogue of symplectic geometry
for every operad.
However, one can do so if the operad is ``reversible''.
The commutative, associative and Lie operads are three examples.
The notion of a reversible operad 
was developed independently by Getzler and Kapranov
under the name of cyclic operad~\cite{\getzlerkap}.
For some information on the history,
see the discussion in (\ref{ss:hn}).

The goal of this paper is to formulate 
the first part of Kontsevich's theorem
in what appears to be its most natural setting -
the world of a reversible operad.
To such an operad $P$,
we associate a Lie algebra $QA_{\infty}$ and a graph complex
$(Q\mathcal G,\partial_E)$
and then show that they have the same homology.
A precise statement is given in Theorem~\ref{t:mt}
in Section~\ref{s:mt}.

\subsection{Organisation of the paper} \label{ss:op}

We now explain the structure of the paper 
and give a guide to how readers 
with different interests and backgrounds can read it.

The two words that dominate the paper are 
species and operads.
These notions are explained in Section~\ref{s:ro}.
In the same section we define reversible operads.
We give plenty of simple examples 
with stress laid on drawing pictures.
More examples are discussed in Sections~\ref{s:me}
and \ref{s:gp}.
We hope that the amount of detail given
along with the examples will be sufficient 
for the reader unfamiliar with these concepts.

In Section~\ref{s:mf}, we construct a 
$$
\text{Mating Functor} \ : \ \{\text{Reversible operads}\}
\longrightarrow \{\text{Species}\}.
$$
Let $P$ be a reversible operad and
$Q$ its image under this functor.
Remember that our first goal is to associate a Lie algebra 
and a graph complex to $P$.
The fact is that the mating functor is the main step
in constructing these two objects,
and the species $Q$ plays a key role.
So we call these two objects $QA_{\infty}$ and
$(Q\mathcal G,\partial_E)$ respectively,
stressing the close connection to the species $Q$.
The letters $A$ and $\mathcal G$ stand for 
algebra and graph respectively.
The subscript $\infty$ refers to the fact that $QA_{\infty}$
will be defined as a direct limit 
of a family of Lie algebras $QA_n$.

Sections~\ref{s:ro} and \ref{s:mf},
on which we have elaborated so far,
form the basic structure
on which the rest of the paper is built.
It is a good idea to avoid technical details
on a first reading and just keep a few concrete examples in mind.
The rest of the paper can be split into three parts.
Part A is about the left hand side of Theorem~\ref{t:mt}
(Lie algebras),
Part B is about the right hand side (graph homology)
and Part C is about stating and proving the theorem.

\subsubsection*{Part A : Sections~\ref{s:sgo}-\ref{s:me}}
This deals with symplectic operad geometry.
We start with an overview in Section~\ref{s:sgo}.
The goal of Section~\ref{s:cm}
is to define the Lie algebra $QA_{\infty}$.
This is done via a simple pictorial calculus 
of cuttings and matings.
The Lie algebra homology $H_*(QA_{\infty})$
has the structure of a Hopf algebra,
which we then explain.
This completes the definition of the left hand side of
Theorem~\ref{t:mt}.

In Section~\ref{s:st},
we explain Kontsevich's symplectic mini-theory.
This is logically not essential to the rest of the paper.
However, it is conceptually the most important section
in Part A.
It shows how the Lie structure on $QA_{\infty}$
emerges naturally from a symplectic form.
This allows us to think of $QA_{\infty}$
as Hamiltonian functions on a ``symplectic operad manifold''.

\subsubsection*{Part B : Sections~\ref{s:gh}-\ref{s:mgh}}

This deals with graph homology.
The reader, who is mainly interested in graph homology,
can directly start with this part even skipping
Sections~\ref{s:ro} and \ref{s:mf}.
We first define a general graph complex $(\mathcal G,\partial_E)$
and then look at other relevant (and smaller) graph complexes,
$(Q\mathcal G,\partial_E)$ being one of them.
This completes the definition of the right hand side of
Theorem~\ref{t:mt}.

Section~\ref{s:gp} is optional.
It contains an interesting example of the theory 
which is based on groups.
In Section~\ref{s:mgh},
we continue with graph homology and develop some tools 
that will be used in the proof of the main theorem.
This involves defining a cochain complex $(\mathcal G,\delta_E)$,
i.e. graph cohomology.
The homology and cohomology are related 
by an interesting and highly non-trivial pairing on graphs
(\ref{ss:pg}).

\subsubsection*{Part C : Sections~\ref{s:mt}-\ref{s:con}}

This deals with the connection between symplectic geometry
and graph homology.
Section~\ref{s:mt} contains a precise formulation
of the main theorem.
To digest the statement completely,
the material on Lie algebras and graph homology
in Sections~\ref{s:cm} and \ref{s:gh} respectively
is a prerequisite.
The next three sections deal with the proof of the theorem.

The main ideas of the proof are already present 
in the commutative case.
Since the proof is quite involved,
we suggest that the reader specialise to this case
on a first reading.
The Lie algebra $QA_{\infty}$
in this case is easy to define directly.
So Sections~\ref{s:ro}-\ref{s:me} are not necessary.
Thus the reader, who mainly wants to see the ideas in the proof,
may start with Part C directly,
skipping Part A and referring back to Part B as necessary.

The proof of Theorem~\ref{t:mt} is given in three parts
(Sections~\ref{s:fin},\ref{s:sta} and~\ref{s:con}).
In each of these sections, we state and prove
a theorem of the form
\[
H_*(\text{Lie algebra}) = H_*(\text{Graph complex}),
\]
gradually getting closer to our goal.
As mentioned earlier,
the main tool in the proof
is classical invariant theory of $\fsp$.
This part is done in Section~\ref{s:fin}.
The second part of the theorem given in
Section~\ref{s:sta} deals with the stability issue
$(n \to \infty)$
and makes strong use of the ideas of 
Section~\ref{s:mgh}.
In Section~\ref{s:con}, we understand the space
of primitive elements of the Hopf algebra $H_*(QA_{\infty})$.
This is largely a matter of unwinding definitions
and there are no new ideas here.
Since Hopf algebras play a minor role in this paper,
the reader unfamiliar with them may simply omit this section
without losing any of the main ideas.

\subsubsection*{Appendix}
We have included two appendices,
which explain how deformation theory
relates to the ideas of this paper.
They are not logically essential
to the understanding of the main theorem.
Appendix~\ref{ss:dq} speculates 
on the deformation quantisation problem for operads
in this setting.
In Appendix~\ref{ss:dg},
we show that a certain Lie bracket 
on graph homology defined in~\cite{\cv} is zero.

\subsection{A historical note} \label{ss:hn}

This work grew out of a seminar
organised by Karen Vogtmann,
devoted to understanding Kontsevich's work
(Fall 2000).
Other participants included 
D. Ciubotaru, F. Gerlits, D. Brown, J. Conant, 
M. Horak, F. Schwartz, M. Cohen and J. West.
My understanding of Kontsevich's ideas,
particularly the proof of the commutative case,
is due to them.
The material in (\ref {ss:og}-\ref {ss:gc})
and (\ref {ss:lah}-\ref {ss:inv}) is based 
on the seminar notes.
Conant and Vogtmann are writing an exposition 
that has some overlap with this.

The notion of a reversible operad and
its role in symplectic operad geometry
was done in July 2001 (Sections~\ref{s:ro}-\ref{s:me}).
The rest of the paper, namely, 
defining graph homology for a mated species and Theorem~\ref{t:mt},
was then relatively easy. 
The result in Appendix~\ref{ss:dg} was proved later
in November 2001.
These ideas, in some form or the other,
had already appeared in the works of 
Getzler and Kapranov~\cite{\getzlerkap,\getzkap},
Ginzburg~\cite{\ginz} and Markl~\cite{\markl}. 
Being ignorant of the operad literature at that time,
this was a gradual discovery for me.
In this paper, I give my own viewpoint of the subject,
which was formed by reading Kontsevich.

In \cite{\getzlerkap},
Getzler and Kapranov introduced the notion of cyclic operads
in order to extend the formalism of cyclic homology
for associative algebras~\cite{\loday} to operad algebras.
This notion coincides with what we call a reversible operad.
For such an operad $P$, they introduced a functor
$$
\lambda(P,_{\ -}) \ : \ \{P\text{-algebras}\}
\longrightarrow \{\text{Vector spaces}\}.
$$
In our notation (see Sections~\ref{s:sgo} and \ref{s:cm}), 
$\lambda(P,PA) = QA$,
and equation~\eqref{e:r} in (\ref {ss:fa})
can be taken as the definition 
of $\lambda$.
The mating functor is similar in spirit to $\lambda$, except that
we deal with operads and species rather than
algebras and vector spaces.

In \cite{\ginz}, Ginzburg explained Kontsevich's
symplectic mini-theory for a cyclic Koszul operad.
He denoted the functor $\lambda$ by the letter $R$.
Since a Koszul operad is quadratic,
the relations in equation~\eqref{e:r} take the simpler form
\begin{equation*}
a \otimes b = b \otimes a  \quad \text{and} \quad
\mu(a,b) \otimes c = a \otimes \mu(b,c) \quad \text{for} \quad \mu \in P[2].
\end{equation*}
These are precisely the relations that Kontsevich 
wrote in~\cite{\god}.
Elements of associative and Lie calculus 
had appeared earlier in the works of Karoubi~\cite{\karoubi}
and Drinfeld~\cite{\drinfeld}.

However, as we show, to do symplectic geometry,
it is sufficient to assume that $P$ is reversible,
not necessarily Koszul or quadratic.
In fact, all the examples in Section~\ref{s:me}
and some earlier ones too are non-quadratic.
We discovered reversibility by requiring that
the notion of a partial derivative make sense
(Proposition~\ref{p:pd}).
So in this sense, it is also a necessary condition
to do symplectic geometry.

In~\cite{\markl}, Markl considered graph complexes
associated to cyclic operads.
And for this, he referred back to the Feynman transform
construction of Getzler and Kapranov~\cite{\getzkap}.

We conclude the introduction with two remarks. 
Koszul operads will play an important role in the second part
of Kontsevich's theorem,
which we do not consider here.
If we switch from Lie to Liebniz algebras~\cite{\lodayo}
then we will end up with another avatar of graph homology.
Now vertices of graphs will have honest labels, 
unlike in usual graph homology,
where vertices are labelled only upto orientation.

\newpage
\section{Species, operads and reversible operads}  \label{s:ro}

In this section, we give a brief introduction 
to species and operads and then define reversible operads.
An excellent exposition on species can be found in the book
by Bergeron, Labelle and Leroux~\cite{\bll}.
The notes by Voronov~\cite{\voronov}
contain a good review on operads.
A more comprehensive reference is the seminal work
of Ginzburg and Kapranov~\cite{\gk}.
For some recent developments on operads, 
see~\cite{\getzler,\getzlerkap,\getzkap,\lsv,\lodayetal}.
For earlier literature,
see~\cite{\joyal,\may,\stasheff,\lazard},
where the concepts of species and operads first originated.

\subsection{Species}  \label{ss:s}

A species $Q$ is a functor from the category of 
(Finite sets, bijections) to the category of 
(Sets, maps).
We denote the image of a set $I$ by $Q[I]$
and say that $Q[I]$ is the set of $Q$-structures on the set $I$.
The set $Q[\phi]$ will always be empty.

\medskip
\begin{center}
\begin{tabular}{|c|c|} \hline 
Species $Q$ & Picture of an element of $Q[I]$ 
for $I=\{a,b,c,d\}$ \\ \hline \hline
$c$ & 
\begin{minipage}{.4 in}
\vspace{.1in}
\begin{center}
\psfrag{a}{\Huge $a$}
\psfrag{b}{\Huge $b$}
\psfrag{c}{\Huge $c$}
\psfrag{d}{\Huge $d$}
\resizebox{.4 in}{.4 in}{\includegraphics{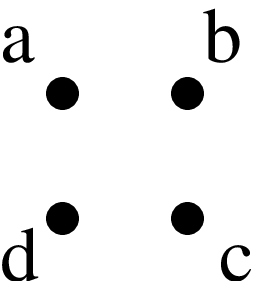}}
\end{center}
\vspace{.1in}
\end{minipage} 
\ \ = \ \
\begin{minipage}{.4 in}
\vspace{.1in}
\begin{center}
\psfrag{a}{\Huge $a$}
\psfrag{b}{\Huge $b$}
\psfrag{c}{\Huge $c$}
\psfrag{d}{\Huge $d$}
\resizebox{.4 in}{.4 in}{\includegraphics{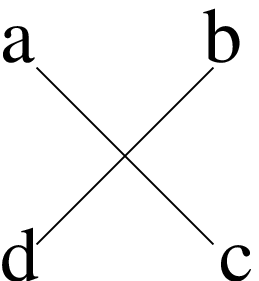}}
\end{center}
\vspace{.1in}
\end{minipage} 
\ \ = \ \
\begin{minipage}{.4 in}
\vspace{.1in}
\begin{center}
\psfrag{a}{\Huge $a$}
\psfrag{b}{\Huge $b$}
\psfrag{c}{\Huge $c$}
\psfrag{d}{\Huge $d$}
\resizebox{.4 in}{.4 in}{\includegraphics{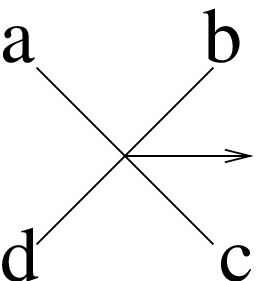}}
\end{center}
\vspace{.1in}
\end{minipage} \\ \hline
$a$ &  
\begin{minipage}{.65 in}
\vspace{.1in}
\begin{center}
\psfrag{a}{\Huge $a$}
\psfrag{b}{\Huge $b$}
\psfrag{c}{\Huge $c$}
\psfrag{d}{\Huge $d$}
\resizebox{.65 in}{.2 in}{\includegraphics{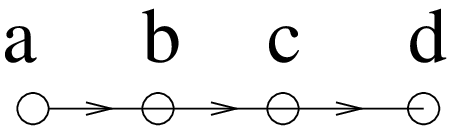}}
\end{center}
\vspace{.1in}
\end{minipage}
\ \ = \ \
\begin{minipage}{.55 in}
\vspace{.1in}
\begin{center}
\psfrag{a}{\Huge $a$}
\psfrag{b}{\Huge $b$}
\psfrag{c}{\Huge $c$}
\psfrag{d}{\Huge $d$}
\resizebox{.55 in}{.5 in}{\includegraphics{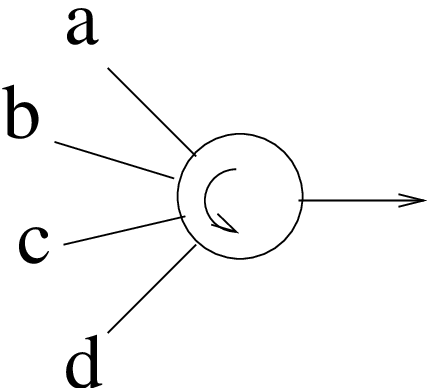}}
\end{center}
\vspace{.1in}
\end{minipage} \\ \hline
$aa$ &  
\begin{minipage}{.5 in}
\vspace{.1in}
\begin{center}
\psfrag{a}{\Huge $a$}
\psfrag{b}{\Huge $b$}
\psfrag{c}{\Huge $c$}
\psfrag{d}{\Huge $d$}
\resizebox{.5 in}{.5 in}{\includegraphics{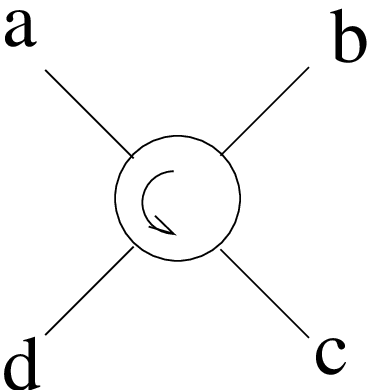}}
\end{center}
\vspace{.1in}
\end{minipage} \\ \hline
$t$ &
\begin{minipage}{.8 in}
\vspace{.1in}
\begin{center}
\psfrag{a}{\Huge $a$}
\psfrag{b}{\Huge $b$}
\psfrag{c}{\Huge $c$}
\psfrag{d}{\Huge $d$}
\resizebox{.75 in}{.4 in}{\includegraphics{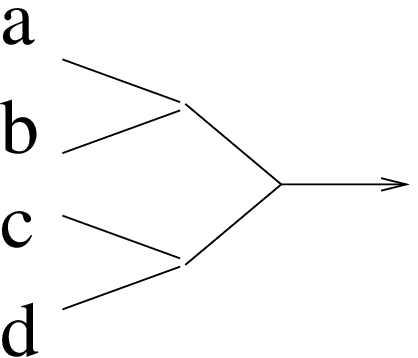}}
\end{center}
\vspace{.1in}
\end{minipage} \\ \hline
$tt$ & 
\begin{minipage}{.8 in}
\vspace{.1in}
\begin{center}
\psfrag{a}{\Huge $a$}
\psfrag{b}{\Huge $b$}
\psfrag{c}{\Huge $c$}
\psfrag{d}{\Huge $d$}
\resizebox{.7 in}{.45 in}{\includegraphics{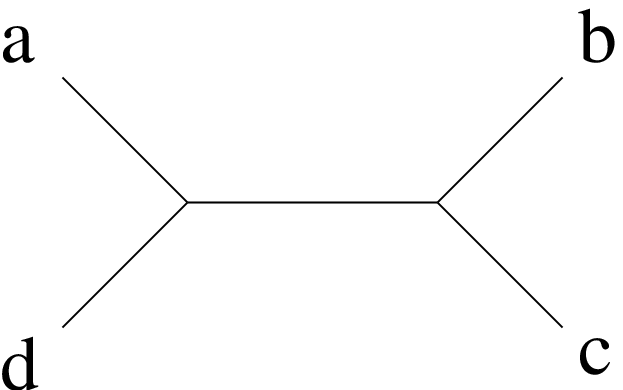}}
\end{center}
\vspace{.1in}
\end{minipage} \\ \hline
\end{tabular}
\end{center}

\vspace{0.3 cm}
\noindent
We now give some examples of species.
The pictures that go along with the examples are shown in the table.
\begin{itemize}

\item
$c[I]=\{I\}$,
that is, there is exactly one $c$-structure on the set $I$.
There are a variety of ways to show this via pictures.

\item
$a[I]=$ the set of linear orders on the set $I$.
We have shown two pictures for it.

\item
$aa[I]=$ the set of cyclic orders on the set $I$
if $|I| \geq 2$ and empty otherwise.

\item
$t[I]=$ the set of rooted trees with leaves labelled 
by elements of the set $I$.

\item
$tt[I]=$ the set of trees with leaves labelled 
by elements of the set $I$
if $|I| \geq 2$ and empty otherwise.

\end{itemize}
An element of $Q[I]$ can be schematically drawn as

\[
\begin{minipage}{.8 in}
\begin{center}
\psfrag{a}{\Huge $a$}
\psfrag{b}{\Huge $b$}
\psfrag{c}{\Huge $c$}
\psfrag{d}{\Huge $d$}
\resizebox{.6 in}{.6 in}{\includegraphics{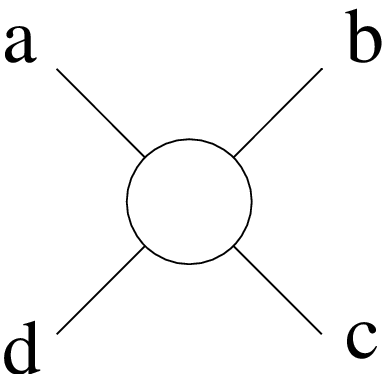}}
\end{center}
%\label{f:q}
\end{minipage}
\in Q[I] \ \ \text{for} \ \ I=\{a,b,c,d\}.
\]
Note that the picture for the species ``$aa$''
fits in with this representation perfectly.
And the remaining examples can also be
made to fit in without difficulty.
For example,
the species ``$c$'' can be drawn as
\[
\begin{minipage}{.8 in}
\begin{center}
\psfrag{a}{\Huge $a$}
\psfrag{b}{\Huge $b$}
\psfrag{c}{\Huge $c$}
\psfrag{d}{\Huge $d$}
\resizebox{.5 in}{.5 in}{\includegraphics{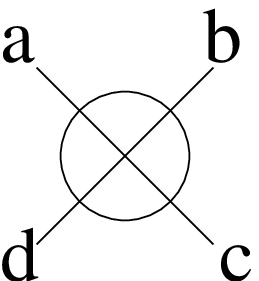}}
\end{center}
\end{minipage},
\]
and so on for the others too.

A species can be equivalently defined as a
sequence $Q[0],Q[1],Q[2], \ldots$,
where $Q[n]$ is a $\Sigma_n$ set,
where $\Sigma_n$ is the symmetric group
on $n$ letters.
Here we have abbreviated our notation so that
$Q[n]$ stands for $Q[\{1,2,\ldots,n\}]$.
In the pictures, 
the group $\Sigma_n$ acts by permuting the $n$ letters.

It is also useful to use a more general definition of species 
where the target category is replaced
by the category of vector spaces.
In this case,
a species is a sequence as above with
$Q[n]$ being a linear representation of $\Sigma_n$.
It is clear that using vector spaces as the target category
is more general because
one can go from sets to vector spaces 
by linearising the representation.

\begin{remark}
Throughout this paper,
we use only sets and vector spaces as the target categories.
However, this restriction is mainly for simplicity.
For example,
one may consider a species in the category of topological spaces.
It would be a sequence $X[0],X[1],X[2], \ldots$
of topological spaces with an action of
$\Sigma_n$ on $X[n]$.
\end{remark}

\subsection{Operads}  \label{ss:o}

Note that for each of the examples $c,a,t$,
we suggested a picture with an arrow ``$\rightarrow$''
drawn in it.
This is possible because these species
have the additional structure of an operad.

We begin with an informal discussion on operads.
An operad $P$ is a species in which there is a substitution rule.
It can be shown schematically as
\[
\begin{minipage}{.8 in}
\begin{center}
\psfrag{a}{\Huge $a$}
\psfrag{b}{\Huge $b$}
\psfrag{c}{\Huge $c$}
\psfrag{d}{\Huge $d$}
\resizebox{.6 in}{.6 in}{\includegraphics{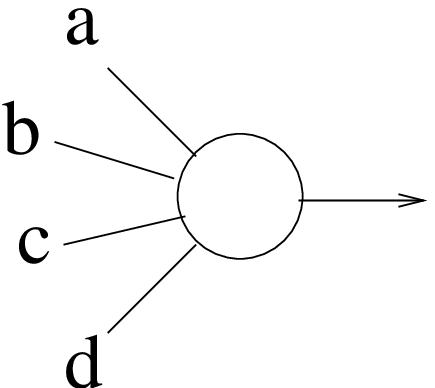}}
\end{center}
%\label{f:p}
\end{minipage}
\in P[I] \ \ \text{for} \ \ I=\{a,b,c,d\}.
\]
One thinks of $a,b,c,d$ as four inputs and the arrow as an output.
The substitution rule allows us to feed 
the output of one object $p_2$ into 
the input of another object $p_1$.
We write this as $p_2 \rightarrow p_1$ or
$p_1 \leftarrow p_2$.
This can be shown as
\[
\begin{minipage}{.6 in}
\begin{center}
\psfrag{a}{\Huge $a$}
\psfrag{b}{\Huge $b$}
\psfrag{c}{\Huge $c$}
\psfrag{d}{\Huge $d$}
\psfrag{p2}{\Huge $p_2$}
\resizebox{.6 in}{.6 in}{\includegraphics{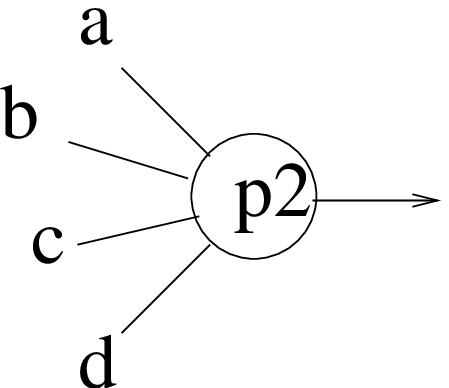}}
\end{center}
%\label{f:psubt}
\end{minipage}
\hfill
\begin{minipage}{.8 in}
\begin{center}
\psfrag{x}{\Huge $x$}
\psfrag{y}{\Huge $y$}
\psfrag{z}{\Huge $z$}
\psfrag{p1}{\Huge $p_1$}
\resizebox{.6 in}{.6 in}{\includegraphics{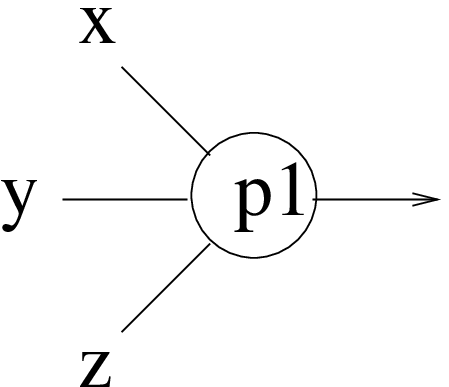}}
\end{center}
%\label{f:psubi}
\end{minipage}
=
\begin{minipage}{.8 in}
\begin{center}
\psfrag{a}{\Huge $a$}
\psfrag{b}{\Huge $b$}
\psfrag{c}{\Huge $c$}
\psfrag{d}{\Huge $d$}
\psfrag{x}{\Huge $x$}
\psfrag{z}{\Huge $z$}
\resizebox{.6 in}{.7 in}{\includegraphics{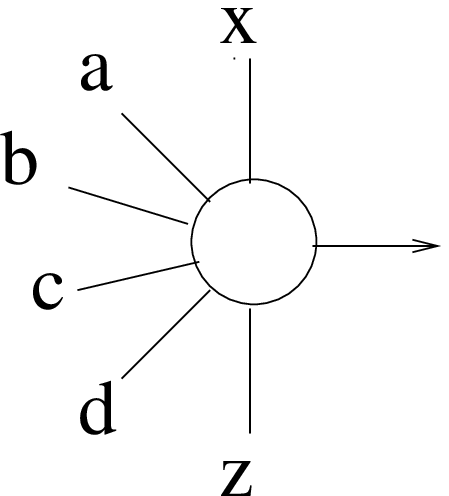}}
\end{center}
%\label{f:psubo}
\end{minipage}.
\]
If we wish to be more specific, we will write
$p_2 \yright p_1$.
This means that we feed $p_2$ to that input of $p_1$
whose label is $y$.

Now we show how the substitution rule
works in each of the examples $c,a,t$.
\[
\begin{minipage}{.6 in}
\begin{center}
\psfrag{a}{\Huge $a$}
\psfrag{b}{\Huge $b$}
\psfrag{c}{\Huge $c$}
\psfrag{d}{\Huge $d$}
\resizebox{.5 in}{.5 in}{\includegraphics{csubt.eps}}
\end{center}
%\label{f:csubt}
\end{minipage}
\hfill
\begin{minipage}{.8 in}
\begin{center}
\psfrag{x}{\Huge $x$}
\psfrag{y}{\Huge $y$}
\psfrag{z}{\Huge $z$}
\resizebox{.6 in}{.6 in}{\includegraphics{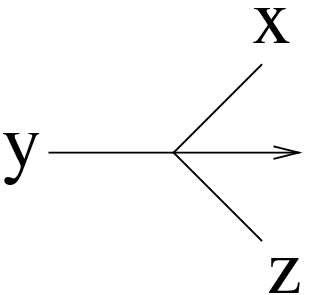}}
\end{center}
%\label{f:csubi}
\end{minipage}
=
\begin{minipage}{.8 in}
\begin{center}
\psfrag{a}{\Huge $a$}
\psfrag{b}{\Huge $b$}
\psfrag{c}{\Huge $c$}
\psfrag{d}{\Huge $d$}
\psfrag{x}{\Huge $x$}
\psfrag{z}{\Huge $z$}
\resizebox{.6 in}{.6 in}{\includegraphics{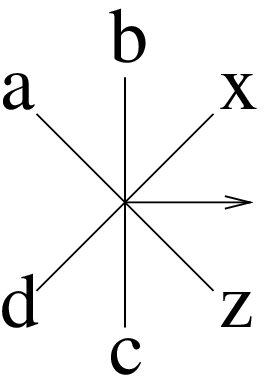}}
\end{center}
%\label{f:csubo}
\end{minipage}.
\]

\[
\begin{minipage}{.6 in}
\begin{center}
\psfrag{a}{\Huge $a$}
\psfrag{b}{\Huge $b$}
\psfrag{c}{\Huge $c$}
\psfrag{d}{\Huge $d$}
\resizebox{.6 in}{.6 in}{\includegraphics{asubt.eps}}
\end{center}
%\label{f:asubt}
\end{minipage}
\hfill
\begin{minipage}{.8 in}
\begin{center}
\psfrag{x}{\Huge $x$}
\psfrag{y}{\Huge $y$}
\psfrag{z}{\Huge $z$}
\resizebox{.6 in}{.6 in}{\includegraphics{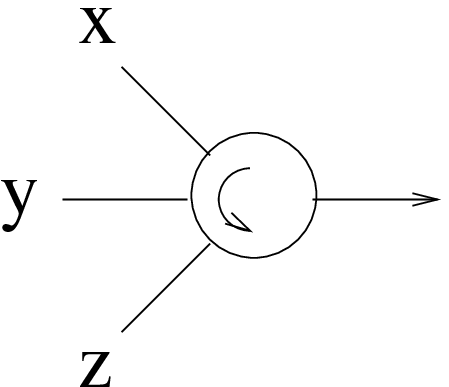}}
\end{center}
%\label{f:asubi}
\end{minipage}
=
\begin{minipage}{.8 in}
\begin{center}
\psfrag{a}{\Huge $a$}
\psfrag{b}{\Huge $b$}
\psfrag{c}{\Huge $c$}
\psfrag{d}{\Huge $d$}
\psfrag{x}{\Huge $x$}
\psfrag{z}{\Huge $z$}
\resizebox{.6 in}{.7 in}{\includegraphics{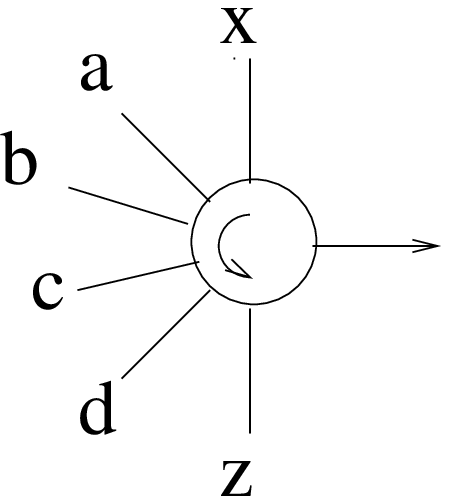}}
\end{center}
%\label{f:asubo}
\end{minipage}.
\]
For the species $t$, we graft the root of the first tree 
into the specified leaf of the second tree.
\[
\begin{minipage}{.7 in}
\begin{center}
\psfrag{a}{\Huge $a$}
\psfrag{b}{\Huge $b$}
\psfrag{c}{\Huge $c$}
\psfrag{d}{\Huge $d$}
\resizebox{.7 in}{.6 in}{\includegraphics{tsubt.eps}}
\end{center}
%\label{f:tsubt}
\end{minipage}
\hfill
\begin{minipage}{.8 in}
\begin{center}
\psfrag{x}{\Huge $x$}
\psfrag{y}{\Huge $y$}
\psfrag{z}{\Huge $z$}
\resizebox{.6 in}{.5 in}{\includegraphics{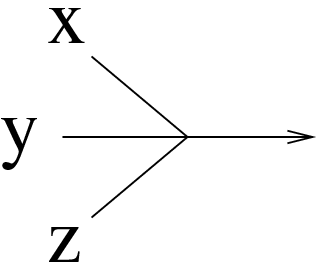}}
\end{center}
%\label{f:tsubi}
\end{minipage}
=
\begin{minipage}{1 in}
\begin{center}
\psfrag{a}{\Huge $a$}
\psfrag{b}{\Huge $b$}
\psfrag{c}{\Huge $c$}
\psfrag{d}{\Huge $d$}
\psfrag{x}{\Huge $x$}
\psfrag{z}{\Huge $z$}
\resizebox{.9 in}{.6 in}{\includegraphics{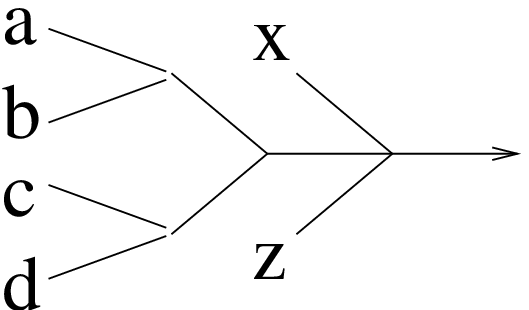}}
\end{center}
%\label{f:tsubo}
\end{minipage}.
\]
Note that there is another way possible.
One could contract the edge at which the grafting took place.
This is what happened for the species ``$c$''.
There are numerous other variations one can do on this example;
for instance, one could use rooted binary trees.

To give a formal definition,
an operad $P$ is a species with a substitution rule
which is
\begin{itemize}
\item
associative, and

\item
compatible with the morphisms in the source category.

\end{itemize}
The first condition says that 
if we perform two substitutions,
one after the other,
then the order in which we do them does not matter.
The second condition says that substitution
commutes with relabelling of the inputs.
For more detail on these conditions,
see the references cited earlier.

For convenience of bookkeeping, it is useful to label
the output as well.
This idea can be formalised as follows.
For any operad $P$, define a bi-functor
\[
P : (\text{Sets,\ bijections}) \times
(\text{Singleton sets,\ bijections}) \to
(\text{Vector spaces,\ maps})
\]
by setting $P[I,U]=P[I]$,
where $U$ is a singleton set.
One thinks of an element of $P[I,U]$
as an element of $P[I]$
whose output is labelled by the element of $U$.
For the substitution rule in this situation, 
we demand that the output label of $p_2$
match the input label of $p_1$.
So the notation $p_2 \yright p_1$
now means that 
the label of the output of $p_2$ is $y$
and it is fed to an input of $p_1$
whose label is again $y$.

Just as for species,
we can define an operad as a sequence $P[n]$
of $\Sigma_n$ modules
equipped with a substitution rule
that is associative and $\Sigma_n$ invariant. 
A unit element of $P$ is an element of $P[1]$ 
which when substituted into an input of any
$p \in P[n]$ gives back $p$.
We assume that our operads have units
and call the unit element $u$.
In the second notation,
we would write $P[n,1]$ instead of $P[n]$
and the unit element would lie in $P[1,1]$.

\subsection{Examples} \label{ss:eg}

We now organise together a basic set of examples.
These include the ones shown in the table (\ref{ss:s}).
The motivation for organising them in this manner will become clear
when we discuss the mating functor (\ref{ss:mf}).
We make two quick remarks.
All examples here are based on sets.
Furthermore, the ones that are operads have unit elements
and satisfy $P[1]=\Q$,
i.e. $P[1]$ is just the span of the unit element $u$ of $P$.

\subsubsection{The unit operad $u$ and species $uu$} \label{sss:u}
Define 
$u[I]=\{I\}$ if $I$ is a singleton and empty otherwise.
And similarly, put
$uu[I]=\{I\}$ if $|I|=2$ and empty otherwise.
These two examples play a fundamental role in this paper
because of their connection to the symplectic Lie algebra
$\fsp$, see (\ref{ss:sl}).
\[
\begin{minipage}{1 in}
\begin{center}
\psfrag{a}{\huge $a$}
\psfrag{b}{\huge $b$}
\psfrag{u}{\Huge $u$}
\resizebox{.7 in}{.25 in}{\includegraphics{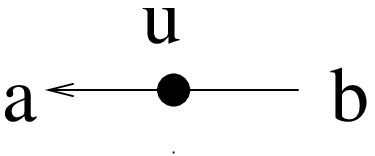}}
\end{center}
\end{minipage}
\hfill
\begin{minipage}{2 in}
\begin{center}
\psfrag{1}{}
\psfrag{2}{}
\psfrag{p1}{\huge $a$}
\psfrag{q2}{\huge $b$}
\psfrag{uu}{\huge $uu$}
\resizebox{1 in}{.2 in}{\includegraphics{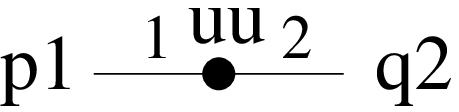}}
\end{center}
\end{minipage}
\]

\subsubsection{The commutative operad $c$ and species $cc$}
\label{sss:c}

Define 
$c[I]=\{I\}$,
that is, there is exactly one $c$-structure on the set $I$.
Also define $cc[I]=\{I\}$ if $|I| \geq 2$ and empty otherwise.
Though this is also an operad (without a unit),
we call it the commutative species.
As far as pictures go,
we will use the one with arrow ``$\rightarrow$'' for ``$c$''
and the one without it for ``$cc$''.

\subsubsection{The associative operad $a$ and species $aa$}  \label{sss:a}
Let $a[I]=$ the set of linear orders on the set $I$.
Also let $aa[I]=$ the set of cyclic orders on the set $I$
if $|I| \geq 2$ and empty otherwise.

\subsubsection{The tree operad $t$ and species $tt$}  \label{sss:t}

Let $t[I]=$ the set of rooted trees with leaves labelled 
by elements of the set $I$.
For the tree species $tt$, we drop the word ``rooted''. 
In other words,
$tt[I]=$ the set of trees with leaves labelled 
by elements of the set $I$
if $|I| \geq 2$ and empty otherwise.
Since there is no root, one cannot define a substitution rule.

\subsubsection{The chord operad $k$ and species $kk$}  \label{sss:k}
The chord species $kk[I]$
is the set of chord diagrams on the set $I$.
That is, an element of $kk[I]$
specifies a way to pair off the elements of $I$.
Clearly, this is non-empty only when the cardinality of $I$ is even.
Define the chord operad $k[I]$
in the same way, except that
one of the chords is left hanging at one end.
Hence this is non-empty only when the cardinality of $I$ is odd.
The substitution rule should be clear.

\medskip
\noindent
The species $aa,tt,kk$ have no natural substitution rule;
so they are not operads.
We will give some more examples later (Sections~\ref{s:me}
and \ref{s:gp}).

\subsection{Reversible operads} \label{ss:ro}

An operad $P$ is reversible if there is a rule that allows
us to switch the output with any given input.
One represents this pictorially as
\[
r_{x,y}
\left(
\begin{minipage}{.8 in}
\begin{center}
\psfrag{a}{\Huge $a$}
\psfrag{b}{\Huge $b$}
\psfrag{x}{\Huge $x$}
\psfrag{y}{\Huge $y$}
\resizebox{.7 in}{.6 in}{\includegraphics{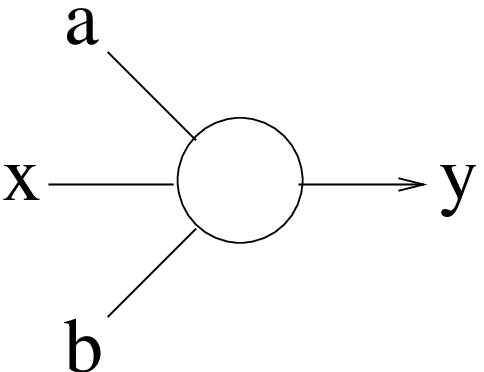}}
\end{center}
%\label{f:prevt}
\end{minipage}
\right)
=
\begin{minipage}{.8 in}
\begin{center}
\psfrag{a}{\Huge $a$}
\psfrag{b}{\Huge $b$}
\psfrag{x}{\Huge $x$}
\psfrag{y}{\Huge $y$}
\resizebox{.7 in}{.6 in}{\includegraphics{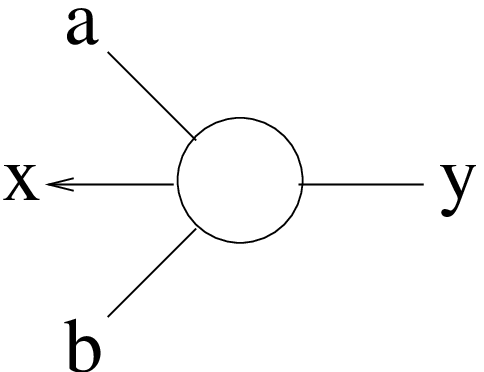}}
\end{center}
%\label{f:prevo}
\end{minipage}
.
\]
To give a formal definition, for each set $I$ and $x \in I$, 
there exist reversal maps
\[
r_{x,y} : P[I,\{y\}] \to P[I \setminus \{x\} \cup \{y\},\{x\}],
\]
which are subject to two conditions.

1. We require
$r_{y,x} \circ r_{x,y} = \ \text{identity}$ and
$r_{x,y} \circ r_{y,z} = r_{x,z}.$ 
It follows from these two relations that the composite
of a sequence of reversals is either the identity map
or a single reversal. 

2. Reversal must be compatible with substitution.
That is,
$$
  r_{x,z} (p_2 \yright p_1) \ = \left\{ \begin{array}{c l}
    p_2 \yright r_{x,z} (p_1) &
\mbox{ if $x$ is an input of $p_1$,} \\
    r_{x,y} (p_2) \yleft r_{y,z} (p_1) &
\mbox{ if $x$ is an input of $p_2$.} 
                    \end{array} \right. $$
Here $y$ labels one of the inputs of $p_1$
as well as the output of $p_2$ 
and $z$ labels the output of $p_1$.

The examples of operads $(P=u,c,a,t,k)$
that we gave in (\ref {ss:eg}) are all reversible.
For each of these examples,
there is a picture of the operad
that fits in with the schematic one.
The reversal map in all cases then works exactly
as shown above in the schematic picture.

We give some intuitive meaning to the conditions
on the reversal maps.
Condition 1 says that reversal is an external operation.
The internal data of the operad remains unchanged.
Condition 2 says that operad substitution
commutes with operad reversal.
This will become clearer from the examples
of non-reversible operads given in (\ref{ss:ne}).

\begin{remark}
The reversal rules may not be unique.
For example, in (\ref{ss:g}), we give two distinct ways
to reverse the graph operad.
\end{remark}

\subsection{Non-examples}  \label{ss:ne}
We now give two examples of operads that are not reversible.
So they play no role in the rest of the paper.
Let Perm$[I] = I$ and
Dias$[I] =$ the set of linear orders with 
a distinguished element on the set $I$.
\[
\begin{minipage}{.8 in}
\begin{center}
\psfrag{a}{\Huge $a$}
\psfrag{b}{\Huge $b$}
\psfrag{c}{\Huge $c$}
\psfrag{d}{\Huge $d$}
\resizebox{.5 in}{.5 in}{\includegraphics{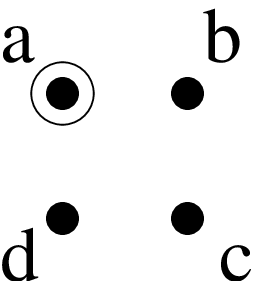}}
\end{center}
\end{minipage}
\in \text{Perm}[I], \qquad \qquad \qquad
\begin{minipage}{.8 in}
\begin{center}
\psfrag{a}{\Huge $a$}
\psfrag{b}{\Huge $b$}
\psfrag{c}{\Huge $c$}
\psfrag{d}{\Huge $d$}
\resizebox{.7 in}{.6 in}{\includegraphics{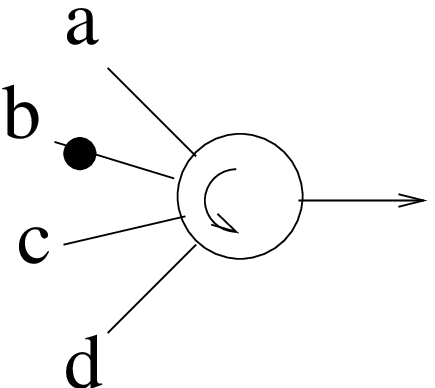}}
\end{center}
\end{minipage}
\in \text{Dias}[I].
\]
The pictures show examples for $I=\{a,b,c,d\}$.

We explain the substitution rule for the operad Dias.
There are two cases depending on whether
the input into which the substitution occurs
is circled or not.
\[
\begin{minipage}{.6 in}
\begin{center}
\psfrag{a}{\Huge $a$}
\psfrag{b}{\Huge $b$}
\psfrag{c}{\Huge $c$}
\psfrag{d}{\Huge $d$}
\resizebox{.6 in}{.6 in}{\includegraphics{diast.eps}}
\end{center}
%\label{f:diast}
\end{minipage}
\hfill
\begin{minipage}{.8 in}
\begin{center}
\psfrag{x}{\Huge $x$}
\psfrag{y}{\Huge $y$}
\psfrag{z}{\Huge $z$}
\resizebox{.6 in}{.6 in}{\includegraphics{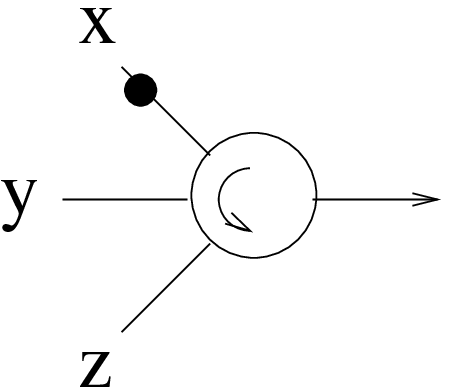}}
\end{center}
%\label{f:diasi}
\end{minipage}
=
\begin{minipage}{.8 in}
\begin{center}
\psfrag{a}{\Huge $a$}
\psfrag{b}{\Huge $b$}
\psfrag{c}{\Huge $c$}
\psfrag{d}{\Huge $d$}
\psfrag{x}{\Huge $x$}
\psfrag{z}{\Huge $z$}
\resizebox{.6 in}{.7 in}{\includegraphics{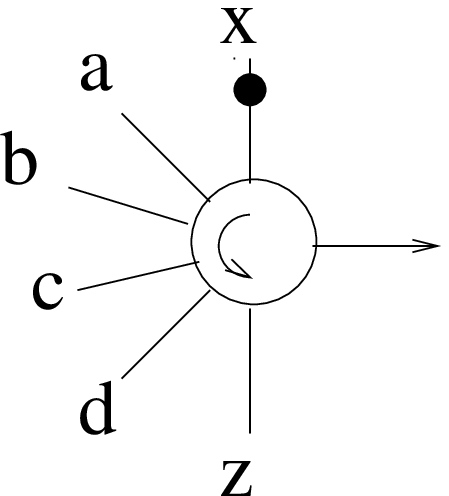}}
\end{center}
%\label{f:diaso}
\end{minipage}.
\]
The picture illustrates the second case.
The first case works as follows.
If the input $y$ were circled instead of $x$
then the result would be a circle on the input $b$.
The substitution rule for Perm is similar.

Observe that Perm and Dias are closely related
to the commutative (\ref{sss:c}) and associative (\ref{sss:a})
operads respectively.
They are binary quadratic and were introduced
by Loday and his collaborators 
with motivation from algebraic K-theory~\cite{\lodayetal}.
For other examples of non-reversible operads,
see~\cite{\getzlerkap}.

\section{The mating functor}  \label{s:mf}

In this section, we introduce the mating functor,
which allows us to construct a mated species 
from a reversible operad.
We then explain the concept of a partial derivative
for a mated species (Proposition~\ref{p:pd}).
This is the most crucial part of the theory and 
shows exactly how the conditions 
on the reversal map given in (\ref{ss:ro}) arise.

\subsection{The mating functor}  \label{ss:mf}

Define a species $Q$ starting with a reversible operad $P$
as follows. 
For any finite set $K$, let
$$Q[K] = \bigoplus_{I \sqcup J = K} P[I,U] \otimes P[J,U],$$
subject to the two relations:

(R1) $p_2 \otimes p_1 = p_1 \otimes p_2,$ and 

(R2) $p_3 \rightarrow p_2 \otimes p_1 = p_3 \otimes r(p_2) \leftarrow p_1.$

\medskip
\noindent
The set $U$ is any singleton,
all choices being considered equivalent.
We interpret the tensor sign as a mating, and say that
$p_2 \otimes p_1$ is a mating of $p_2$ and $p_1$.
We show it as 
\begin{minipage}{1.2 in}
\begin{center}
\psfrag{p1}{\Huge $p_1$}
\psfrag{p2}{\Huge $p_2$}
\resizebox{1 in}{.3 in}{\includegraphics{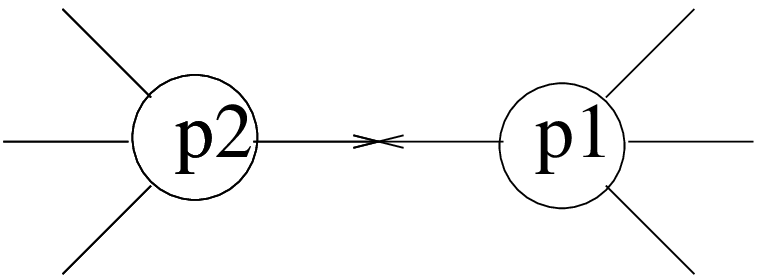}}
\end{center}
\end{minipage}.
For simplicity of notation,
we are suppressing the labels.
The symmetry of relation $(R1)$ is thus built into the picture.
If either of $p_2$ or $p_1$ is the unit element $u$
then we say that the mating 
$p_2 \otimes p_1$ is trivial.
We show it as 
\begin{minipage}{1 in}
\begin{center}
\psfrag{u}{\Huge $u$}
\psfrag{p2}{\Huge $p_2$}
\resizebox{.8 in}{.3 in}{\includegraphics{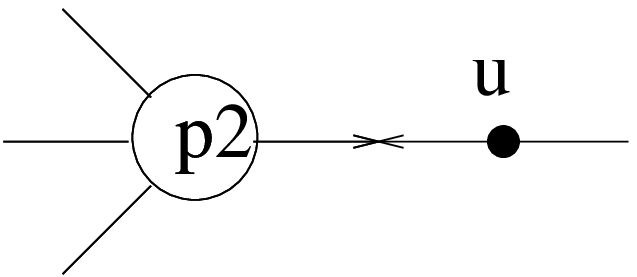}}
\end{center}
\end{minipage}.
In a nontrivial mating,
we will call the segment joining $p_2$ and $p_1$ an \emph{ideal edge}.
It has two opposing arrowheads in the centre.

\begin{remark}
Ideal edges will play an important role later
in the definition of graph cohomology (\ref{ss:stacb}).
\end{remark}

The second relation $(R2)$ is written as under
\[
\begin{minipage}{1.9 in}
\begin{center}
\psfrag{p1}{\Huge $p_3$}
\psfrag{p2}{\Huge $p_2$}
\psfrag{p3}{\Huge $p_1$}
\resizebox{1.7 in}{.4 in}{\includegraphics{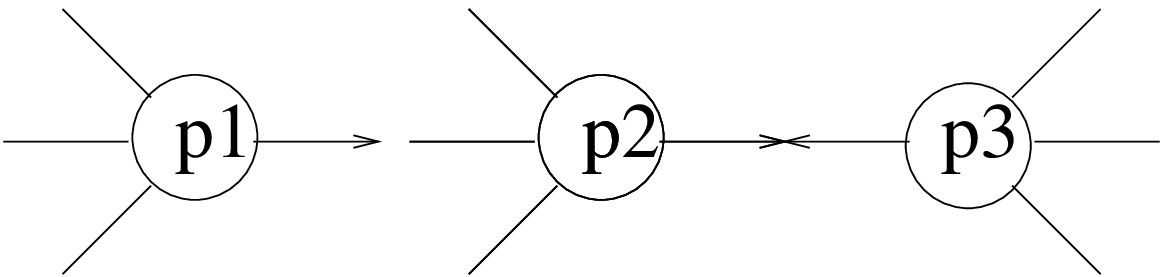}}
\end{center}
\end{minipage}
=
\begin{minipage}{1.9 in}
\begin{center}
\psfrag{p1}{\Huge $p_3$}
\psfrag{pr}{\Huge $\overline{p_2}$}
\psfrag{p3}{\Huge $p_1$}
\resizebox{1.7 in}{.4 in}{\includegraphics{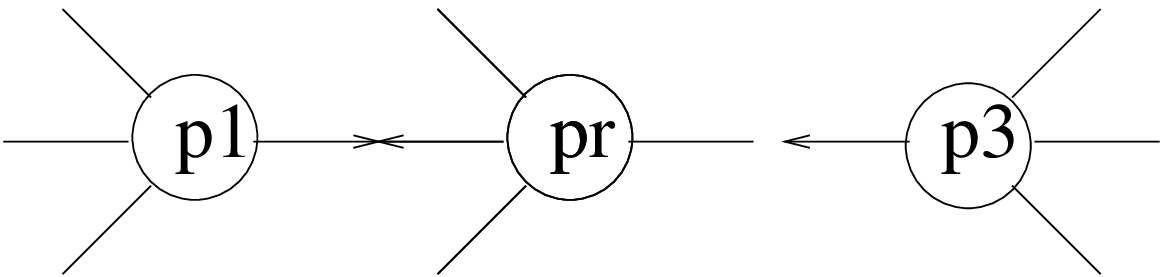}}
\end{center}
\end{minipage}
.
\]
The notation $\overline{p_2}$ stands for $r(p_2)$.

Since we assume that the operad $P$ has a unit,
we can use relation $(R2)$ to express any mating
$p_2 \otimes p_1$ as a trivial mating as follows.
Write $p_1 = p_1 \leftarrow u$,
where $u$ is the unit element of $P$. Then
\begin{equation}  \label{e:tm}
p_2 \otimes p_1 =
p_2 \otimes p_1 \leftarrow u =
p_2 \rightarrow r(p_1) \otimes u.
\end{equation}
This allows us to think of a mating $p_2 \otimes p_1$
roughly
as a reversal of $p_1$ followed by a substitution of
$p_2$ into $p_1$ or vice-versa.

We will refer to a species obtained from a reversible operad
by the above procedure as a \emph{mated species}.
For obvious reasons,
$Q[K]$ is non-empty only if $|K| \geq 2$.

\subsection{How elements of a reversible operad mate}  \label{ss:am}
In examples, substitutions and reversals 
usually have a pictorial description.
And hence so do matings. 
We now show the mating in the associative case.
\[
\begin{minipage}{1.3 in}
\begin{center}
\psfrag{a}{\Huge $a$}
\psfrag{b}{\Huge $b$}
\psfrag{c}{\Huge $c$}
\psfrag{d}{\Huge $d$}
\psfrag{x}{\Huge $x$}
\psfrag{y}{\Huge $y$}
\psfrag{z}{\Huge $z$}
\resizebox{1.2 in}{.6 in}{\includegraphics{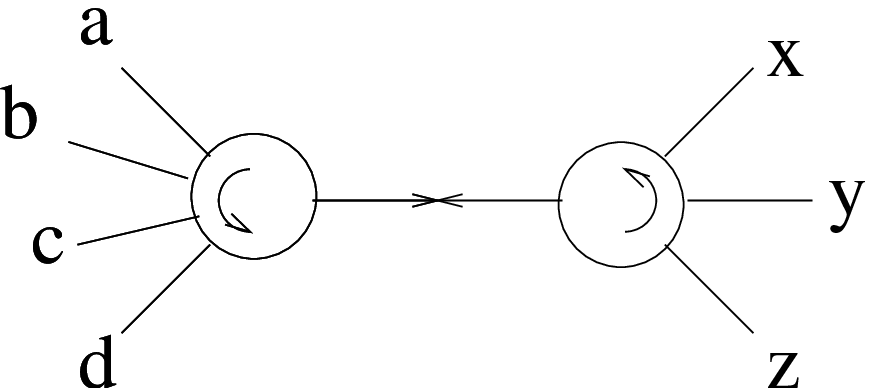}}
\end{center}
\end{minipage}
=
\begin{minipage}{1.1 in}
\begin{center}
\psfrag{a}{\Huge $a$}
\psfrag{b}{\Huge $b$}
\psfrag{c}{\Huge $c$}
\psfrag{d}{\Huge $d$}
\psfrag{x}{\Huge $x$}
\psfrag{y}{\Huge $y$}
\psfrag{z}{\Huge $z$}
\resizebox{1 in}{.6 in}{\includegraphics{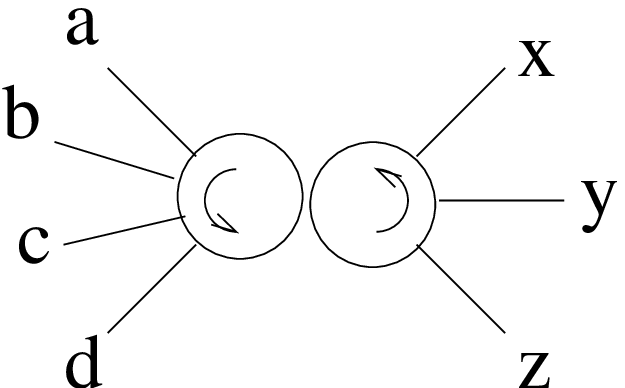}}
\end{center}
\end{minipage}
=
\begin{minipage}{1.1 in}
\begin{center}
\psfrag{a}{\Huge $a$}
\psfrag{b}{\Huge $b$}
\psfrag{c}{\Huge $c$}
\psfrag{d}{\Huge $d$}
\psfrag{x}{\Huge $x$}
\psfrag{y}{\Huge $y$}
\psfrag{z}{\Huge $z$}
\resizebox{1 in}{.6 in}{\includegraphics{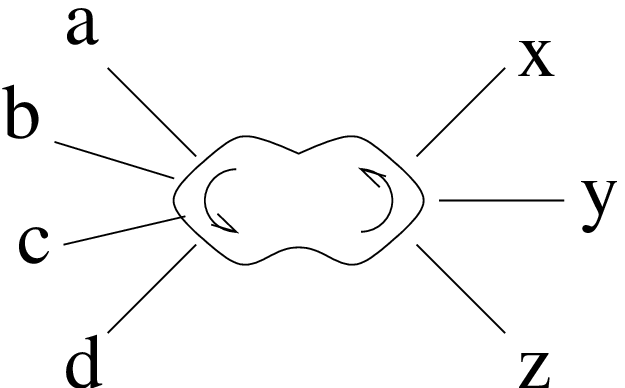}}
\end{center}
\end{minipage}
=
\begin{minipage}{.9 in}
\begin{center}
\psfrag{a}{\Huge $a$}
\psfrag{b}{\Huge $b$}
\psfrag{c}{\Huge $c$}
\psfrag{d}{\Huge $d$}
\psfrag{x}{\Huge $x$}
\psfrag{y}{\Huge $y$}
\psfrag{z}{\Huge $z$}
\resizebox{.8 in}{.6 in}{\includegraphics{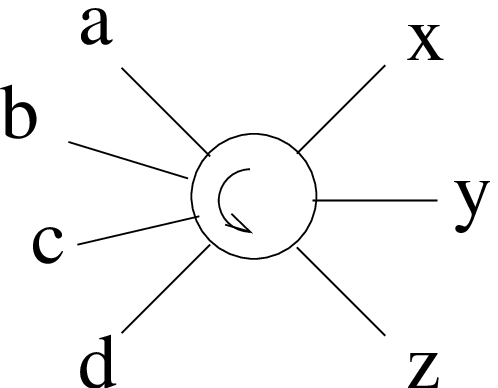}}
\end{center}
\end{minipage}.
\]
The  picture shows that the mating functor maps
the associative operad $a$ to the associative species $aa$.

For the commutative case,
mating occurs as follows.
\[
\begin{minipage}{1.2 in}
\begin{center}
\psfrag{a}{\Huge $a$}
\psfrag{b}{\Huge $b$}
\psfrag{c}{\Huge $c$}
\psfrag{d}{\Huge $d$}
\psfrag{x}{\Huge $x$}
\psfrag{z}{\Huge $z$}
\resizebox{1 in}{.5 in}{\includegraphics{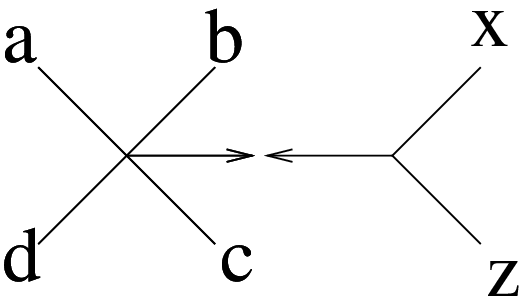}}
\end{center}
\end{minipage}
=
\begin{minipage}{.8 in}
\begin{center}
\psfrag{a}{\Huge $a$}
\psfrag{b}{\Huge $b$}
\psfrag{c}{\Huge $c$}
\psfrag{d}{\Huge $d$}
\psfrag{x}{\Huge $x$}
\psfrag{z}{\Huge $z$}
\resizebox{.5 in}{.5 in}{\includegraphics{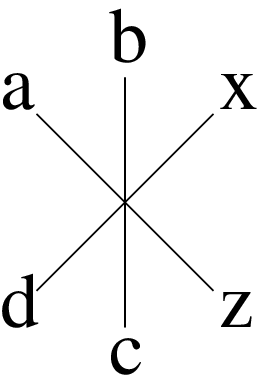}}
\end{center}
\end{minipage}.
\]
Thus the commutative operad $c$ 
maps to the commutative species $cc$.

Similarly, one sees that an operad $P$ 
maps to the species $Q=PP$ when $P=u,t,k$,
see (\ref{ss:eg}).
The letter repetition is used to indicate that
the species is obtained as a mating.
A more correct notation would be to write
$Q=P\reflectbox{$P$}$,
but we will not bother with that.

\subsection{The special role of $u$ and $uu$}  \label{ss:uu}

In all examples so far (\ref{ss:eg}),
we have $P[1]=\Q$.
In other words, 
it is just the span of the unit element $u$ of $P$.
Hence $Q[2]=\Q$ is the span of the element
obtained by mating the unit element of $P$
with itself, which is $uu$.
And so $cc[2]=aa[2]=tt[2]=kk[2]$.
For example,
there is only one cyclic order or
unrooted tree or chord diagram on $2$ letters.
Thus if we ignore the pieces of degree $>$ 1 in $P$
and $>$ 2 in $Q$
then we are left precisely with the unit operad $u$
and species $uu$.

In general, we always assume that $P$ has a unit.
Hence $u$ and $uu$ are always a part of $P$ and $Q$ respectively.
For convenience, we use the letter $u$
for both the unit operad and 
the unit element in an operad.
Similarly for $uu$.

\subsection{The partial derivative}  \label{ss:pd}

Now we introduce the notion of a partial derivative 
for a mated species.

\begin{proposition}  \label{p:pd}
Let $P$ be a reversible operad with unit $u$ and
$Q$ be its mated species.
For $a \in K$ and $q \in Q[K]$,
there is a unique element
$p = \frac{\partial q}{\partial a} \in P[K \setminus a]$
such that $q = \frac{\partial q}{\partial a} \otimes u$. 
In terms of pictures, 
\[
\begin{minipage}{1 in}
\begin{center}
\psfrag{a}{\Huge $a$}
\psfrag{b}{\Huge $b$}
\psfrag{c}{\Huge $c$}
\psfrag{d}{\Huge $d$}
\psfrag{q}{\Huge $q$}
\resizebox{.8 in}{.6 in}{\includegraphics{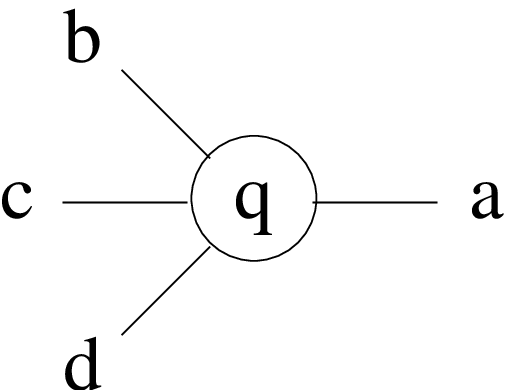}}
\end{center}
\end{minipage}
= 
\begin{minipage}{1.3 in}
\begin{center}
\psfrag{a}{\Huge $a$}
\psfrag{b}{\Huge $b$}
\psfrag{c}{\Huge $c$}
\psfrag{d}{\Huge $d$}
\psfrag{p}{\Huge $p$}
\psfrag{u}{\Huge $u$}
\resizebox{1.1 in}{.6 in}{\includegraphics{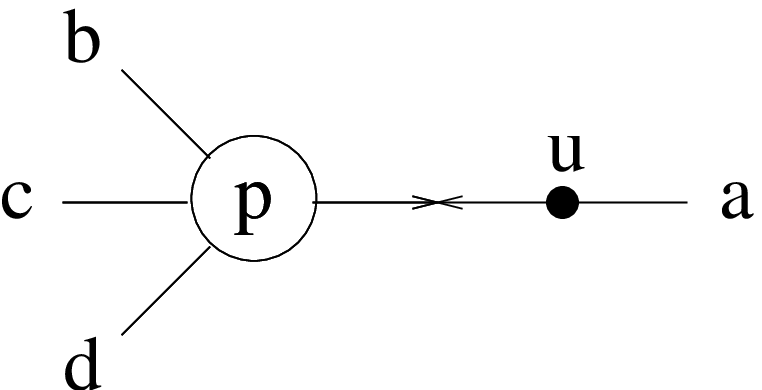}}
\end{center}
\end{minipage}.
\]
\end{proposition}
\begin{proof}
There are two parts to the proposition.
The first one is the existence of
$\frac{\partial q}{\partial a}$.
This is same as saying that $q$ 
can be written as a trivial mating
(at a specified input).
We have already derived this in
equation~\eqref{e:tm}
as a simple consequence of relation $(R2)$,
see (\ref{ss:mf}).
 
The second part is to show the uniqueness of 
$\frac{\partial q}{\partial a}$.
Suppose that $q= p_1 \otimes u = p_2 \otimes u$.
This means that one can obtain
$p_2 \otimes u$ from 
$p_1 \otimes u$ by successive applications of relation $(R2)$.
Now the reduction lemma below says that
this can be done in one step.
This implies that $p_1 = p_2$.
\end{proof}

To complete the proof of the above proposition, 
we prove a reduction lemma
that analyses the relation $(R2)$
in the definition of a mated species (\ref{ss:mf}).
The proof of the lemma will use the two conditions
imposed on the reversal maps
in the definition of a reversible operad (\ref{ss:ro}).

\subsection*{Local notation}

We will follow the convention that
subscripts increase from right to left.
If $p = p_2 \rightarrow p_1$  
then we say that $p_2 \rightarrow p_1$ is a splitting of $p$.
Similarly, $p = p_3 \rightarrow (p_2 \rightarrow p_1) = 
(p_3 \rightarrow p_2) \rightarrow p_1$ 
is a splitting of $p$ into three parts.
Since operad substitution is associative,
we may also write
$p = p_3 \rightarrow p_2 \rightarrow p_1$.
For simplicity of notation, we are suppressing labels.
\vanish{
Now say $p = p_2 \rightarrow p_1$.
Hence
$p \otimes p_0 = (p_2 \rightarrow p_1) \otimes p_0$.
If we apply relation $(R2)$, we get
$p_2 \otimes r(p_1) \leftarrow p_0$.
We say that
$p_2 \otimes r(p_1) \leftarrow p_0$
is obtained from $p \otimes p_0$ by
(splitting $p$ and) applying relation $(R2)$.
}
\begin{reduction} 
Let $P$ be a reversible operad 
with reversal rule $r$ and unit $u$.
The result obtained by two successive applications
of relation $(R2)$ to $p \otimes u$ can, in fact, 
be obtained by a single application of $(R2)$.
\end{reduction}
\begin{proof}
We may assume that the first application of relation $(R2)$ to
$p \otimes u$ involves splitting $p$.
If we split $u$ then nothing changes
and the lemma is proved directly.
Hence say
$p = p^{\prime}_2 \rightarrow p_1$
and applying $(R2)$, we obtain
$$p \otimes u \ = \ p^{\prime}_2 \rightarrow p_1 \otimes u
\ \ \app \ \ p^{\prime}_2 \otimes r(p_1) \leftarrow u \ = \
p^{\prime}_2 \otimes p^{\prime}_1,$$
where $p^{\prime}_1 = r(p_1) \leftarrow u$.

For the second application, we apply relation $(R2)$ to 
$p^{\prime}_2 \otimes p^{\prime}_1$.
This involves two cases.

1. We split $p^{\prime}_2$ as say $p^{\prime}_2 = p_3 \rightarrow p_2$.

2. We split $p^{\prime}_1$.\\
\noindent
Let us look at the first case. We have
$$p^{\prime}_2 \otimes p^{\prime}_1 \ = \ p_3 \rightarrow p_2 \otimes p^{\prime}_1
\ \ \app \ \ p_3 \otimes r(p_2) \leftarrow p^{\prime}_1.$$
Substituting $p^{\prime}_1 = r(p_1) \leftarrow u$
and using the associativity of substitution and
the compatibility of reversal with substitution,
we get
$$p_3 \otimes r(p_2) \leftarrow p^{\prime}_1 \ = \
p_3 \otimes r(p_2) \leftarrow (r(p_1) \leftarrow u) \ = \
p_3 \otimes r(p_2 \rightarrow p_1) \leftarrow u.$$
Hence the net effect of two applications of relation $(R2)$ has been
$$p \otimes u \ \ \rightsquigarrow \ \
p_3 \otimes r(p_2 \rightarrow p_1) \leftarrow u,$$
where
$p = p^{\prime}_2 \rightarrow p_1 =
(p_3 \rightarrow p_2) \rightarrow p_1 =
p_3 \rightarrow (p_2 \rightarrow p_1)$.
It is clear that the same effect is achieved by
$$p \otimes u \ = \ p_3 \rightarrow (p_2 \rightarrow p_1) \otimes u
\ \ \app \ \ p_3 \otimes r(p_2 \rightarrow p_1) \leftarrow u,$$
which is just one application of relation $(R2)$
to $p \otimes u$.
This completes the analysis for the first case.

For the second case, we just point out
a simple subcase and leave the rest out.
If we split $p^{\prime}_1$ as $p^{\prime}_1 = r(p_1) \leftarrow u$
and apply relation $(R2)$ then
we just reverse the first step and get back $p \otimes u$.
And here we used the relation
$r_{y,x} \circ r_{x,y} = \ \text{identity}$.

\end{proof}

\section{Overview of symplectic operad geometry} \label{s:sgo}

In this section, we explain
some of the philosophy of symplectic operad geometry.
The reader can just glance through it on a first reading.
Complete details are given in 
Sections~\ref{s:cm}-\ref{s:me}.
The reader more interested in graph homology
may directly go to Sections~\ref{s:gh}-\ref{s:mgh}.

\subsection{What is an operad manifold?} \label{ss:om}

An ordinary manifold $X$ is a commutative object 
in the following sense.
Functions on $X$ form a commutative algebra,
the product being pointwise multiplication.
To go a little further,
differential forms on $X$ form 
a graded or supercommutative algebra
with the usual wedge product.
Many geometric notions associated to $X$
can be captured by these commutative algebras.
For example, vector fields on $X$
are derivations of the algebra of functions on $X$.
Among all manifolds, 
$\R^n$ plays a special role.
And the space of polynomial functions on $\R^n$
is the free commutative algebra on $n$ generators.

Now for a general operad $P$,
there is a notion of a $P$-algebra, a $P$-superalgebra, 
a free $P$-algebra, etc.
For example,
for the associative operad, 
one has associative algebras;
for the Lie operad, one has Lie algebras and so on.
We would like to think of a $P$-algebra 
as the algebra of functions on a ``$P$-manifold'',
and of the free $P$-algebra on $n$ generators
as the algebra of functions on 
the $n$ dimensional ``$P$-manifold''
which is the analogue of $\R^n$.
We emphasise that a ``$P$-manifold'' as a geometric object 
does not make any sense.
However, thinking of a $P$-algebra in this manner
is useful because it allows us to make analogies
with the commutative case.
 
In this paper, we will only deal with free $P$-algebras.
 
\subsection{What is symplectic operad geometry?} \label{s:sog}

Let us start with the commutative case.
For an introduction to classical symplectic topology, 
see the book by McDuff and Salamon~\cite{\mcduff}.
A symplectic manifold $(X,\omega)$ is a manifold $X$
with a closed non-degenerate 2 form $\omega$.
The existence of such a 2 form implies that $X$
is even dimensional.
Further, the algebra of smooth functions
$C^{\infty}(X)$ has the structure of a Lie algebra.
The standard example of a symplectic manifold is
$(\R^{2n}, \omega_0)$.
Here $\omega_0$ is the standard symplectic form
$\sum dp_i \wedge dq_i$,
where $p_1,\ldots,p_n,q_1,\ldots,q_n$
are coordinates on $\R^{2n}$.
The Lie algebra structure on $C^{\infty}(\R^{2n})$
is given by the usual Poisson bracket
\begin{equation} \label{e:pb}
\{F,H\}=\sum_{i=1}^n \frac{\partial F}{\partial p_i}\frac{\partial H}{\partial q_i} - \frac{\partial F}{\partial q_i}\frac{\partial H}{\partial
p_i} \ \ \text{for} \ \ F,H \in C^{\infty}(\R^{2n}).
\end{equation}
Instead of smooth functions, 
it is easier to deal with polynomial functions
and henceforth we will always do so.

\subsubsection{The general case}
We want to do something similar for a general operad $P$.
Let $X$ be a ``$P$-manifold''.
As explained in (\ref{ss:om}),
this simply means that we have a $P$-algebra.
And since we are only going to deal with the free case,
we further assume that we have a free $P$-algebra.
Call it $PA$.
Since one thinks of $PA$ as polynomial functions on $X$,
one would expect a Lie algebra structure on $PA$,
in analogy with the commutative case.
However, this is not true.

It turns out that one needs to consider 
another algebraic object $QA$
(Section~\ref{s:cm}).
And this object can be constructed only if $P$ is reversible.
Then if $PA$ is free on an even number of generators,
it is true that $QA$ is a Lie algebra.
The two important operations are 
$$\frac{\partial }{\partial x_i} : QA \longrightarrow PA 
\ \ \ \text{and} \ \ \ 
\{\ ,\ \} : QA \otimes QA \longrightarrow QA.$$
The Poisson bracket $\{\ ,\ \}$ equips $QA$ 
with the structure of a Lie algebra.
It is given by a formula that is almost identical 
to equation~\eqref{e:pb}, see 
equation~\eqref{e:ptb} in (\ref{ss:mate}).

Just as $PA$ is related to the operad $P$,
the space $QA$ is related to an object $Q$.
The object $Q$ is not an operad but a simpler object called a species.
To give an analogy, operads are like algebras and
species are like vector spaces.
The species $Q$ is the image of $P$ under the mating functor
constructed in Section~\ref{s:mf}.
The space $QA$ can be directly defined 
in terms of $Q$.

So far, we have not talked about any symplectic form
in the general case.
The point is that the Lie structure on $QA$ 
emerges naturally from a symplectic form.
This is explained in Section~\ref{s:st}.
A ``symplectic $P$-manifold'' must be a ``$P$-manifold''
which has a symplectic form.
The only example we give is that of a free $P$-algebra
on an even number of generators, for $P$ reversible.
When specialised to the commutative case,
this says that we are only considering the symplectic manifolds
$(\R^{2n}, \omega_0)$.

\subsubsection{Back to the classical case}  \label{ss:tcc}

We now explain how the general discussion above
specialises to the commutative situation.
That is, $P$ is the commutative operad $c$, see (\ref{sss:c}).
For $P=c$, we have the free (non-unital) $P$-algebra
$PA=$ polynomial functions with no constant terms
on the symplectic manifold $(\R^{2n}, \omega_0)$.
The mated species $Q$ in this case,
which is $cc$, is almost the same.
Hence the corresponding algebraic object $QA$ is also similar.
We have $QA=$ polynomial functions in $2n$ variables 
with no constant or linear terms. 
In this case, the partial derivative has the usual meaning 
and the Poisson bracket on $QA$ is given by the formula
in equation~\eqref{e:pb}.
Thus the content of the general discussion in quite simple
in this case.

We now give a pictorial description of this calculus. 
If $F \in QA$ is a monomial, say
$F = x_1^2 x_2 x_3$,
then we represent it as 
$F = 
\mbox{
\begin{minipage}{.4 in}
\begin{center}
\psfrag{a}{\Huge $x_2$}
\psfrag{b}{\Huge $x_1$}
\psfrag{c}{\Huge $x_3$}
\psfrag{d}{\Huge $x_1$}
\resizebox{.4 in}{.4 in}{\includegraphics{cc.eps}}
\end{center}
\end{minipage} 
}
$.
And if $F$ is a polynomial rather than a monomial
then we represent it as a formal sum of pictures.
Though $PA$ is almost the same as $QA$,
we will represent its elements slightly differently.
Namely, if $x_1^2 x_2 x_3 \in PA$
then we write it as
$F = 
\mbox{
\begin{minipage}{.4 in}
\begin{center}
\psfrag{a}{\Huge $x_2$}
\psfrag{b}{\Huge $x_1$}
\psfrag{c}{\Huge $x_3$}
\psfrag{d}{\Huge $x_1$}
\resizebox{.4 in}{.4 in}{\includegraphics{csubt.eps}}
\end{center}
\end{minipage} 
}
$.

We explain how the partial derivative
$\frac{\partial }{\partial x_i} : QA \longrightarrow PA$ 
works by an example.
\[
\frac{\partial}{\partial x_1}
\left(
\begin{minipage}{.6 in}
\begin{center}
\psfrag{a}{\Huge $x_2$}
\psfrag{b}{\Huge $x_1$}
\psfrag{c}{\Huge $x_3$}
\psfrag{d}{\Huge $x_1$}
\resizebox{.4 in}{.4 in}{\includegraphics{cc.eps}}
\end{center}
\end{minipage} 
\right) = 
\begin{minipage}{.6 in}
\begin{center}
\psfrag{a}{\Huge $x_2$}
\psfrag{b}{\Huge $$}
\psfrag{c}{\Huge $x_3$}
\psfrag{d}{\Huge $x_1$}
\resizebox{.4 in}{.4 in}{\includegraphics{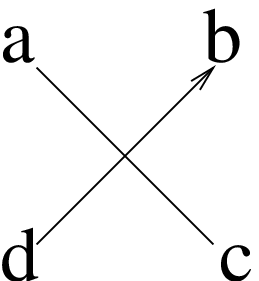}}
\end{center}
\end{minipage} 
\ \ + \ \
\begin{minipage}{.6 in}
\begin{center}
\psfrag{a}{\Huge $x_2$}
\psfrag{b}{\Huge $x_1$}
\psfrag{c}{\Huge $x_3$}
\psfrag{d}{\Huge $$}
\resizebox{.4 in}{.4 in}{\includegraphics{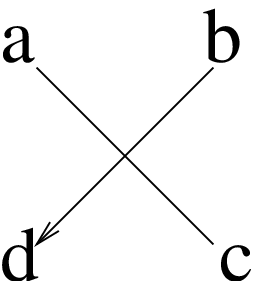}}
\end{center}
\end{minipage} 
.
\]
In other words,
cut all inputs with label $x_1$, one at a time.
The above pictorial equation just says that
$\frac{\partial}{\partial x_1}
(x_1^2 x_2 x_3) = 2 x_1 x_2 x_3$.

Now we illustrate the Poisson bracket on $QA$
by an example.
\[
\left\{
\begin{minipage}{.6 in}
\begin{center}
\psfrag{x}{\Huge $p_1$}
\psfrag{y}{\Huge $q_1$}
\psfrag{z}{\Huge $p_2$}
\resizebox{.5 in}{.4 in}{\includegraphics{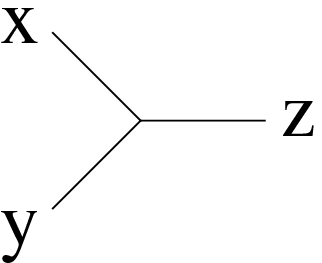}}
\end{center}
\end{minipage}
\ \ , \ \
\begin{minipage}{.8 in}
\begin{center}
\psfrag{p1}{\Huge $q_2$}
\psfrag{p2}{\Huge $p_2$}
\psfrag{q1}{\Huge $q_1$}
\psfrag{q2}{\Huge $p_2$}
\psfrag{1}{}
\psfrag{2}{}
\psfrag{uu}{}
\resizebox{.7 in}{.15 in}{\includegraphics{sf.eps}}
\end{center}
\end{minipage}
\right\}
\ \ = \ \
\begin{minipage}{1 in}
\begin{center}
\psfrag{x}{\Huge $p_1$}
\psfrag{y}{\Huge $q_1$}
\psfrag{z}{\Huge $p_2$}
\psfrag{p2}{\Huge $p_2$}
\resizebox{1 in}{.4 in}{\includegraphics{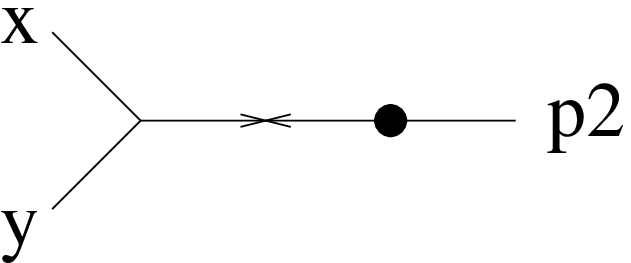}}
\end{center}
\end{minipage}
\ \ = \ \
\begin{minipage}{.6 in}
\begin{center}
\psfrag{x}{\Huge $p_1$}
\psfrag{y}{\Huge $q_1$}
\psfrag{z}{\Huge $p_2$}
\resizebox{.5 in}{.4 in}{\includegraphics{cca.eps}}
\end{center}
\end{minipage}
.
\]
In usual terms, this says that
$\{p_1 q_1 p_2, p_2 q_2 \} = p_1 q_1 p_2$.
So to compute $\{F,H\}$,
cut a $p_i$ from $F$, and
a $q_i$ from $H$ (or vice-versa),
do a mating and then sum over all possibilities.
Compare the picture above
with the mating suggested for the commutative operad 
in (\ref{ss:am}).

The pictorial calculus may look silly but
the point is that it generalises very nicely
to reversible operads (Section~\ref{s:cm}).
The pictorial way of thinking can also be extended 
to the computation of the homology of the Lie algebra $QA$.
This involves the building of a chain complex of graphs,
which leads to graph homology.
This connection is made in the first part of the proof
of the main theorem (Section~\ref{s:fin}).

We now record these ideas in a nutshell as follows.
\medskip

\begin{center}
\begin{tabular}{|c|c|} \hline 
Pictorial & Formal \\ \hline \hline
Cutting of a vertex & Partial derivative \\ \hline
Mating two cut vertices & Poisson bracket \\ \hline
Building a graph & Homology \\ \hline 
\end{tabular}
\end{center}

\section{Cuttings and Matings}  \label{s:cm}

The goal of this section
is to define the Lie algebra $QA_{\infty}$
that occurs in the left hand side of the main theorem.
Let $V$ be the vector space with basis
$x_1, \ldots,x_n$.
Later to do symplectic theory,
we will use the basis
$p_1, \ldots,p_n,q_1, \ldots,q_n$.
If we want to show the dependence of $n$ explicitly,
we will write $V_n$ for $V$.

Let $P$ be a reversible operad and $Q$ its mated species
obtained by applying the mating functor (Section~\ref{s:mf}).
We will first define the algebraic objects $QA$ and $PA$.
Using a calculus of cuttings and matings,
we will give $QA$ a Lie algebra structure.
The example of the commutative operad $(P=c)$
was discussed in (\ref{ss:tcc}).
A trivial looking but important example
is the unit operad $(P=u)$.
The Lie algebra in this case is 
the usual symplectic Lie algebra $\fsp$.
This is explained in (\ref{ss:sl}).

As in the case of the underlying vector space $V$, 
we will sometimes write $QA_n$ and $PA_n$,
to show the dependence of $n$ explicitly.
The Lie algebra $QA_{\infty}$ will be defined
as a direct limit $\underset{n \to \infty}{\lim} QA_n$.
After this, we will briefly explain 
the Hopf algebra structure on $H_*(QA_{\infty})$.

\subsection{The free algebras} \label{ss:fa}

We now define the free algebra $QA$ of a species $Q$
on generators $x_1, \ldots,x_n$.
This also includes the case of $PA$,
since every operad is a species.

\medskip
\begin{center}
\begin{tabular}{|c|c|c|c|} \hline 
Operad $P$ & An element of $PA$ & 
Species $Q$ & An element of $QA$ \\ \hline \hline
$c$ & 
\begin{minipage}{.4 in}
\vspace{.1in}
\begin{center}
\psfrag{a}{\Huge $x_2$}
\psfrag{b}{\Huge $x_1$}
\psfrag{c}{\Huge $x_3$}
\psfrag{d}{\Huge $x_1$}
\resizebox{.4 in}{.4 in}{\includegraphics{csubt.eps}}
\end{center}
\vspace{.1in}
\end{minipage} &
$cc$ &
\begin{minipage}{.4 in}
\vspace{.1in}
\begin{center}
\psfrag{a}{\Huge $x_2$}
\psfrag{b}{\Huge $x_1$}
\psfrag{c}{\Huge $x_3$}
\psfrag{d}{\Huge $x_1$}
\resizebox{.4 in}{.4 in}{\includegraphics{cc.eps}}
\end{center}
\vspace{.1in}
\end{minipage} \\ \hline
$a$ &  
\begin{minipage}{.4 in}
\vspace{.1in}
\begin{center}
\psfrag{a}{\Huge $x_2$}
\psfrag{b}{\Huge $x_1$}
\psfrag{c}{\Huge $x_3$}
\psfrag{d}{\Huge $x_1$}
\resizebox{.4 in}{.4 in}{\includegraphics{asubt.eps}}
\end{center}
\vspace{.1in}
\end{minipage} &
$aa$ &
\begin{minipage}{.4 in}
\vspace{.1in}
\begin{center}
\psfrag{a}{\Huge $x_2$}
\psfrag{b}{\Huge $x_1$}
\psfrag{c}{\Huge $x_3$}
\psfrag{d}{\Huge $x_1$}
\resizebox{.4 in}{.4 in}{\includegraphics{aa.eps}}
\end{center}
\vspace{.1in}
\end{minipage} \\ \hline
$t$ &
\begin{minipage}{.6 in}
\vspace{.1in}
\begin{center}
\psfrag{a}{\Huge $x_2$}
\psfrag{b}{\Huge $x_1$}
\psfrag{c}{\Huge $x_3$}
\psfrag{d}{\Huge $x_1$}
\resizebox{.6 in}{.4 in}{\includegraphics{tsubt.eps}}
\end{center}
\vspace{.1in}
\end{minipage} &
$tt$ &
\begin{minipage}{.6 in}
\vspace{.1in}
\begin{center}
\psfrag{a}{\Huge $x_2$}
\psfrag{b}{\Huge $x_1$}
\psfrag{c}{\Huge $x_3$}
\psfrag{d}{\Huge $x_1$}
\resizebox{.6 in}{.4 in}{\includegraphics{tt.eps}}
\end{center}
\vspace{.1in}
\end{minipage} \\ \hline
\end{tabular}
\end{center}
\vspace{0.3 cm}
\noindent
To get an element of $QA$,
take any element of $Q[I]$,
for some finite non-empty set $I$
and replace each element of $I$ by one of the generators
$x_1, \ldots,x_n$.
Call such an element a monomial.
To get a general element of $QA$,
take linear combinations of monomials.
\[
\begin{minipage}{.8 in}
\begin{center}
\psfrag{a}{\Huge $a$}
\psfrag{b}{\Huge $b$}
\psfrag{c}{\Huge $c$}
\psfrag{d}{\Huge $d$}
\resizebox{.6 in}{.6 in}{\includegraphics{q.eps}}
\end{center}
\end{minipage}
\in Q[I] \qquad \rightsquigarrow \qquad
\begin{minipage}{.8 in}
\begin{center}
\psfrag{x1}{\Huge $x_1$}
\psfrag{x2}{\Huge $x_2$}
\psfrag{x4}{\Huge $x_4$}
\resizebox{.8 in}{.6 in}{\includegraphics{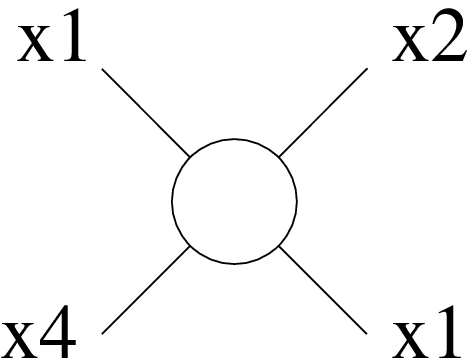}}
\end{center}
\end{minipage}
\in QA^4.
\]
We may also write this as
$QA = \bigoplus_{j\geq 2} (Q[j] \otimes V^{\otimes j})_{\Sigma_j}$,
where $\Sigma_j$ acts on $Q[j]$
by permuting letters and
on $V^{\otimes j}$ by permuting the tensor factors.
Strictly speaking, the grading starts in degree $1$.
But since we are mainly interested in mated species,
we started the grading in degree $2$.
We will sometimes write $QA^j$ for the degree $j$ piece.
Though we call $QA$ the free algebra of $Q$,
it only has the structure of a graded vector space.

If the species $Q$ has the structure of an operad $P$
then we emphasise it by drawing instead
\[
\begin{minipage}{.8 in}
\begin{center}
\psfrag{a}{\Huge $a$}
\psfrag{b}{\Huge $b$}
\psfrag{c}{\Huge $c$}
\psfrag{d}{\Huge $d$}
\resizebox{.7 in}{.6 in}{\includegraphics{p.eps}}
\end{center}
\end{minipage}
\in P[I] \qquad \rightsquigarrow \qquad
\begin{minipage}{.8 in}
\begin{center}
\psfrag{x1}{\Huge $x_1$}
\psfrag{x2}{\Huge $x_2$}
\psfrag{x4}{\Huge $x_4$}
\resizebox{.7 in}{.6 in}{\includegraphics{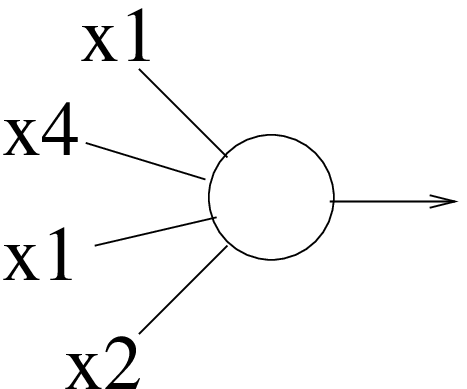}}
\end{center}
\end{minipage}
\in PA^4.
\]
Just as for a species, we get
$PA = \bigoplus_{j\geq 1} (P[j] \otimes V^{\otimes j})_{\Sigma_j}.$
However, the situation for operads is more interesting.
The substitution rule for operads allows us
to define binary (and higher) products on $PA$.
Hence $PA$ looks more like the algebras 
that we are used to.
For example,
$cA=$ free non-unital commutative algebra and
$aA=$ free non-unital associative algebra
on the generators $x_1, \ldots,x_n$.

As this discussion suggests,
it is possible to give a more abstract definition of $PA$.\
To every operad $P$, one can associate a category of $P$-algebras
(\ref{ss:om});
the object $PA$ is then the free object,
with respect to $V$, in this category.

\begin{remark}
Recall that, roughly speaking, 
the mated species
$Q$ was defined as a quotient of $P \otimes P$,
see (\ref{ss:mf}).
In the same way, one can define
$QA$ as $PA \otimes PA$ subject to the two relations:
\begin{equation}  \label{e:r}
a \otimes b = b \otimes a  \quad \text{and} \quad
a \rightarrow b \otimes c = a \otimes r(b) \leftarrow c.
\end{equation}
\end{remark}

\subsubsection*{Examples}  

Refer to the table in (\ref{ss:fa}).
In the commutative case, the pictures for $c$ and $cc$
contain the same amount of information,
they represent a monomial in the commutative variables
$x_1,x_2,\ldots,x_n$.
In the associative case, the pictures for $a$ and $aa$
are different.
They represent linear and cyclic orders respectively
in non-commutative variables
$x_1,x_2,\ldots,x_n$.
In the tree case, the objects $PA$ and $QA$ are 
best described using pictures rather than words.

Following the analogy with the commutative case
discussed in (\ref{ss:tcc}),
we now discuss the calculus of cuttings and matings 
in the general case.
The reader may keep the above examples in mind.

\subsection{Cutting}  \label{ss:cut}
We define
$\frac{\partial }{\partial x_i} : QA \longrightarrow PA$
by showing how it works on a schematic example.
\[
\begin{minipage}{.9 in}
\begin{center}
\psfrag{x1}{\Huge $x_1$}
\psfrag{x2}{\Huge $x_2$}
\psfrag{x4}{\Huge $x_4$}
\resizebox{.8 in}{.6 in}{\includegraphics{qa.eps}}
\end{center}
\end{minipage}
\qquad \buildrel \frac{\partial}{\partial x_1}\over \longmapsto \qquad
\begin{minipage}{.9 in}
\begin{center}
\psfrag{x1}{\Huge $x_1$}
\psfrag{x2}{\Huge $x_2$}
\psfrag{x4}{\Huge $x_4$}
\resizebox{.8 in}{.6 in}{\includegraphics{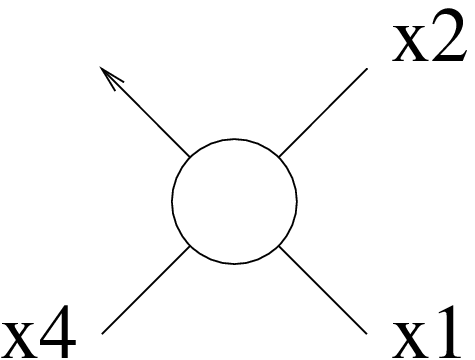}}
\end{center}
\end{minipage}
\ \ + \ \
\begin{minipage}{.9 in}
\begin{center}
\psfrag{x1}{\Huge $x_1$}
\psfrag{x2}{\Huge $x_2$}
\psfrag{x4}{\Huge $x_4$}
\resizebox{.8 in}{.6 in}{\includegraphics{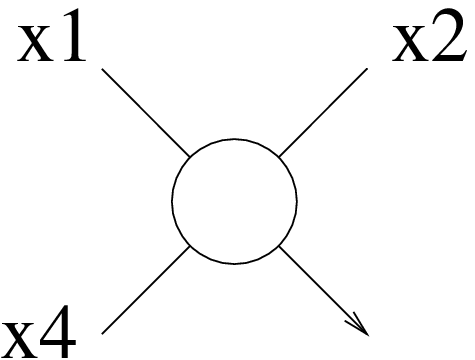}}
\end{center}
\end{minipage}
.
\]
Namely, to define
$\frac{\partial}{\partial x_1}$,
cut all inputs with label $x_i$,
one at a time.
The justification for such a definition 
is given by Proposition~\ref{p:pd}
and the above remark.

As an example, in the associative case,
the derivative with respect to $x_1$ of
the cyclic word $x_1 x_4 x_1 x_2$ is the sum of
two linear words, $x_4 x_1 x_2$ and $x_2 x_1 x_4$.

\subsection{Mating}  \label{ss:mate}
In order to perform matings, 
one requires an even number of characters.
So assume that
$V$ is even dimensional and has a basis
$p_1, \ldots,p_n,q_1, \ldots,q_n$.
Define the Poisson bracket
$\{\ ,\ \} : QA \otimes QA \longrightarrow QA$
using the formula
\begin{equation} \label{e:ptb}
\{F,H\}=\sum_{i=1}^n \frac{\partial F}{\partial p_i} \otimes \frac{\partial H}{\partial q_i} - 
\frac{\partial F}{\partial q_i} \otimes \frac{\partial H}{\partial p_i}.
\end{equation}
In (\ref{ss:am}), we explained how the tensor sign 
can be interpreted as a mating.
The only difference now is that the objects being mated,
instead of being labelled by elements of a set $I$,
are labelled by elements of a vector space $V$.
So one can continue with the same interpretation.
\[
\left\{
\begin{minipage}{.8 in}
\begin{center}
\psfrag{a}{\Huge $p_2$}
\psfrag{c}{\Huge $p_3$}
\psfrag{d}{\Huge $q_3$}
\psfrag{b}{\Huge $q_1$}
\resizebox{.7 in}{.6 in}{\includegraphics{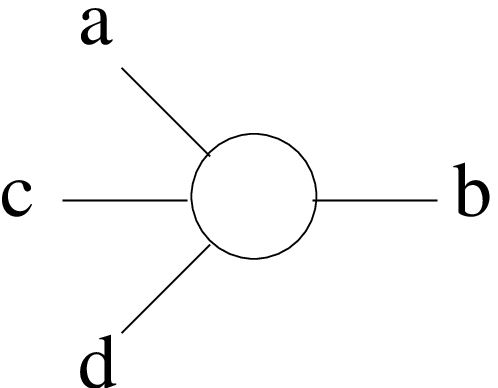}}
\end{center}
\end{minipage}
\ \ , \ \
\begin{minipage}{.9 in}
\begin{center}
\psfrag{p1}{\Huge $p_1$}
\psfrag{p2}{\Huge $p_2$}
\psfrag{q1}{\Huge $q_1$}
\psfrag{q2}{\Huge $q_2$}
\resizebox{.8 in}{.6 in}{\includegraphics{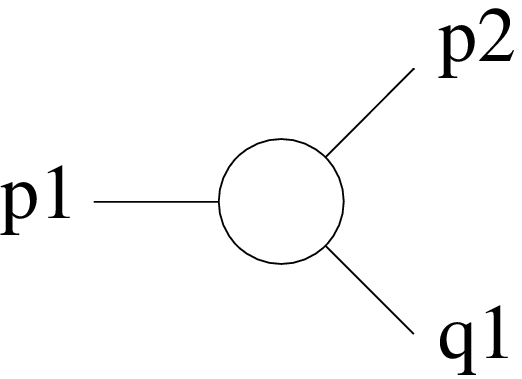}}
\end{center}
\end{minipage}
\right\}
\ \ = \ \
-
\begin{minipage}{1.4 in}
\begin{center}
\psfrag{a}{\Huge $p_2$}
\psfrag{c}{\Huge $p_3$}
\psfrag{d}{\Huge $q_3$}
\psfrag{x}{\Huge $p_2$}
\psfrag{z}{\Huge $q_1$}
\resizebox{1.2 in}{.6 in}{\includegraphics{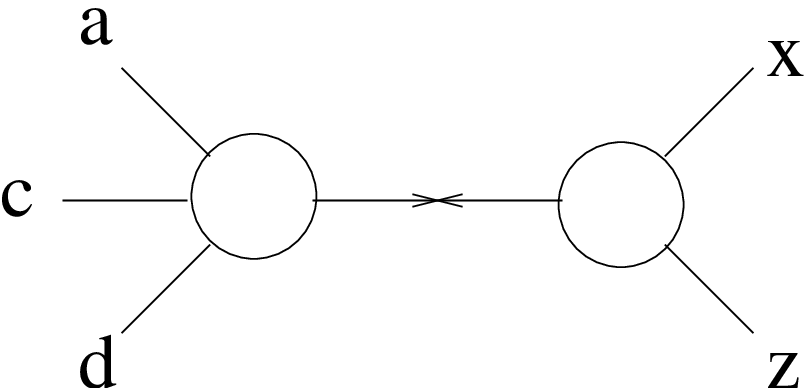}}
\end{center}
\end{minipage}
.
\]
In the above example,
there is only one mating possible.
So there is only one term on the right hand side.
The minus sign reflects the fact that 
the ``$q$'' was cut from the first term and
the ``$p$'' from the second term.

\begin{proposition}  \label{p:pb}
The bracket $\{\ ,\ \}$ equips $QA$ 
with the structure of a Lie algebra.
\end{proposition}

\noindent
This can be verified directly.
It will also follow from the symplectic operad theory
that we will present in Section~\ref{s:st}.
So we assume this result for the moment.
Note that the Lie algebra $QA$ is graded
and the bracket $\{\ ,\ \}$ is a map of degree $-2$.

\begin{remark}
If the reversible operad $P$ is based on sets
then $QA$ also has a Lie coalgebra structure 
defined using breakups (rather than matings).
This parallels the definition of the coboundary operator
on graphs (\ref{ss:stacb}).
\end{remark}

\subsection{Connection to the symplectic Lie algebra $\fsp$}  \label{ss:sl}

We now look at an important example
that is already present in the classical theory.
This is the simplest case, i.e.,
$P=u$ and $Q=uu$,
see (\ref{ss:uu}).
It is easy to see that $uA_n=V$, 
the underlying vector space and
$uuA_n=$ homogeneous polynomials of degree 2
in the variables $p_1, \ldots,p_n,q_1, \ldots,q_n$.
In this case, $\{\ ,\ \}$ is the usual Poisson bracket
with the tensor sign in equation~\eqref{e:ptb}
being ordinary multiplication.

Going to the general case,
note that $uuA$ is always a subalgebra of $QA$ 
that sits inside $QA^2$, the piece of degree 2.
This is because $uu$ is a subspecies of $Q$,
see (\ref{ss:uu}).
Moreover, since the bracket $\{\ ,\ \}$ has degree $-2$,
it restricts to a map
$QA^j \otimes uuA \longrightarrow QA^j$.
Hence each graded piece $QA^j$ is a right $uuA$ module.
As an example,
\[
\left\{
\begin{minipage}{.9 in}
\begin{center}
\psfrag{p1}{\Huge $p_1$}
\psfrag{p2}{\Huge $p_2$}
\psfrag{q1}{\Huge $q_1$}
\psfrag{q2}{\Huge $q_2$}
\resizebox{.8 in}{.6 in}{\includegraphics{qas.eps}}
\end{center}
\end{minipage}
\ \ , \ \
\begin{minipage}{.8 in}
\begin{center}
\psfrag{p1}{\Huge $p_1$}
\psfrag{p2}{\Huge $p_2$}
\psfrag{q1}{\Huge $q_1$}
\psfrag{q2}{\Huge $q_2$}
\psfrag{1}{}
\psfrag{2}{}
\psfrag{uu}{}
\resizebox{.7 in}{.15 in}{\includegraphics{sf.eps}}
\end{center}
\end{minipage}
\right\}
\ \ = \ \
\begin{minipage}{.9 in}
\begin{center}
\psfrag{p1}{\Huge $p_1$}
\psfrag{p2}{\Huge $p_2$}
\psfrag{q1}{\Huge $q_1$}
\psfrag{q2}{\Huge $q_2$}
\resizebox{.8 in}{.6 in}{\includegraphics{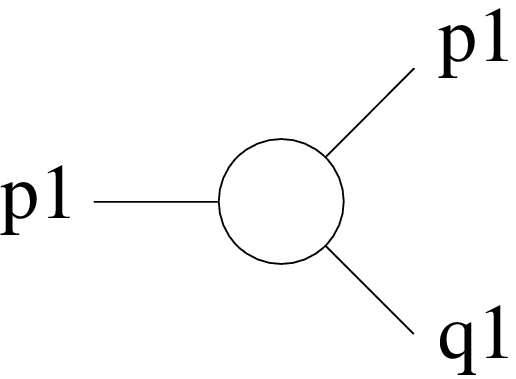}}
\end{center}
\end{minipage}
\ \ - \ \
\begin{minipage}{.9 in}
\begin{center}
\psfrag{p1}{\Huge $p_1$}
\psfrag{p2}{\Huge $p_2$}
\psfrag{q1}{\Huge $q_1$}
\psfrag{q2}{\Huge $q_2$}
\resizebox{.8 in}{.6 in}{\includegraphics{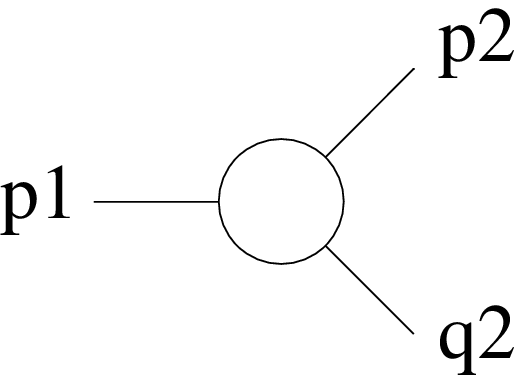}}
\end{center}
\end{minipage}
.
\]
Observe that the action only involves changing labels.
It does not change the internal structure 
of the vertex.
This is because the matings are all trivial.

The symplectic Lie algebra $\fsp$ is defined
as the space of linear maps on $V$
that kill the symplectic form
$\sum dp_i \wedge dq_i$
(the action is as a derivation);
the bracket being the usual commutator.
Note that
$QA^j = (Q[j] \otimes V^{\otimes j})_{\Sigma_j}$
is a left $\fsp$ module with the action 
induced by the left action on $V$ and
the trivial action on $Q[j]$.

We now explain the relation between the Lie algebras
$uuA_n$ and $\fsp$.

\begin{proposition}  \label{p:degtwo}
There is a Lie algebra anti-isomorphism
$uuA_n \to \fsp$ given by $H \mapsto \xi_H$,
where
$\xi_H(p_i)=\frac{\partial H}{\partial q_i}$ and
$\xi_H(q_i)=-\frac{\partial H}{\partial p_i}$.
\end{proposition}

\noindent
Again this can be checked directly
but will be a consequence of the general theory
(Section~\ref{s:st}).
The definition of the module structures on $QA$
together with the above proposition gives us the following.

\begin{corollary}
The right $uuA_n$ and the left $\fsp$ 
module structures on $QA^j$ are compatible
via the above anti-isomorphism $uuA_n \to \fsp$.
In other words,
$\{F,H\} = \xi_H \cdot F$ for $F \in QA^j$
and $H \in uuA_n$.
\end{corollary}

\subsection{The Lie algebra $QA_{\infty}$}
\label{ss:ha}

We have the underlying vector space 
$V_{n} \subset V_{n+1}$, with
the basis of $V_{n+1}$ extending the basis on $V_n$ to
$p_1, \ldots,p_{n+1},q_1, \ldots,q_{n+1}$.
This gives us a sequence of Lie algebra inclusions
$$QA_1 \subset \ldots \subset QA_n \subset QA_{n+1} \subset \ldots$$
We denote the direct limit by $QA_{\infty}$.

A Hopf algebra is a vector space 
which has compatible algebra and coalgebra structures
and an inverse which is called the antipode.
For basic information on Hopf algebras,
see~\cite{\sweedler,\kassel}.
The graded vector space $H_*(QA_{\infty})$
has the structure of a Hopf algebra,
which we now explain.

To describe the product on $H_*(QA_{\infty})$,
we start with a morphism of Lie algebras
$QA_n \oplus QA_m \into QA_{n+m}$, defined by
$(F,H) \mapsto F+H$,
where we shift up the indices of $H$ by $n$.
In other words, we think of
$QA_n$ and $QA_m$ as disjoint Lie subalgebras of $QA_{n+m}$.
This defines a map
$H_*(QA_{n}) \otimes H_*(QA_{m}) \to H_*(QA_{n+m})$.
Taking direct limits, we get a product on $H_*(QA_{\infty})$.
We point out that the diagram
\begin{equation*}
\begin{CD}
QA_{n} \oplus QA_{m} @>>> QA_{n+m} \\
@VVV @VVV \\
QA_{n+1} \oplus QA_{m} @>>> QA_{n+m+1} \\
\end{CD}
\end{equation*}
commutes only upto index shifting.
However, this cannot be detected on the homology level.

The coproduct is relatively straightforward to define.
We start with the diagonal map
$QA_{\infty} \to QA_{\infty} \oplus QA_{\infty}$,
which is a morphism of Lie algebras.
This induces a cocommutative coproduct on
$H_*(QA_{\infty})$.

\section{Symplectic operad theory}  \label{s:st}

In this section, we give a self-contained treatment of 
Kontsevich's symplectic mini-theory.
But in order to appreciate it,
it is important to be familiar with
classical symplectic theory.
This involves the principles of Hamiltonian mechanics
\cite{\mcduff} and
basic notions of differential topology such as
vector fields, differential forms,
Lie derivatives and contraction operators
\cite{\boothby,\spivak}.

For some recent ideas and applications 
in the associative case, see the papers of
Ginzburg~\cite{\ginz} and Bocklandt and Le Bruyn~\cite{\bl}. 
For a detailed exposition,
see Chapter $7$ of Le Bruyn's book ``Non-commutative geometry'',
which is available on his homepage.

The basic objects of interest are summarised in the table.
We have already seen $PA$ and $QA$ before.
The remaining ones will be defined in this section.
\medskip
\begin{center}
\begin{tabular}{|c|c|c|} \hline 
Geometric objects & \multicolumn2{c|}{Algebraic objects} \\ \hline \hline
$X$ & $PA$ & $QA$ \\ 
$\mathfrak X(X)$ & $\Der(PA)$ & $\Der(QA)$ \\ 
$\Pi TX$ & $\Omega(PA)$ & $\Omega(QA)$ \\ 
$\mathfrak X(\Pi TX)$ & $\Der(\Omega(PA))$ & $\Der(\Omega(QA))$ \\ \hline
\end{tabular}
\end{center}

\medskip
\begin{intuition}   
The geometric objects in the first column 
do not exist except in the commutative case $(P=c)$
in which case $X=\R^n$.
One may say that they are defined 
via their algebraic counterparts.
It is useful to have the following dictionary:

$X=$ operad manifold,
$\mathfrak X(X)=$ vector fields on $X$,
$\Pi TX=$ total space of the odd tangent bundle to $X$,
$\mathfrak X(\Pi TX)=$ vector fields on $\Pi TX$.
\end{intuition}

Let $X$ be an operad manifold with
the free $P$-algebra $PA$ as its algebra of functions
(\ref{ss:om}).
Starting with this purely algebraic data,
we want to construct an algebra of
``differential forms on X'', 
which is defined as some differential $P$-superalgebra
$\Omega(PA)$.
One would guess that the symplectic form $\omega$ should be 
a suitable element of $\Omega^2(PA)$.
However, this is not true.

It turns out that one needs to consider 
$\Omega(QA)$, which is the corresponding object for $QA$.
And the symplectic form $\omega$ will lie in $\Omega^2(QA)$.
So in this sense,
it is $QA$ and $\Omega(QA)$, which should be regarded
as functions and forms on $X$.

\subsection{Algebraic objects for an operad}  \label{ss:ao}

We first look more generally at operad manifolds
(\ref{ss:om}).
To every operad $P$, one can associate a category of $P$-algebras.
Similarly, one can consider the graded version,
which is the category of $P$-superalgebras.
And one can talk of derivations of $P$-algebras
or superderivations of $P$-superalgebras~\cite{\gk}.

In what follows, we will freely use this language.
For example,
$PA$ is the free object in the category
of $P$-algebras.
However, just as we did for $PA$,
we will also give explicit descriptions
of all the objects we deal with.
This should help the reader 
who is unfamiliar with the above language.
We also recommend that the reader keep
a concrete case, say commutative or associative,
in mind.

\begin{definition}   
Let $\Der(PA)$ be the space of derivations of $PA$.
It is a Lie algebra with the bracket being the commutator.
\[
\begin{minipage}{.8 in}
\begin{center}
\psfrag{x1}{\Huge $x_1$}
\psfrag{x2}{\Huge $x_2$}
\resizebox{.7 in}{.7 in}{\includegraphics{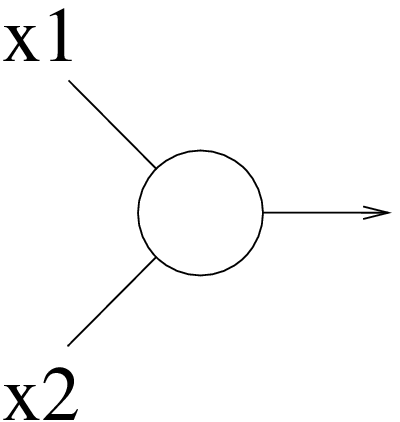}}
\end{center}
\end{minipage}
\qquad \buildrel {\xi}\over \longmapsto \qquad
\begin{minipage}{.8 in}
\begin{center}
\psfrag{x1}{\Huge $\xi(x_1)$}
\psfrag{x2}{\Huge $x_2$}
\resizebox{.7 in}{.7 in}{\includegraphics{pad.eps}}
\end{center}
\end{minipage}
\ \ + \ \
\begin{minipage}{.8 in}
\begin{center}
\psfrag{x1}{\Huge $x_1$}
\psfrag{x2}{\Huge $\xi(x_2)$}
\resizebox{.7 in}{.7 in}{\includegraphics{pad.eps}}
\end{center}
\end{minipage}
.
\]
\end{definition}
\noindent
A derivation ${\xi} \in \Der(PA)$ is
an operator on $PA$.
It is uniquely determined by its values 
${\xi}(x_1),\ldots,{\xi}(x_n) \in PA$ on the generators
$x_1,\ldots,x_n$.
An example of how this works is shown above;
for each input, we replace
the value $x_i$ on it by $\xi(x_i)$.
In the first term on the right,
$\xi(x_1)$ is written on one input.
This is to be understood as a substitution of $\xi(x_1)$
into that input.

\begin{definition}   
The algebra of forms $\Omega(PA)$ 
is the free differential superalgebra 
on the generators $x_1, \ldots,x_n$.
The term ``differential'' means that
there is a superderivation $d:\Omega(PA) \to \Omega(PA)$ 
of odd degree such that $d^2=0$.
\end{definition}
%\noindent
The algebra of forms $\Omega(PA)$ is freely generated 
as a $P$-superalgebra by the symbols 
$x_1, \ldots,x_n, dx_1, \ldots,dx_n$.
Hence a superderivation on $\Omega(PA)$ is uniquely determined
by its values on these generators.
Note that the differential $d$ is specified by
$d(x_i)=dx_i$ and $d(dx_i)=0$.

Just as for $PA$,
one can give an explicit description of $\Omega(PA)$.
To get a monomial in $\Omega(PA)$,
take any element of $P[I]$,
for some finite non-empty set $I$
and replace each element of $I$ by one of the generators
$x_1, \ldots,x_n,dx_1, \ldots,dx_n$.
Since we are now dealing with a super-object,
we need to be more careful.
We order the inputs with differential symbols 
in the sense of orientation as follows.
\begin{equation}  \label{e:dpa}
\begin{minipage}{.9 in}
\begin{center}
\psfrag{x1}{\Huge $x_1$}
\psfrag{x2}{\Huge $x_2$}
\psfrag{x3}{\Huge $x_3$}
\psfrag{dx2}{\Huge $dx_2$}
\psfrag{dx1}{\Huge $dx_1$}
\psfrag{1}{\Huge $1$}
\psfrag{2}{}
\psfrag{3}{\Huge $2$}
\psfrag{4}{}
\resizebox{.9 in}{.7 in}{\includegraphics{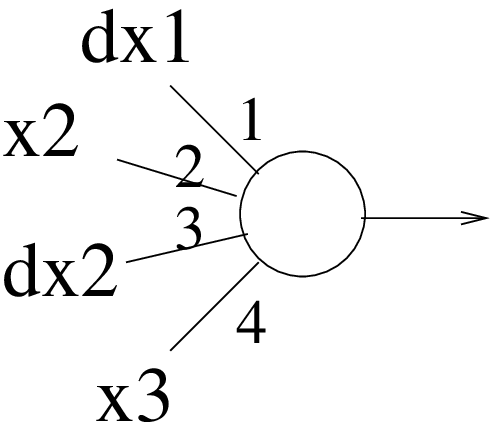}}
\end{center}
\end{minipage}
\qquad = \qquad - \ \ \ \ 
\begin{minipage}{.9 in}
\begin{center}
\psfrag{x1}{\Huge $x_1$}
\psfrag{x2}{\Huge $x_2$}
\psfrag{x3}{\Huge $x_3$}
\psfrag{dx2}{\Huge $dx_2$}
\psfrag{dx1}{\Huge $dx_1$}
\psfrag{1}{\Huge $2$}
\psfrag{2}{}
\psfrag{3}{\Huge $1$}
\psfrag{4}{}
\resizebox{.9 in}{.7 in}{\includegraphics{dpa.eps}}
\end{center}
\end{minipage}
.
\end{equation}
In other words, 
an even permutation of the order leaves an element unchanged
while an odd permutation gives its negative.
We now explain how the differential $d$ works on this example.
\[
\begin{minipage}{.9 in}
\begin{center}
\psfrag{x1}{\Huge $x_1$}
\psfrag{x2}{\Huge $x_2$}
\psfrag{x3}{\Huge $x_3$}
\psfrag{dx2}{\Huge $dx_2$}
\psfrag{dx1}{\Huge $dx_1$}
\psfrag{1}{\Huge $1$}
\psfrag{2}{}
\psfrag{3}{\Huge $2$}
\psfrag{4}{}
\resizebox{.9 in}{.7 in}{\includegraphics{dpa.eps}}
\end{center}
\end{minipage}
\qquad \buildrel {d}\over \longmapsto \qquad
\begin{minipage}{.9 in}
\begin{center}
\psfrag{x1}{\Huge $x_1$}
\psfrag{x2}{\Huge $dx_2$}
\psfrag{x3}{\Huge $x_3$}
\psfrag{dx2}{\Huge $dx_2$}
\psfrag{dx1}{\Huge $dx_1$}
\psfrag{1}{\Huge $2$}
\psfrag{2}{\Huge $1$}
\psfrag{3}{\Huge $3$}
\psfrag{4}{}
\resizebox{.9 in}{.7 in}{\includegraphics{dpa.eps}}
\end{center}
\end{minipage}
\ \ + \ \
\begin{minipage}{.9 in}
\begin{center}
\psfrag{x1}{\Huge $x_1$}
\psfrag{x2}{\Huge $x_2$}
\psfrag{x3}{\Huge $dx_3$}
\psfrag{dx2}{\Huge $dx_2$}
\psfrag{dx1}{\Huge $dx_1$}
\psfrag{1}{\Huge $2$}
\psfrag{2}{}
\psfrag{3}{\Huge $3$}
\psfrag{4}{\Huge $1$}
\resizebox{.9 in}{.7 in}{\includegraphics{dpa.eps}}
\end{center}
\end{minipage}
.
\]
In general, we replace each $x_i$ by a $dx_i$,
one at a time.
And the input at which the replacement occurs
is given label $1$ and the remaining labels
are shifted up by $1$.

It is clear that $\Omega(PA)$ is in fact $\Z$ graded
with the grading given by the number of differential symbols
and $\Omega^0(PA)=PA$.

\begin{definition}   
Let $\Der(\Omega(PA))$ be the space 
of superderivations of $\Omega(PA)$.
It is a Lie superalgebra with the bracket being the supercommutator.
\end{definition}

To act by a superderivation $L$
on an element of $\Omega(PA)$,
one replaces the value $v$ at each input by $L(v)$,
one at a time.
If $v=x_i$ then we shift the labels by $\degg(L)$.
If $v=dx_i$ then we first reorder such that the $dx_i$ 
being substituted into has label $1$ and then 
shift the remaining labels by $\degg(L)$.
The special case when $L$ is the differential $d$ was shown above.

\begin{definition}   
For any derivation ${\xi} \in \Der(PA)$,
we define the Lie derivative $L_{\xi} \in \Der(\Omega(PA))$ 
and the contraction operator $i_{\xi} \in \Der(\Omega(PA))$ 
by specifying them on the generators as under
$$L_{\xi}(x_i)={\xi}(x_i), \ L_{\xi}(dx_i)=d{\xi}(x_i) \ \ \text{and} \ \
i_{\xi}(x_i)=0, \ i_{\xi}(dx_i)={\xi}(x_i).$$
They have degrees $0$ and $-1$ respectively.
\end{definition} 

\begin{lemma} \label{l:rel}
For derivations ${\xi}, {\eta} \in \Der(PA)$,
the superderivations $L_{\xi}, i_{\xi}, d$ satisfy
the following relations.
\begin{itemize}
\item[(1)]
$[i_{\xi},d] = L_{\xi}$ \ \ \ (Cartan's formula). 

\item[(2)]
$[i_{\xi},i_{\eta}] = 0$.

\item[(3)]
$[L_{\xi},i_{\eta}] = i_{[{\xi},{\eta}]}$.

\item[(4)]
$[L_{\xi},L_{\eta}] =  L_{[{\xi},{\eta}]}$.
\end{itemize}

\end{lemma}
\begin{proof}
As mentioned earlier, superderivations are preserved
by the supercommutator.
Hence in all cases, both sides are superderivations,
so we need to check $(1)-(4)$ only on the generators
$x_1, \ldots,x_n, dx_1, \ldots,dx_n$.
This is a straightforward check that is
independent of the operad $P$.
\end{proof}

\begin{lemma} \label{l:exact}
The algebra of forms $\Omega(PA)$ is exact.
\end{lemma}
\begin{proof}
Consider the Euler vector field $e \in \Der(PA)$,
defined by $e(x_i)=x_i$ for all $i$.
Then $L_e(x_i)=x_i$ and $L_e(dx_i)=dx_i$
and hence $L_e$ maps every monomial
to a non-zero multiple of itself.
Therefore it is invertible on $\Omega(PA)$.
So it induces an isomorphism on cohomology.
On the other hand, Cartan's formula
(item $(1)$ in Lemma~\ref{l:rel}) shows that
it induces the zero map on cohomology.
So we conclude that $\Omega(PA)$ is exact.
\end{proof}

For the operad $P=c$,
the lemma just says that the algebra of polynomial forms on $\R^n$,
that vanish at the origin, is exact.

\subsection{Algebraic objects for a species}  \label{ss:as}

Now assume that the operad $P$ is reversible;
so we can talk of its mated species $Q$.
Analogues for $Q$ of the objects in (\ref{ss:ao})
can be defined without difficulty.
As far as pictures are concerned,
we draw the same ones as for operads
except that the output arrow is now omitted.

We have seen that $QA$ can be written
as a quotient of $PA \otimes PA$ (equation~\eqref{e:r}).
Similarly, one can define $\Omega(QA)$
as a quotient of $\Omega(PA) \otimes \Omega(PA)$ 
subject to the super-relations:
\begin{equation}  \label{e:sr}
a \otimes b = (-1)^{|a||b|} b \otimes a
\quad \text{and} \quad
a \rightarrow b \otimes c = a \otimes r(b) \leftarrow c.
\end{equation}
Just like $\Omega(PA)$,
we see that $\Omega(QA)$ is $\Z$ graded
and $\Omega^0(QA)=QA$.
We also define $\Der(QA)$, $\Der(\Omega(QA))$ 
and for $\xi \in \Der(PA)$,
the operators
$L_{\xi}, i_{\xi} \in \Der(\Omega(QA))$.
The relations in Lemma~\ref{l:rel} continue to hold.
Hence it follows that $\Omega(QA)$ 
is also exact.

Note that any element of $\Omega^1(QA)$
can be uniquely written as
$\sum_{i=1}^{n} f_i \otimes dx_i$, with
$f_i \in PA$.
And for $H \in QA$,
we get $dH = \sum_{i=1}^{n} \frac{\partial H}{\partial x_i} \otimes dx_i$.
These facts are again a consequence of Proposition~\ref{p:pd}
and will be crucial in what follows.

\subsection{The symplectic form $\omega$}  

Since $P$ has a unit element $u$,
its mated species $Q$ has the special element $uu$
in degree 2, see (\ref{ss:uu}).
It allows us to define the symplectic form
$\omega = \sum dp_i \otimes dq_i \in \Omega^2(QA)$.
For a picture, we may draw 
$\sum
\begin{minipage}{.8 in}
\begin{center}
\psfrag{p1}{\Huge $dp_i$}
\psfrag{q2}{\Huge $dq_i$}
\psfrag{1}{\Huge $1$}
\psfrag{2}{\Huge $2$}
\psfrag{uu}{}
\resizebox{.7 in}{.15 in}{\includegraphics{sf.eps}}
\end{center}
\end{minipage}
$.
Note that if we switch the order of the factors,
we pick a minus sign (equation~\eqref{e:sr}).
This can also be seen from the picture and
the analogue of equation~\eqref{e:dpa} for species.

Also $dw=0$ and hence the form is closed.
This gives us a symplectic operad manifold $X$.
We think of $\Der(PA)$ as vector fields on $X$
that vanish at a point and 
$\Omega(QA)$ as differential forms on $X$ 
with no constant or linear terms,
with the $0$-forms $QA$ being Hamiltonian functions on $X$.

\begin{lemma} \label{l:iso}
The map $\Der(PA) \to \Omega^1(QA)$ which sends
${\xi}$ to $i_{\xi} w$ is an isomorphism.
\end{lemma}
\begin{proof}
We have the chain of equalities
$$i_{\xi} w = \sum_{i=1}^{n} i_{\xi} (dp_i \otimes dq_i) =
\sum_{i=1}^{n} {\xi}(p_i) \otimes dq_i - {\xi}(q_i) \otimes dp_i.$$
This shows that $i_{\xi} w$ is determined by 
${\xi}(p_i)$ and ${\xi}(q_i)$, which also determine ${\xi}$.
Since one can represent any element of $\Omega^1(QA)$ 
uniquely by a sum of the form
$\sum_{i=1}^{n} ( \ ) \otimes dp_i + ( \ ) \otimes dq_i$,
the map in the lemma is an isomorphism.
\end{proof}

\begin{definition}   
Let $\Der(PA, w) = \{ {\xi} \in \Der(PA) \ | \ L_{\xi} w = 0 \}$.
In analogy with the classical case, 
we call these the symplectic vector fields.
\end{definition}

It follows from Cartan's formula (Lemma~\ref{l:rel})
that 
$${\xi} \in \Der(PA, w) \Longleftrightarrow d(i_{\xi} w) = 0.$$
Hence under the isomorphism of Lemma~\ref{l:iso},
one sees that symplectic vector fields correspond 
to closed 1 forms.
Define Hamiltonian vector fields to be the ones that correspond 
to exact 1 forms.
However, since $\Omega(QA)$ is exact, 
the Hamiltonian and symplectic vector fields coincide
in this case.
And there is an isomorphism
$QA \buildrel d \over \to \Omega^1(QA)_{\text{closed}}$.
Putting all facts together, 
there is an isomorphism 
$\xi : QA \isoto \Der(PA, w)$, 
where $H \mapsto {\xi}_H$
is defined by the equation 
$dH = i_{{\xi_H}}(w)$.
Writing both sides in the unique form
$\sum_{i=1}^{n} ( \ ) \otimes dp_i + ( \ ) \otimes dq_i$,
we derive the Hamiltonian equations
\begin{equation}  \label{e:ham}
{\xi_H}(p_i) = \frac{\partial H}{\partial q_i} \ \ \text{and} \ \ 
{\xi_H}(q_i) = - \frac{\partial H}{\partial p_i}.
\end{equation}
We now derive the isomorphism between
Hamiltonian functions and Hamiltonian vector fields.
By Kontsevich kimaya,
the proof is same as in the classical case;
we give it here for completeness.

\begin{proposition}  
The map $QA \to \Der(PA, w)$ defined above
is an anti-isomorphism of Lie algebras,
with the bracket on $QA$ being the Poisson bracket
and on $\Der(PA, w)$ 
being the usual commutator.
\end{proposition}
\begin{proof}

First note that $\Der(PA, w)$ is closed under
taking commutators
(item $(4)$ in Lemma~\ref{l:rel}).
This gives it a Lie algebra structure.
So to prove the proposition,
we trace through the isomorphism 
$QA \isoto \Der(PA, w)$,
compute the induced Lie structure on $QA$
and see that upto a minus sign 
it coincides with the Poisson bracket
(\ref{ss:mate}).
We do the computation in two steps.

\subsection*{Claim $1$}
$i_{[{\xi_F},{\xi_H}]}w = d(i_{\xi_F}(dH))$.
$$
\begin{array}{r c l l}
i_{[{\xi_F},{\xi_H}]}w & = & [L_{\xi_F},i_{\xi_H}]w & \mbox{(item (3) in Lemma~\ref{l:rel})} \\
             & = & L_{\xi_F}(i_{\xi_H} w) & (L_{\xi_F} w=0 \ \mbox{since} \ \xi_F \ \mbox{is symplectic}) \\
             & = & L_{\xi_F}(dH) & \\
             & = & [i_{\xi_F},d](dH) & \mbox{(item (1) in Lemma~\ref{l:rel})} \\
             & = & d(i_{\xi_H} (dH)). &  \\
\end{array}
$$

\subsection*{Claim $2$}
$i_{\xi_F}(dH) = - \{F,H\}$.\\
$$
\begin{array}{r c l l}
i_{\xi_F}(dH) & = & i_{\xi_F}(\sum_{i=1}^{n} \frac{\partial H}{\partial p_i} \otimes dp_i+ \frac{\partial H}{\partial q_i} \otimes dq_i) & \\
             & = & \sum_{i=1}^{n} \frac{\partial H}{\partial p_i} \otimes {\xi_F}(p_i) + \frac{\partial H}{\partial q_i} \otimes {\xi_F}(q_i) & \\
             & = & \sum_{i=1}^n \frac{\partial H}{\partial p_i} \otimes \frac{\partial F}{\partial q_i} - \frac{\partial H}{\partial q_i} \otimes \frac{\partial F}{\partial p_i} & \mbox{(equation~\eqref{e:ham})} \\
             & = & - \{F,H\}. & \mbox{(equation~\eqref{e:ptb})} \\
                    
\end{array}
$$
The proposition now follows from the two claims.
\end{proof}

As a corollary, we obtain Proposition~\ref{p:pb}.
Also, if we specialise to $P=u$ and $Q=uu$,
we obtain
Proposition~\ref{p:degtwo}.

\section{Examples motivated by PROPS}  \label{s:me}

We have, so far, only defined the notion of a derivative.
If one wants to define higher order derivatives
then one should consider a more general structure,
namely, a PROP.
For some motivation for PROPs
from the viewpoint of homotopy theory, see~\cite{\adams}. 

\subsection{PROP}  \label{ss:prop}

We denote a PROP by $Pr$.
It consists roughly of objects 
with many inputs and many outputs.
For a precise definition, 
see the notes by Voronov~\cite{\voronov}. 
We write $Pr[I,J]$ for the ``set of $Pr$-structures 
on the set $I \sqcup J$ that have $I$ inputs and $J$ outputs''.
There is also a composition law
which allows to feed the outputs of one object
to the inputs of another.
\[
\begin{minipage}{2 in}
\begin{center}
\resizebox{1.8 in}{.6 in}{\includegraphics{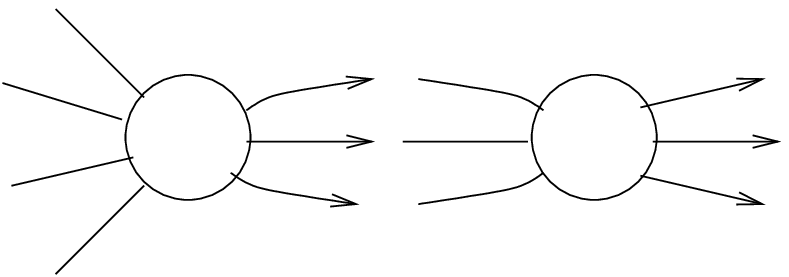}}
\end{center}
\end{minipage}
=
\begin{minipage}{1 in}
\begin{center}
\resizebox{.85 in}{.6 in}{\includegraphics{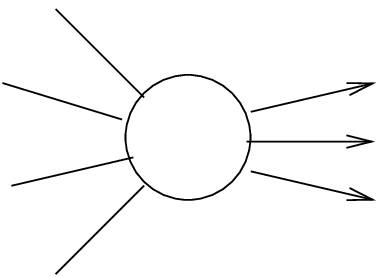}}
\end{center}
\end{minipage}.
\]

It is useful to think of species and operads
as parts of a PROP corresponding to 
no output ($J$ empty) and
a single output ($J$ singleton) respectively.
We saw that a partial derivative (Proposition~\ref{p:pd})
sends a species element to an operad element.
Stated differently, it just converts 
an input to an output.
If this is to be the meaning of a derivative
then to define higher derivatives
one is forced to consider objects with many outputs.
Thus one sees that PROPs provide a natural framework
to talk of higher derivatives.
To make the theory of higher matings work,
one would need a ``reversible PROP''.
We do not make these ideas precise.
Instead, we look at some examples
that fit this pattern.
They may be of independent interest.

\subsection{The graph operad $g$ and species $gg$} \label{ss:g}

Let $g[I]=$ set of graphs with a specified source and sink vertex
and directed edges labelled by elements of the set $I$.
The figure shows an element of $g[\{a,b,c,d\}]$.
\[
\begin{minipage}{1.6 in}
\begin{center}
\psfrag{a}{\Huge $a$}
\psfrag{b}{\Huge $b$}
\psfrag{c}{\Huge $c$}
\psfrag{d}{\Huge $d$}
\psfrag{so}{\Huge $\text{source}$}
\psfrag{si}{\Huge $\text{sink}$}
\resizebox{1 in}{.6 in}{\includegraphics{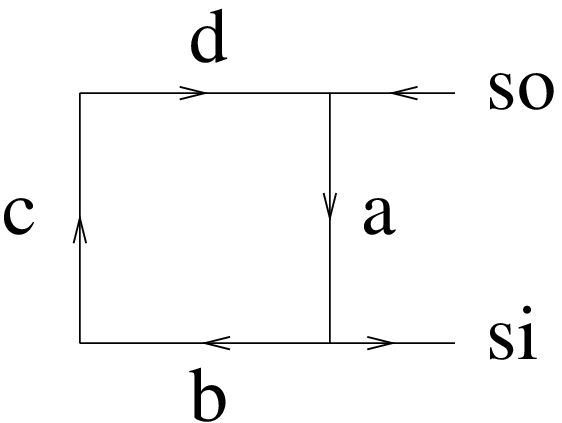}}
\end{center}
\end{minipage}
\]
We indicate the source (resp. sink) vertex 
by a half-edge with an arrow going in (resp. out).
We regard the two half-edges together as 
constituting the output edge of the operad.
Substitution works as under.
\[
\begin{minipage}[b]{.9 in}
\begin{center}
\psfrag{a}{\Huge $a$}
\psfrag{b}{\Huge $b$}
\psfrag{c}{\Huge $c$}
\psfrag{d}{\Huge $d$}
\psfrag{so}{\Huge $$}
\psfrag{si}{\Huge $$}
\resizebox{.9 in}{.7 in}{\includegraphics{gsubt.eps}}\\
$p_2$
\end{center}
%\label{f:gsubt}
\end{minipage}
\hfill
\begin{minipage}[b]{1.2 in}
\begin{center}
\psfrag{x}{\Huge $x$}
\psfrag{y}{\Huge $y$}
\psfrag{z}{\Huge $z$}
\psfrag{w}{\Huge $$}
\psfrag{so}{\Huge $\text{source}$}
\psfrag{si}{\Huge $\text{sink}$}
\resizebox{.9 in}{.5 in}{\includegraphics{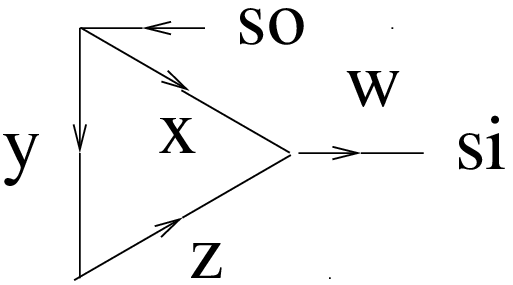}}\\
\vspace{.1 in}
$p_1$
\end{center}
%\label{f:gsubi}
\end{minipage}
\hspace{1 in}
\begin{minipage}[b]{1.4 in}
\begin{center}
\psfrag{a}{\Huge $a$}
\psfrag{b}{\Huge $b$}
\psfrag{c}{\Huge $c$}
\psfrag{d}{\Huge $d$}
\psfrag{x}{\Huge $x$}
\psfrag{z}{\Huge $z$}
\psfrag{w}{\Huge $$}
\psfrag{so}{\Huge $\text{source}$}
\psfrag{si}{\Huge $\text{sink}$}
\resizebox{1.3 in}{.7 in}{\includegraphics{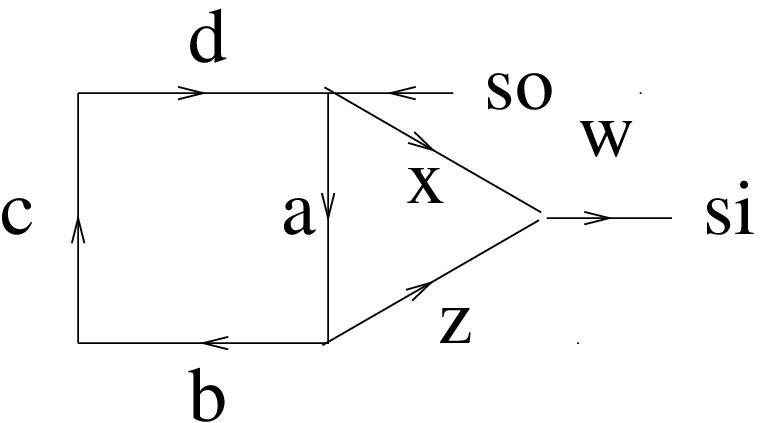}}\\
$p_2 \yright p_1$
\end{center}
%\label{f:gsubo}
\end{minipage}.
\]
In other words,
to obtain
$p_2 \yright p_1$,
we identify the source (resp. sink) vertex of $p_2$
with the initial (resp. terminal) vertex of the edge $y$ in $p_1$.

\subsubsection{Sub-examples}
The commutative, associative, chord 
(\ref {ss:eg}) and many other operads
can be obtained as special cases
by simply restricting the type of graphs allowed.
There is seemingly no end to the number of interesting examples
one can obtain in this way.
We illustrate the commutative and associative cases.
\[
\begin{minipage}[b]{1.8 in}
\begin{center}
\vspace{.1 in} 
\psfrag{a}{\Huge $a$}
\psfrag{b}{\Huge $b$}
\psfrag{c}{\Huge $c$}
\psfrag{d}{\Huge $d$}
\psfrag{x}{\Huge $x$}
\psfrag{z}{\Huge $z$}
\psfrag{w}{\Huge $$}
\psfrag{so}{\Huge $\text{source}$}
\psfrag{si}{\Huge $\text{sink}$}
\resizebox{1.8 in}{.4 in}{\includegraphics{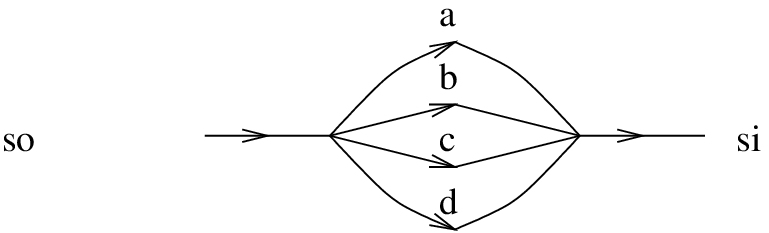}}\\
commutative
\end{center}
\end{minipage}
\hspace{1 in}
\begin{minipage}[b]{1.6 in}
\begin{center}
\vspace{.1 in} 
\psfrag{a}{\Huge $a$}
\psfrag{b}{\Huge $b$}
\psfrag{c}{\Huge $c$}
\psfrag{d}{\Huge $d$}
\psfrag{x}{\Huge $x$}
\psfrag{z}{\Huge $z$}
\psfrag{w}{\Huge $$}
\psfrag{so}{\Huge $\text{source}$}
\psfrag{si}{\Huge $\text{sink}$}
\resizebox{1.6 in}{.4 in}{\includegraphics{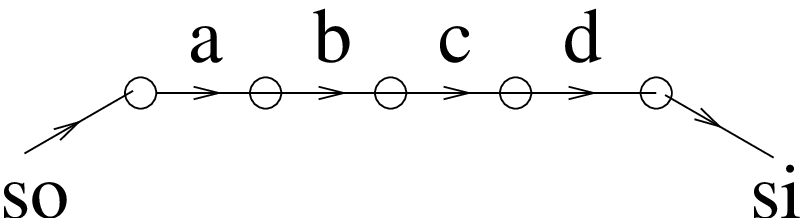}}\\
%\vspace{.1 in} 
associative
\end{center}
\end{minipage}
\]
Thus to get the commutative operad,
we only allow graphs with two vertices (source and sink)
and all edges are directed from the source to the sink.
For the associative operad,
we allow graphs whose edges form a directed path 
from the source to the sink.

There is a variation on the associative case
(and the commutative),
where one allows the arrows on the edges to point
in either direction.
This is a binary quadratic operad.
Axioms for an algebra over this operad can be written explicitly.
We do not know whether they have been considered before.

\subsubsection{The reversal maps}
The graph operad is reversible 
with the reversal map as shown.
\[
r_{x,y}
\left(
\begin{minipage}{1.4 in}
\begin{center}
\psfrag{x}{\Huge $a$}
\psfrag{y}{\Huge $x$}
\psfrag{z}{\Huge $b$}
\psfrag{w}{\Huge $y$}
\psfrag{so}{\Huge $\text{source}$}
\psfrag{si}{\Huge $\text{sink}$}
\resizebox{1 in}{.5 in}{\includegraphics{gsubi.eps}}
\end{center}
\end{minipage}
\right)
\ \ = \ \
\begin{minipage}{1.4 in}
\begin{center}
\psfrag{x}{\Huge $a$}
\psfrag{y}{\Huge $x$}
\psfrag{z}{\Huge $b$}
\psfrag{w}{\Huge $y$}
\psfrag{so}{\Huge $\text{source}$}
\psfrag{si}{\Huge $\text{sink}$}
\resizebox{1.2 in}{.5 in}{\includegraphics{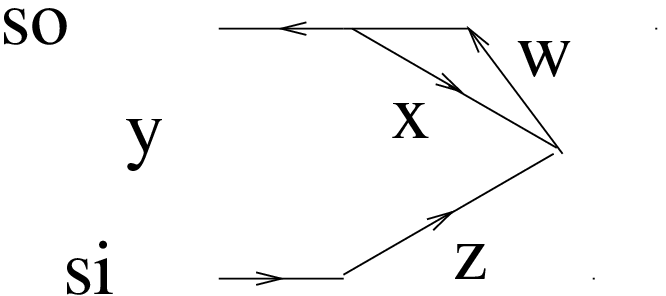}}
\end{center}
\end{minipage}
.
\]
Namely,
we open the input edge $x$ and close the output edge $y$.

There is another way to reverse the graph operad,
where in addition to the above,
we also switch the directions on $x$ and $y$.
The two methods of reversal are compatible
with the associative and commutative cases respectively.
Since the two ways are analogous, we will work with the first definition.

\subsubsection{How graph objects mate}
So far, we have only defined the graph operad.
But the definition of the graph PROP should be clear;
we simply allow more output edges.
In other words, graphs can now have many pairs
of source and sink vertices.
And for the composition law in the PROP,
we perform simultaneous substitutions.
The graph species $gg$ is the part of the graph PROP
with no outputs.
So these are graphs with no source and sink vertex.
We now show how graph objects mate.
\[
\begin{minipage}{1.5 in}
\begin{center}
\psfrag{a}{\Huge $a$}
\psfrag{b}{\Huge $b$}
\psfrag{c}{\Huge $c$}
\psfrag{d}{\Huge $d$}
\psfrag{x}{\Huge $x$}
\psfrag{y}{\Huge $y$}
\psfrag{z}{\Huge $z$}
\resizebox{1.4 in}{.7 in}{\includegraphics{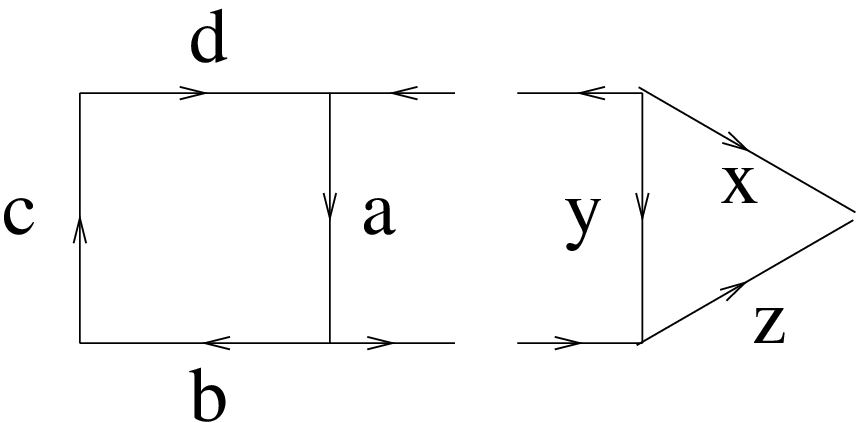}}
\end{center}
\end{minipage}
=
\begin{minipage}{1.5 in}
\begin{center}
\psfrag{a}{\Huge $a$}
\psfrag{b}{\Huge $b$}
\psfrag{c}{\Huge $c$}
\psfrag{d}{\Huge $d$}
\psfrag{x}{\Huge $x$}
\psfrag{y}{\Huge $y$}
\psfrag{z}{\Huge $z$}
\resizebox{1.3 in}{.7 in}{\includegraphics{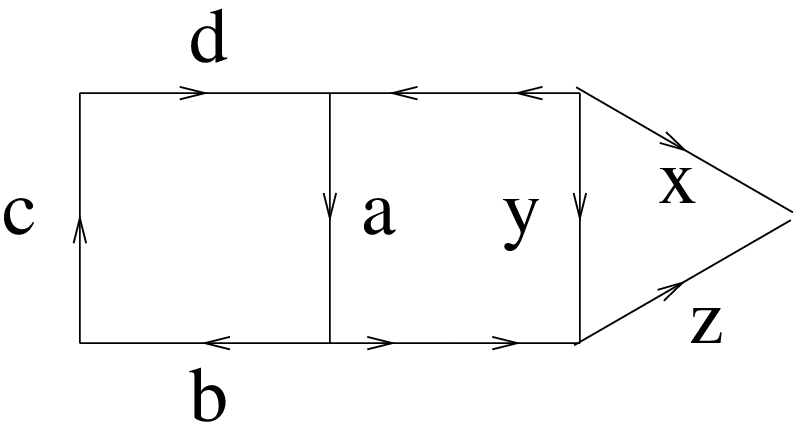}}
\end{center}
\end{minipage}
=
\begin{minipage}{1.1 in}
\begin{center}
\psfrag{a}{\Huge $a$}
\psfrag{b}{\Huge $b$}
\psfrag{c}{\Huge $c$}
\psfrag{d}{\Huge $d$}
\psfrag{x}{\Huge $x$}
\psfrag{y}{\Huge $y$}
\psfrag{z}{\Huge $z$}
\resizebox{.9 in}{.7 in}{\includegraphics{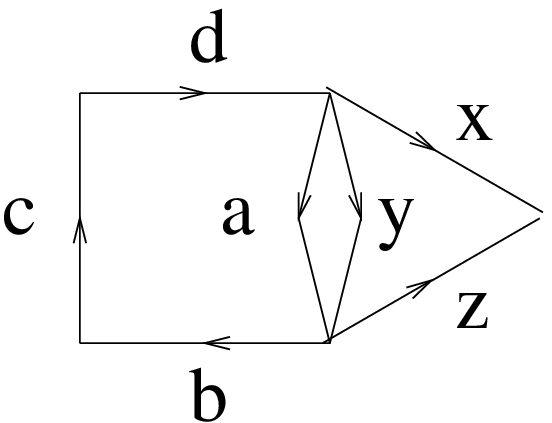}}
\end{center}
\end{minipage}.
\]
Namely, the sink of one merges (mates) with the source of the other
and vice-versa.
The result is a graph with no source and sink vertex.
Thus the mating functor maps the graph operad $g$
to the graph species $gg$.

\subsubsection{The algebraic objects}
The free algebras $gA$ and $ggA$ consist of graphs
with directed edges labelled by one of the generators
$x_1, \ldots,x_n$.
Of course, for elements of $gA$,
graphs also have a specified source and sink vertex.
The partial derivative 
$\frac{\partial }{\partial x_i} : ggA \longrightarrow gA$
works very nicely.
\[
\begin{minipage}{1 in}
\begin{center}
\psfrag{c}{\Huge $x_1$}
\psfrag{z}{\Huge $x_1$}
\psfrag{b}{\Huge $x_2$}
\psfrag{a}{\Huge $x_2$}
\psfrag{d}{\Huge $x_4$}
\psfrag{y}{\Huge $x_4$}
\psfrag{x}{\Huge $x_3$}
\resizebox{.8 in}{.7 in}{\includegraphics{gmd.eps}}
\end{center}
\end{minipage}
\qquad \buildrel \frac{\partial}{\partial x_1}\over \longmapsto \qquad
\begin{minipage}{1 in}
\begin{center}
\psfrag{c}{\Huge $x_1$}
\psfrag{z}{\Huge $x_1$}
\psfrag{b}{\Huge $x_2$}
\psfrag{a}{\Huge $x_2$}
\psfrag{d}{\Huge $x_4$}
\psfrag{y}{\Huge $x_4$}
\psfrag{x}{\Huge $x_3$}
\resizebox{.9 in}{.7 in}{\includegraphics{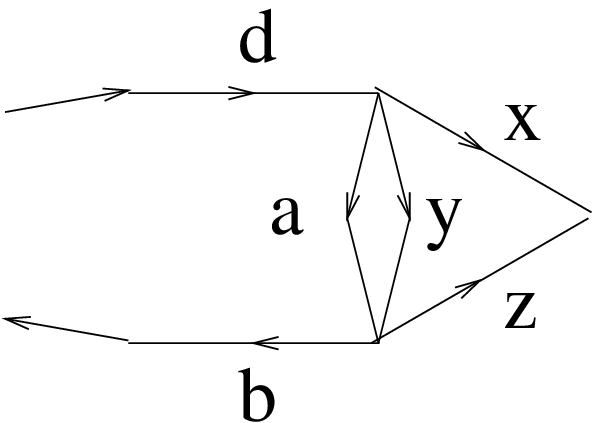}}
\end{center}
\end{minipage}
\ \ + \ \
\begin{minipage}{1 in}
\begin{center}
\psfrag{c}{\Huge $x_1$}
\psfrag{z}{\Huge $x_1$}
\psfrag{b}{\Huge $x_2$}
\psfrag{a}{\Huge $x_2$}
\psfrag{d}{\Huge $x_4$}
\psfrag{y}{\Huge $x_4$}
\psfrag{x}{\Huge $x_3$}
\resizebox{.8 in}{.7 in}{\includegraphics{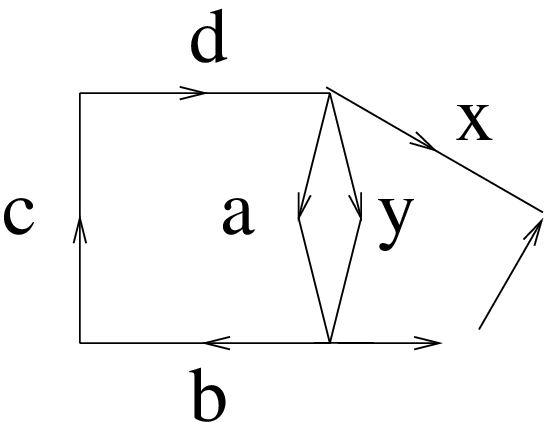}}
\end{center}
\end{minipage}
.
\]
In words,
$\frac{\partial F}{\partial x_i}$ 
is the sum of graphs obtained by cutting edges of $F$
that are labelled $x_i$, one at a time.
The half-edges that result from a cut get induced directions
and provide the source and sink vertex
of the graph created.
We leave it to the reader to explicitly describe
the algebra of forms $\Omega(gA)$, $\Omega(ggA)$
and other notions described in Section~\ref{s:st}.

For higher order derivatives,
we consider the free algebra of the graph PROP.
The order of the derivative determines 
the number of edges that are cut.
And for higher order matings,
we take two graphs with the same number of cuts (outputs)
and the result is a graph with no cuts.
This generalises the first order mating 
given by the mating functor.

\subsection{A generalisation of the tree example
(\ref {sss:t})} \label{ss:gt}

We now present a second example based on graphs.
It complements the earlier example in the sense
that (univalent) vertices rather than edges
play the role of inputs.
Let $hh[I]=$ set of graphs with univalent vertices
labelled by elements of the set $I$.
The operad $h$ is defined as the species $hh$
except that we use one of the univalent vertices
as the output or root.
And for the PROP,
we do not place any restriction
on the number of outputs or roots.
Operad substitution and reversal and matings work
exactly as in the tree case.
\[
\begin{minipage}[b]{1.8 in}
\begin{center}
\psfrag{a}{\Huge $a$}
\psfrag{b}{\Huge $b$}
\psfrag{y}{\Huge $y$}
\psfrag{z}{\Huge $z$}
\resizebox{1.7 in}{.5 in}{\includegraphics{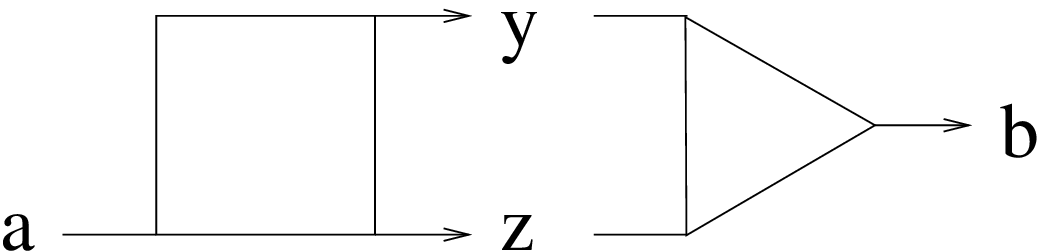}}
\end{center}
\end{minipage}
\hspace{1 in}
\begin{minipage}[b]{1.6 in}
\begin{center}
\psfrag{a}{\Huge $a$}
\psfrag{b}{\Huge $b$}
\psfrag{y}{\Huge $y$}
\psfrag{z}{\Huge $z$}
\resizebox{1.6 in}{.5 in}{\includegraphics{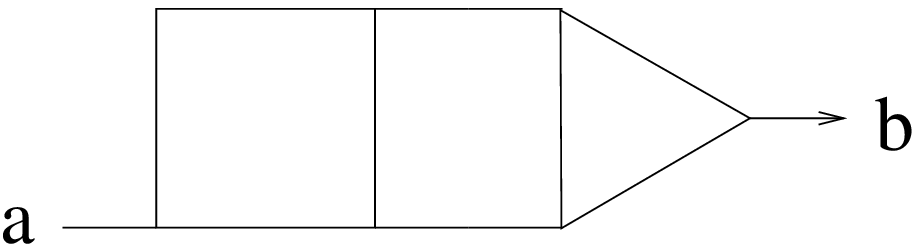}}
\end{center}
\end{minipage}.
\]
Note that for the above composition law
of the PROP to make sense,
we necessarily have to work with graphs and not just trees.
We mention that the graph complexes that arise from this species
have been considered in~\cite{\bm}.

\subsection{The surface operad $s$ and species $ss$} \label{ss:ss}

Let $ss[I]=$ set of compact orientable surfaces
(not necessarily connected) with boundary circles
labelled by elements of the set $I$.
For the operad $s[I]$, 
we use one of the boundary circles as the output.
Operad substitution works by gluing.
For the PROP, we allow more outputs as usual.
The composition law in the PROP is simultaneous gluing.
We leave the details to the reader.
Some pictures may be found in~\ref{ss:ae}.

\section{Graph homology} \label{s:gh}

In this section, we define a generalisation of 
Kontsevich's graph homology to mated species.
This coincides with the graph complexes
considered by Markl~\cite{\markl}.
When one specialises to the commutative species
(\ref{sss:c}),
one recovers the usual definition.
The case of the associative species (\ref{sss:a}) is 
also explained in~\cite{\voronov}.
By a graph, we mean a finite 1-dimensional CW complex.

In (\ref {ss:graph}-\ref {ss:gc}),
we define a chain complex $(\mathcal G,\partial_E)$
based on graphs.
Here $\mathcal G$ refers to the chain groups
and $\partial_E$ to the boundary operator.
The boundary operator $\partial_E$
can be regarded as a limit of a family of boundary operators
$\partial_n$.
This family is discussed in 
(\ref {ss:finb}).
In (\ref {ss:cc}),
we list the subcomplexes of $(\mathcal G,\partial_E)$
that are relevant to us.
The graph complex $(Q\mathcal G,\partial_E)$,
which appears in the main theorem, is one of them.

\subsection{Graphs}  \label{ss:graph}

Let $Q$ be any species.
Define a $Q$-graph to be a graph whose vertices are fattened 
by elements of the species $Q$.
To be more precise,
it is a graph $\Gamma$ such that
for every vertex $v$ of $\Gamma$,
a $Q$-structure is specified 
on the set of half-edges incident to $v$.
\[
\begin{minipage}{1.8 in}
\begin{center}
\psfrag{1}{\Huge $1$}
\psfrag{2}{\Huge $2$}
\psfrag{3}{\Huge $3$}
\psfrag{4}{\Huge $4$}
\resizebox{.9 in}{.9 in}{\includegraphics{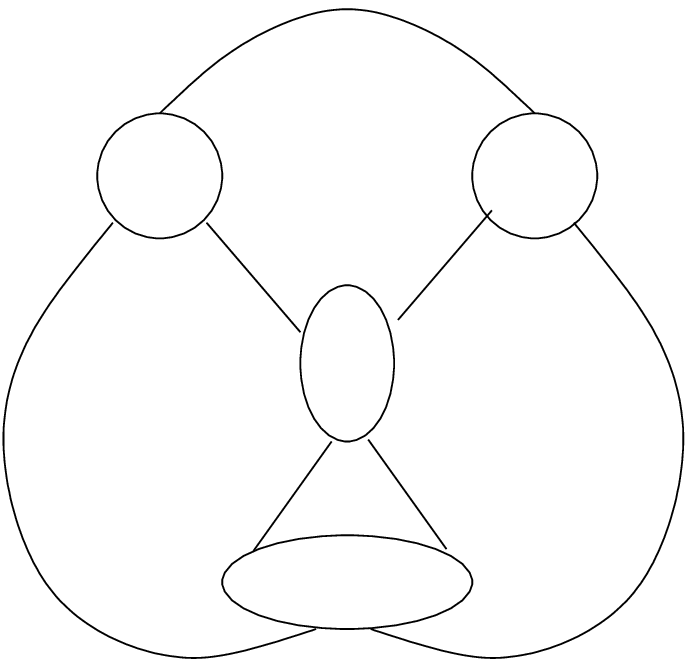}}\\
\vspace{.1 in} Greetings from a $Q$-graph.
\end{center}
\end{minipage}
\]

\noindent
We now give some examples;
see the table of pictures
and examples of species in (\ref{ss:s}) and (\ref{ss:eg}).
\begin{itemize}
\item
A commutative graph $(Q=cc)$  is simply a graph.

\item
An associative graph $(Q=aa)$ is a ribbon graph;
that is a cyclic ordering of
the half-edges is fixed at each vertex.
%And in the Lie graph $(Q=ll)$, a bracket expression of the half-edges is given.

\item
The figure below shows a typical graph for a species based on trees
$(Q=tt)$.
\[
\begin{minipage}{1.8 in}
\begin{center}
\psfrag{1}{\Huge $1$}
\psfrag{2}{\Huge $2$}
\psfrag{3}{\Huge $3$}
\psfrag{4}{\Huge $4$}
\resizebox{1.3 in}{1 in}{\includegraphics{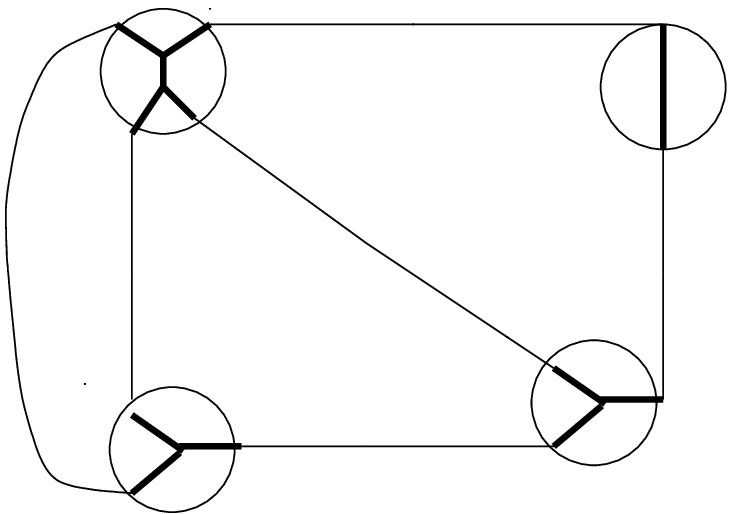}}
\end{center}
\end{minipage}
\]

\item
For the chord species $(Q=kk)$,
the vertices instead of being trees are chord diagrams
(\ref{sss:k}).
In our terminology, these are chord graphs.

\item
For an example with a somewhat different flavour,
consider the surface species $Q=ss$, see (\ref{ss:ss}).
A surface graph is a compact orientable surface
with a collection of disjoint loops.
In this case,
the loops play the role of edges and
the pieces obtained 
by cutting along the loops play the role of vertices.
This example was first mentioned to me by J. Conant.
\end{itemize}
The interested reader can also work out the graphs
for the species $hh$ and the graph species $gg$
in Section~\ref{s:me}.
Note that all the species above have been denoted
by a letter repetition.
This is to indicate the fact that they are
mated species (\ref{ss:mf}-\ref{ss:am}).

The reader, who skipped Sections~\ref{s:ro} and \ref{s:mf}
on species and operads, 
can read this section and Section~\ref{s:mgh}
by keeping the above concrete examples in mind.
For simplicity of notation, from now on,
we will just write ``graph'' instead of ``$Q$-graph''.

\subsection{Oriented graphs}  \label{ss:og}

For a graph $\Gamma$, we shall denote the set of vertices by
$V(\Gamma)$, the set of edges by $E(\Gamma)$
and the set of ends of
an edge $e$ by $V(e)$.
%, and the set of edges emanating from a
%vertex $v$ by $E(v)$.  
Finally, $\mathbb{R}X$ will
denote the real vector space which has the elements of $X$ as a
basis.

\begin{definition}   
An orientation $\sigma$ of a graph $\Gamma$ is an
orientation of the one-dimensional real vector space $\det
\mathbb RV(\Gamma) \otimes\bigotimes_{e\in E(\Gamma)}\det \mathbb
RV(e)$. 
For a vector space $W$ of dimension $n$,
we are using the notation $\det W = \Lambda^n W$. 
We say that $(\Gamma, \sigma)$ is an oriented graph.
There is another notion of orientation of a graph
equivalent to the above; see \cite{\thurston} for details.
\end{definition}

A way to represent an orientation is to order the
vertices and orient each edge of the graph.  An odd permutation
of the labels on the vertices reverses the orientation, and a
single change of the orientation of one edge reverses it as
well.  An even number of these transformations produces an
orientation equivalent to the original one.
So a graph $\Gamma$ has two orientations,
which we may call $\sigma$ and $- \sigma$.
\[
\begin{minipage}{\linewidth}
\begin{center}
\psfrag{1}{\Huge $1$}
\psfrag{2}{\Huge $2$}
\psfrag{3}{\Huge $3$}
\psfrag{4}{\Huge $4$}
\resizebox{1.3 in}{1 in}{\includegraphics{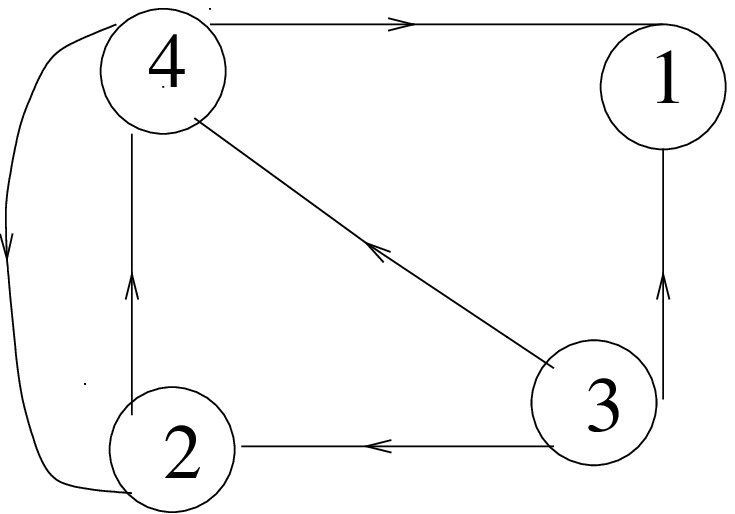}}\\
\vspace{.1 in} A representative of an orientation on a $Q$-graph.
\end{center}
\end{minipage}
\]

We say that $(\Gamma, \sigma )$ and $(\Gamma^\prime,
\sigma^\prime )$ are isomorphic if there is a graph isomorphism
$h:\Gamma \to \Gamma^\prime$ such that $h$ carries the
orientation $\sigma$ to $\sigma^\prime$.

Later for simplicity,
we will sometimes suppress the orientation
from our notation and picture.
It is understood that from now on,
all graphs are oriented.

\subsection{The graph complex}  \label{ss:gc}

From now on, we restrict to those species $Q$ that arise 
from the symplectic theory.
Namely, we assume that $Q$ is a mated species (\ref{ss:mf}-\ref{ss:am}).
The examples in (\ref{ss:graph}) were all of this type.
Hence the vertices necessarily have valence at least 2.
We now define a chain complex $(\mathcal G,\partial_E)$ for such species.

\begin{definition}  
The $k$th chain group of $\mathcal G$, 
which we denote $\mathcal G_k$, 
is the vector space over $\mathbb Q$ generated by all oriented
graphs $(\Gamma,\sigma)$ with $k$ vertices
subject to the relations:
\begin{equation} \label{e:grel}
(\Gamma,\sigma) = -(\Gamma,-\sigma),\ \ \text{and} \ \
(\Gamma,\sigma)=(\Gamma^\prime,\sigma^\prime) \ \text{if} \
(\Gamma,\sigma)\cong(\Gamma^\prime,\sigma^\prime).
\end{equation}
\end{definition}

\noindent
A graph can have an orientation reversing automorphism, that is
$(\Gamma,\sigma)\cong (\Gamma,-\sigma)$, where $\sigma$ and $-\sigma$ are
the two orientations of $\Gamma$. 
Note that the first relation then makes the class of such a graph 
equal to zero.

\begin{definition}  
The boundary map $\partial_E : \mathcal G_k\to \mathcal G_{k-1}$
is defined using edge contractions, i.e. matings
(\ref{ss:am}).
This is where one uses that $Q$ is a mated species.
We do not contract loops.
More precisely, we have

$$\partial_E (\Gamma,\sigma) = \sum_{e \in E(\Gamma)} (\Gamma /e,
\sigma /e),$$
where $\Gamma /e$ is the graph $\Gamma$ with the edge $e$ contracted
(see figure), and $\sigma /e$ is obtained the following way: 
choose a
representative of $\sigma$ where $e$ points from vertex 1 to vertex 2,
give the new vertex arising from the contraction of $e$ the label 1, and
subtract 1 from the label of each of the other vertices; 
finally, keep the
orientations on the edges (other than $e$) unchanged.
\[
\begin{minipage}{\linewidth}
\begin{center}
\begin{minipage}{1 in}
\begin{center}
\psfrag{1}{\Huge $1$}
\psfrag{2}{\Huge $2$}
\psfrag{3}{\Huge $3$}
\psfrag{4}{\Huge $4$}
\resizebox{.8 in}{.8 in}{\includegraphics{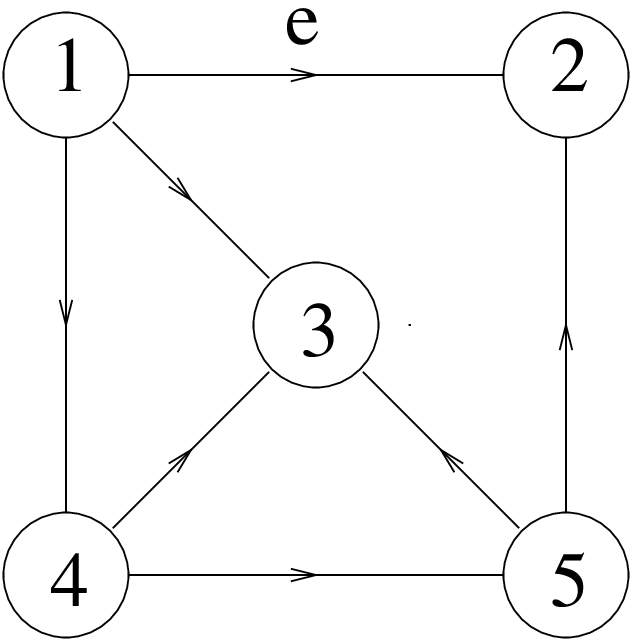}}
\end{center}
\end{minipage}
$\qquad \buildrel {\text{contract} \ e}\over \longmapsto \qquad$
\begin{minipage}{1 in}
\begin{center}
\psfrag{1}{\Huge $1$}
\psfrag{2}{\Huge $2$}
\psfrag{3}{\Huge $3$}
\psfrag{4}{\Huge $4$}
\psfrag{5}{\Huge $5$}
\resizebox{.8 in}{.8 in}{\includegraphics{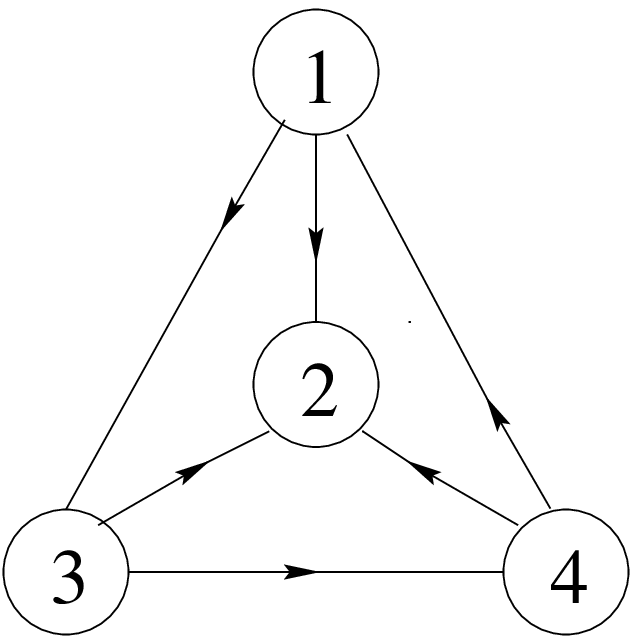}}
\end{center}
\end{minipage}.\\
\vspace{.1 in}
A mating of two vertices along an edge $e$.
\end{center}
\end{minipage}
\]
\end{definition}

\noindent
An equivalent way to describe $\sigma /e$ is the following: if the labels
on the endpoints of $e$ are $i < j$, collapse $e$, label the resulting
vertex $i$, decrease the labels greater than $j$ by one, and multiply this
orientation by $(-1)^j$ if $e$ points from $i$ to $j$, and by
$(-1)^{j+1}$ if it points from $j$ to $i$.

Note that $\sigma /e$ is well-defined and that $\partial_E$ respects the
relations in equation~\eqref{e:grel}. 

\subsubsection*{Edge contractions in examples}

We explain how edge contractions work in the examples
in (\ref{ss:graph}).
A commutative graph is an ordinary graph.
And edge contraction has the usual meaning.
For associative or ribbon graphs,
a (local) picture for an edge contraction
is given in (\ref{ss:am}).
For tree (resp. chord) graphs,
we merge the trees (resp. chords)
at the two vertices along the edge being contracted.
For a surface graph,
an edge is a loop on the surface.
To contract an edge, simply delete the loop.

\begin{lemma} 
$\partial_E^2 =0.$
\end{lemma}

\begin{proof}
If we collapse first the edge $e = (i, j)$ and then the edge $e' =
(i',j')$, we get the same graph as by collapsing $e'$ first and
then $e$, but with the opposite orientation.  This is easiest to
see using the second definition of $\sigma /e$ above; we may
assume that $j < j'$ and $i < j$, $i' < j'$. We get the same
orientation either way, but in the first case with a coefficient
$(-1)^{j}(-1)^{j'-1}$, in the second case with
$(-1)^{j'}(-1)^{j}$.
\end{proof}

Thus $(\mathcal G,\partial_E)$ is a chain complex.
We will denote its homology by $H_*(\mathcal G,\partial_E)$.

\subsection{A family of boundary operators}  \label{ss:finb}

The reader, who has never seen graph homology before,
may omit this section on a first reading.
We now define a second boundary operator $\partial_H$
that plays an important part in this theory.

Let $E_q(\Gamma)$ be the set of all \emph{quasi-edges} of $\Gamma$.
These are unordered pairs of distinct half-edges of $\Gamma$.
A quasi-edge is specified by two vertices, say $v_1$ and $v_2$, 
and two edges $e_1$ and $e_2$ incident to $v_1$ and $v_2$ respectively.
We say that
$e$ is a \emph{quasi-loop} if $v_1 = v_2$.
Observe that $E(\Gamma) \subset E_q(\Gamma)$. 

\begin{definition}  
The boundary map $\partial_H : \mathcal G_k\to \mathcal G_{k-1}$
is defined by contracting quasi-edges, which are not edges.
We do not contract quasi-loops. We have
$$\partial_{H} (\Gamma,\sigma) = \sum_{e \in E_q(\Gamma) \setminus E(\Gamma)} 
(\Gamma /e,\sigma /e),$$
where 
$\Gamma /e$ and $\sigma /e$ 
are defined as follows.

1. Choose a representative for $\sigma$ where
$v_1$ and $v_2$ have labels $1$ and $2$ respectively
and the arrow on edge $e_1$
(resp. $e_2$) points out of $v_1$ (resp. into $v_2$).

2. Cut the edges $e_1$ and $e_2$ and join them 
so that the quasi-edge $e$ and its partner $f$ turn into edges. 
In effect,
$e_1$ and $e_2$ have been replaced by two new edges 
(which get an induced orientation).
\[
\begin{minipage}{\linewidth}
\begin{center}
\begin{minipage}{1.6 in}
\begin{center}
\psfrag{v1}{\Huge $v_1$}
\psfrag{v2}{\Huge $v_2$}
\psfrag{e1}{\Huge $e_1$}
\psfrag{e2}{\Huge $e_2$}
\psfrag{e}{\Huge $e$}
\psfrag{f}{\Huge $f$}
\psfrag{pi}{\Huge $$}
\psfrag{qi}{\Huge $$}
\resizebox{1.4 in}{.9 in}{\includegraphics{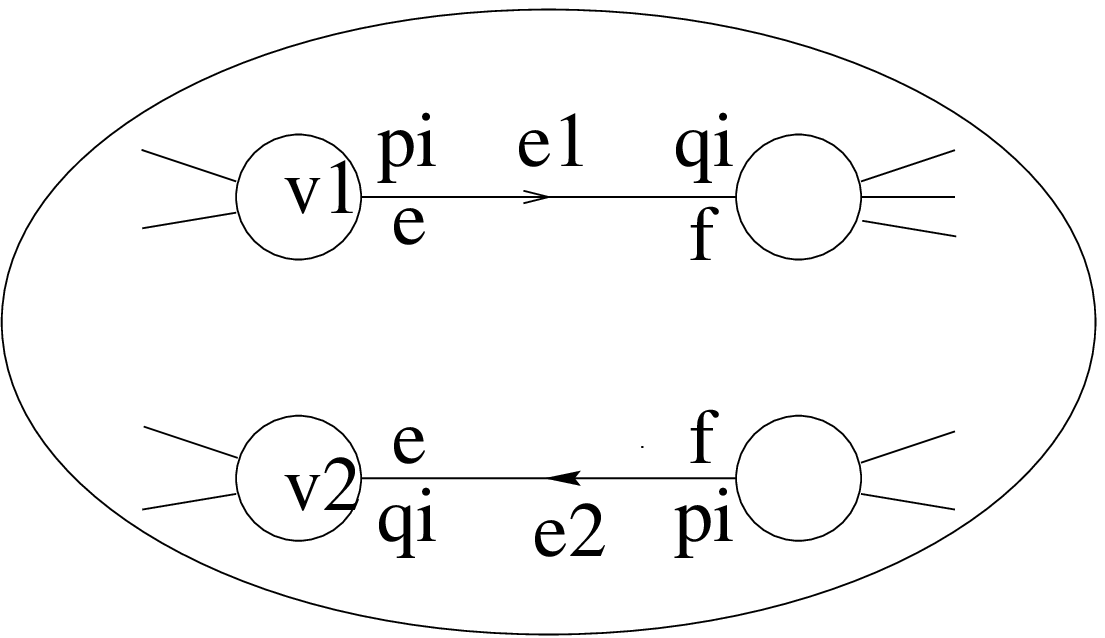}}
\end{center}
\end{minipage}
$\qquad \buildrel {\text{contract} \ e}\over \longmapsto \qquad$
\begin{minipage}{1.6 in}
\begin{center}
\psfrag{v1}{\Huge $v_1$}
\psfrag{v2}{\Huge $v_2$}
\psfrag{e}{\Huge $e$}
\psfrag{f}{\Huge $f$}
\resizebox{1.4 in}{.9 in}{\includegraphics{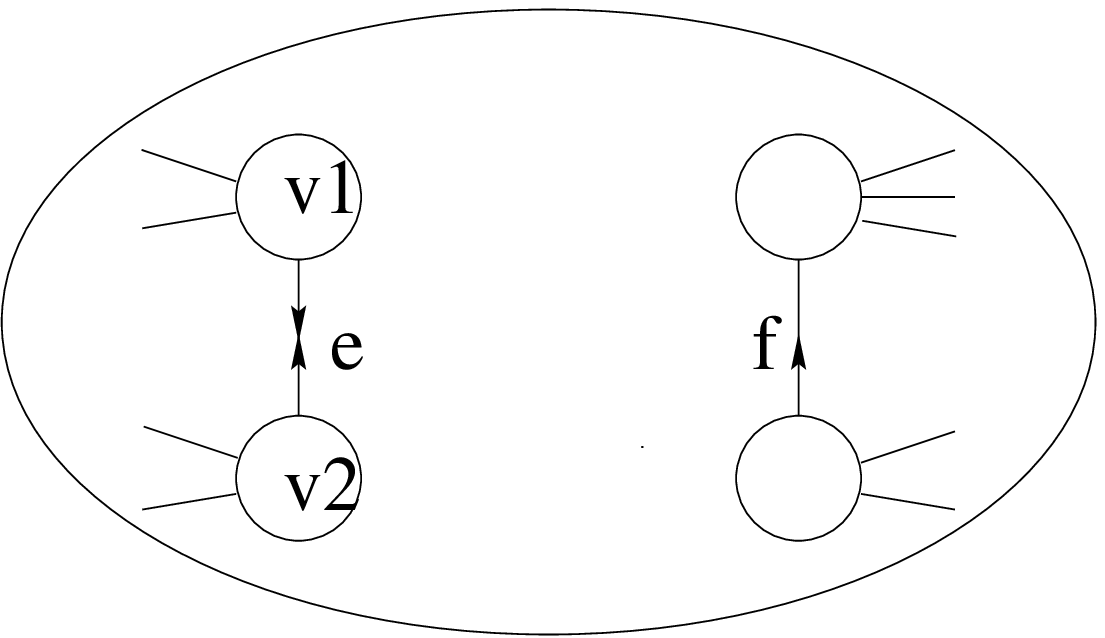}}
\end{center}
\end{minipage}.\\
\vspace{.1 in} A mating along a quasi-edge $e$.
\end{center}
\end{minipage}
\]

3. Collapse the edge joining $v_1$ and $v_2$. 
Thus $v_1$ and $v_2$ become a single vertex, say $v$. 
The resulting graph is $\Gamma /e$. 
Note that the edges of $\Gamma /e$ are already oriented. 
To get an order on the vertices, 
label the new vertex $v$ by $1$ and
subtract $1$ from the labels of the other vertices. 
The resulting orientation is $\sigma /e$. 
\end{definition}

\begin{definition}  
We now define a family of boundary maps by the formula
$$\partial_n = 2n \partial_E + \partial_H.$$
To put in words, $\partial_n(\Gamma)$
is obtained by contracting quasi-edges of $\Gamma$.
And if the quasi-edge is an edge 
then we multiply by a factor of $2n$.
\end{definition}

The fact
that $\partial_n = 2n \partial_E + \partial_H$ 
is a boundary operator for all $n$
will follow from the proof 
of Theorem~\ref{t:fin} (Section~\ref{s:fin}).
This automatically implies that $\partial_E^2=\partial_H^2=0$
and $\partial_E \partial_H = - \partial_H \partial_E$.
In other words, $\partial_E$ and $\partial_H$ together
span an abelian Lie superalgebra
of operators on $\mathcal G$.
We will call $\partial_n$
the finite boundary operator
and write the chain complex as $(\mathcal G,\partial_n)$.

\begin{remark}
If the boundary operator is not stated explicitly,
then we take it to be the stable boundary operator $\partial_E$.
\end{remark}

\subsection{An important simplification}  \label{ss:is}

Note that the boundary maps $\partial_E$ and $\partial_H$ preserve
the Euler characteristic $\chi$ of any graph $\Gamma$.
This is because they decrease the number of edges and vertices
both by 1.
So the graph complex,
irrespective of the boundary map used,
splits as a direct sum of chain complexes
$\mathcal G^{(r)}$:
\begin{gather*}
\cdots \to \mathcal G^{(1)}_3 \to \mathcal G_2^{(1)} \to \mathcal G^{(1)}_1 \to 0\\
\cdots \to \mathcal G^{(2)}_3 \to \mathcal G_2^{(2)} \to \mathcal G^{(2)}_1 \to 0\\
\cdots \to \mathcal G^{(3)}_3 \to \mathcal G_2^{(3)} \to \mathcal G^{(3)}_1 \to 0\\
\hspace{1em}\vdots \hspace{3em} \vdots \hspace{3em}\vdots
\end{gather*}
The useful point is that the chain groups 
$\mathcal G_k^{(r)}$ are finite dimensional.
Instead of $\chi$,
we have used the superscript $r = 1- \chi$ for convenience.
If the graph $\Gamma$ is connected
then $r=\dim_{\Q}H_1(\Gamma, \Q)$ 
is just the rank of $\Gamma$.

\subsection{Other relevant graph complexes}  \label{ss:cc}

From now on, assume that the species $Q$
is based on sets rather than vector spaces.
All the examples
in this paper are of this type.
We first give a preliminary definition.

\begin{definition}  
We say a vertex in a $Q$-graph is fake 
if it corresponds to the element $uu$ in $Q[2]$,
the degree 2 piece of $Q$, see (\ref{ss:uu}).
Else we say that the vertex is real.
\[
\begin{minipage}{4 in}
\begin{center}
\resizebox{1 in}{.6 in}{\includegraphics{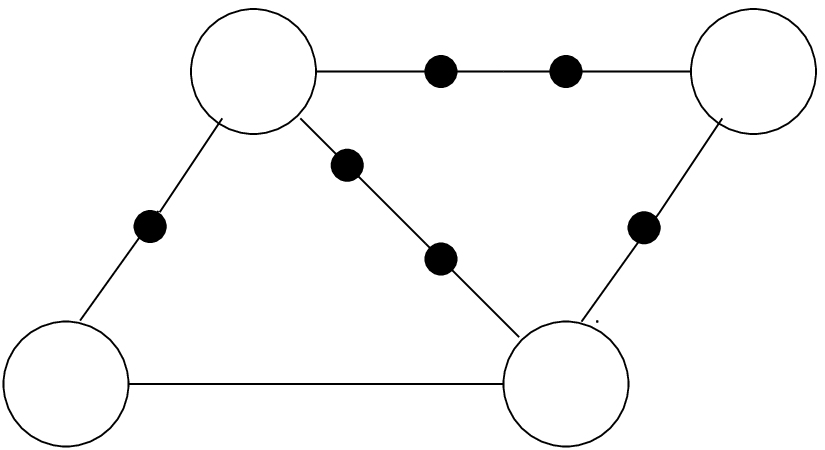}}\\
\vspace{.1 in} The dots stand for the fake vertices.
\end{center}
\end{minipage}
\]
\end{definition}
Note that the fake vertices are always bivalent
and they all look identical.
Furthermore, they behave as unit elements 
in contractions (matings).
So in this sense, they are extremely inert.

For all the examples in (\ref{ss:graph}),
except the surface species $ss$,
$Q[2]$ is a singleton.
For example,
there is only one cyclic order or chord diagram or unrooted tree
on 2 letters.
Hence, in these cases, 
the fake vertices are precisely the bivalent vertices.
And real vertices are exactly those with valence greater than 2.
For the surface species, 
$ss[2]$ consists of compact orientable surfaces 
whose boundary is two disjoint circles.
And there are clearly many such.
Among them, cylinders are precisely the fake vertices,
while the rest are real.

It can even happen that
contracting an edge joining two real vertices 
produces a fake vertex.
An example is given in Section~\ref{s:gp}.
In this case, the two real vertices must be necessarily bivalent.
In any case, the boundary map can reduce the number
of real vertices of a graph by at most 2.

We now look at the various subcomplexes of $(\mathcal G,\partial_E)$,
that appear in the main theorem and its proof.
They are all defined by putting some restriction
on the type of graphs that are allowed.
\begin{itemize}

\item
$\mathcal T =$ graphs all of whose vertices are bivalent.

\item
$\mathcal F =$ graphs all of whose vertices are fake.

\item
$\mathcal C =$ connected graphs.

\item
$Q\mathcal G =$ connected graphs all of whose vertices are real
and with at least one vertex of degree greater than 2.

\item
$\mathcal B =$ connected graphs all of whose vertices are bivalent.

\item
$\mathcal P =$ connected graphs all of whose vertices are fake.
\end{itemize}

It is clear that $\partial_E$ preserves the above chain groups
except $Q\mathcal G$.
The trouble with $Q\mathcal G$ is that an edge contraction
can result in a fake vertex.
Hence for $Q\mathcal G$, modify the definition of $\partial_E$
such that edge contractions that create fake vertices are ignored.  
So, strictly speaking, it is not a subcomplex of $(\mathcal G,\partial_E)$.
In keeping with Kontsevich's original definition,
we will call $H_*(Q\mathcal G,\partial_E)$
the graph homology of the mated species $Q$.

Note that $\mathcal F$ and $\mathcal P$
do not depend on the species $Q$.
The complex $\mathcal P$ just consists of polygons and
$\mathcal F$ consists of disjoint unions of polygons.
Further, since fake vertices are always bivalent,
the complex $\mathcal F$ (resp. $\mathcal P$)
is a subcomplex of $\mathcal T$ (resp. $\mathcal B$).
The containment is strict if $Q[2]$ has more than 1 element,
as for the surface or graph species.
We discuss another example of this kind in Section~\ref{s:gp}. 

\subsection{A reduction step}  

The chain complex that appears in the statement
of the main theorem is $(Q\mathcal G,\partial_E)$.
However, the one that appears naturally in the proof is
$(\mathcal C,\partial_E)$.
The reduction from the homology 
of this larger complex to 
graph homology $H_*(Q\mathcal G,\partial_E)$
is the content of Proposition~\ref{p:ss}.
It will be the last step in the proof
of the main theorem.
Though the statement of the proposition 
is quite intuitive,
its proof is a little technical and may be skipped
on a first reading.

\begin{proposition}  \label{p:ss}
$H_*(\mathcal C) = H_*(Q\mathcal G) \bigoplus H_*(\mathcal B)$.
\end{proposition}
\begin{proof}
The bivalent graph complex $\mathcal B$ is a direct summand 
of the graph complex $\mathcal C$.
The complement is the subcomplex consisting of connected graphs
containing at least one vertex of degree $\geq 3$.
Denote this complex by $\mathcal D$, graded as usual by the number of
vertices. 
The difference between $\mathcal D$ and $Q\mathcal G$ is that
graphs in $\mathcal D$ are allowed to contain fake vertices.
To finish the proof, we need to show that 
$H_*(\mathcal D) = H_*(Q\mathcal G)$.
This can be done by a spectral sequence argument.
For a treatment of spectral sequences, see \cite{\mccleary}.

We have already defined real and fake vertices.
In addition, call an edge fake if it is incident to a fake vertex.
Next, define a filtration on $\mathcal D$ by
$F_p \mathcal D_m =$ graphs in $\mathcal D_m$ with upto $p$ real vertices.
Recall that for the spectral sequence associated to a filtration,
$E_{p q}^0 = F_p \mathcal D_{p+q} \big/ F_{p-1} \mathcal D_{p+q}$.
In our case, $E_{p q}^0$ can be described as the span of oriented
graphs with $p$ 
real and $q$ fake vertices.
The vertical maps on the $E^0$ page, namely,
$d_0 : E_{p q}^0 \mapsto E_{p q-1}^0$
are defined exactly as 
the boundary map for an oriented graph, 
except that now we are
allowed to collapse only the fake edges.
Note that the positive $X$-axis of the $E^0$ page consists exactly of
graphs with no fake vertices. In other words, they are 
the chain groups of the graph complex $Q\mathcal G$.

\medskip
\noindent{\bf Claim}. The homology on the $E^0$ page is trivial except
at $q = 0$, that is, the $X$-axis. Equivalently, 
$E_{p q}^1 = 0$ for $q > 0$. 
In addition, we claim that $E_{p 0}^1 = Q\mathcal G_p$.

Assuming the claim,
it is clear that the induced boundary map on the $E^1$ page 
$d_1 : E_{p 0}^1 \mapsto E_{p-1 0}^1$ coincides with the boundary map
defined on $Q\mathcal G$. 
Namely, edge contractions that create fake vertices are ignored.
The rest of the $E^1$ page is zero.
Thus the spectral sequence associated to our filtration of $\mathcal D$
collapses at the first term to the graph complex $Q\mathcal G$.
This shows that $H_*(\mathcal D) = H_*(Q\mathcal G)$.

\medskip
\noindent{\bf Proof of the claim}. 
To every connected graph containing at least one real vertex, 
one can
associate a new graph all of whose vertices are real. 
The new graph is just the old graph with all the fake vertices removed. 
This gives us
an equivalence relation on the set of all connected graphs containing
a real vertex. 
The equivalence classes are indexed by graphs all of
whose vertices are real.
Now consider the complex 
$$\ldots\to E_{p q}^0\buildrel d_0\over \to E_{p q-1}^0\buildrel 
d_0\over \to \ldots \buildrel d_0\over \to E_{p 0}^0 = Q\mathcal G_p 
\buildrel d_0\over \to 0$$ 
whose homology we want to compute.
This splits as a direct sum of subcomplexes one for each graph
$\Gamma$ in $E_{p 0}^0 = Q\mathcal G_p$. 
Call the subcomplex corresponding
to $\Gamma$ the standard complex of $\Gamma$. 
Note that the graphs
that occur in the standard complex of $\Gamma$ are precisely those
that lie in the same equivalence class as $\Gamma$. 
We need to show that the
only non-trivial homology occurs in degree $0$.

To understand this complex, 
first look at the standard complex of a single edge. 
It has dimension $1$ in each degree $k \geq 0$, 
because there is a unique (upto isomorphism) way to put $k$ points in
the interior of an edge. 
The differential in this standard complex
is zero for $k=0$ and all odd $k$ and 
an isomorphism for positive even $k$. 
Hence it has non-trivial homology only in degree $0$.

Now the standard complex of $\Gamma$ is 
the tensor product of the
standard  complexes of all its edges
modulo the action of the finite group $\Aut(\Gamma$). 
Hence by the Kunneth formula 
and the fact that finite groups have trivial rational homology
\cite{\kbrown},
we are done.
\end{proof}

\begin{corollary}
Let $Q$ be a mated species, 
where all bivalent vertices in any graph are fake. 
Or equivalently,
let $Q$ be obtained from a reversible operad
satisfying $P[1]=\Q$.
Then
$H_*(\mathcal C) = H_*(Q\mathcal G) \bigoplus H_*(\mathcal P)$.
\end{corollary}

\section{Graph homology for groups}  \label{s:gp}

We now digress to give an example of the theory
that is based on groups.
It can be defined more generally for algebras with involution;
we then recover dihedral homology \cite{\loday}.
We hope that this example will further clarify the concepts
discussed in previous sections.
An attraction of this example is that
it is tractable to computations.
Some known results on dihedral homology 
can be obtained this way.

\subsection{Groups as reversible operads} 
Let $K$ be any finite group.
Let $P$ be the operad whose elements 
are given by 
\begin{minipage}{.6 in}
\begin{center}
\psfrag{u}{\Huge $g$}
\psfrag{a}{\Huge $a$}
\psfrag{b}{\Huge $b$}
\resizebox{.5 in}{.2 in}{\includegraphics{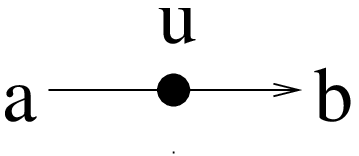}}
\end{center}
\end{minipage}
for $g \in K$,
with substitution being group multiplication.
In other words, $P[1]=\Q K$,
the group algebra of $K$ over $\Q$, and
$P[2]=P[3]=\ldots=0$.
If $K$ is the trivial group
then $P$ is just the unit operad $u$,
see (\ref {sss:u}).
More generally, one can take $P[1]$ 
to be any algebra with a unit.

For this example,
$P$ is reversible (\ref{ss:ro}),
if there exists a map 
$^{*}: K \to K$ which satisfies
$1^*=1$ and $(g h)^* = h^* g^*$
for any $g,h \in K$.
Natural candidates for the $^{*}$ map are 
the identity
(if $K$ is abelian)
and the inverse map.
Thus, for $K$ abelian,
there are two distinct ways of reversal.
The reversal map is given by 
$$
r_{a,b} \left(
\begin{minipage}{.6 in}
\begin{center}
\psfrag{u}{\Huge $g$}
\psfrag{a}{\Huge $a$}
\psfrag{b}{\Huge $b$}
\resizebox{.5 in}{.2 in}{\includegraphics{ur.eps}}
\end{center}
\end{minipage}
\right)
=
\begin{minipage}{.6 in}
\begin{center}
\psfrag{a}{\Huge $a$}
\psfrag{b}{\Huge $b$}
\psfrag{u}{\Huge $g^*$}
\resizebox{.5 in}{.2 in}{\includegraphics{u.eps}}
\end{center}
\end{minipage}
.$$

Applying the mating functor to $P$ gives us a mated species $Q$.
This species $Q$ lives entirely in degree 2.
We leave it to the reader to check that an element of 
$Q[\{a,b\}]$ can be specified by the picture
\begin{equation}  \label{e:gr}
\begin{minipage}{.7 in}
\begin{center}
\psfrag{g}{\Huge $g$}
\psfrag{a}{\Huge $a$}
\psfrag{b}{\Huge $b$}
\resizebox{.7 in}{.3 in}{\includegraphics{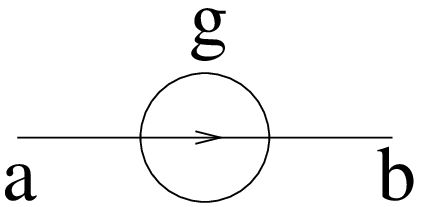}}
\end{center}
\end{minipage}
\ \ \ = \ \ \
\begin{minipage}{.7 in}
\begin{center}
\psfrag{g}{\Huge $g^*$}
\psfrag{a}{\Huge $a$}
\psfrag{b}{\Huge $b$}
\resizebox{.7 in}{.3 in}{\includegraphics{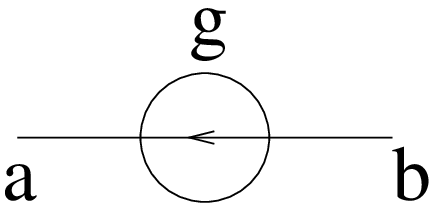}}
\end{center}
\end{minipage}
.
\end{equation}

\subsection{Graph homology for groups}  \label{ss:gp}

For the above species $Q$, the corresponding graphs
are necessarily bivalent (polygons),
with each vertex labelled (more or less) by a group element.
The vertices labelled by the unit element 1 are fake
while the rest are real.

We now explain how an edge contraction (mating) works.
\[
\begin{minipage}{1.2 in}
\begin{center}
\psfrag{g}{\Huge $g$}
\psfrag{h}{\Huge $h$}
\psfrag{e}{\Huge $e$}
\resizebox{1.2 in}{.6 in}{\includegraphics{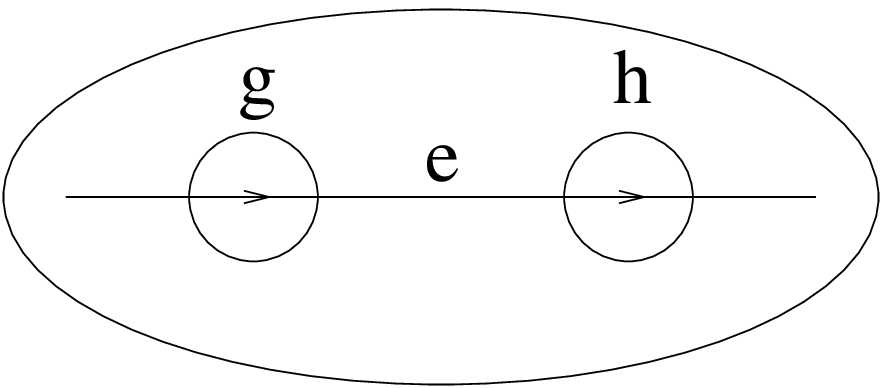}}
\end{center}
\end{minipage}
\qquad \buildrel {\text{contract} \ e}\over \longmapsto \qquad
\begin{minipage}{1.1 in}
\begin{center}
\psfrag{gh}{\Huge $gh$}
\resizebox{1 in}{.6 in}{\includegraphics{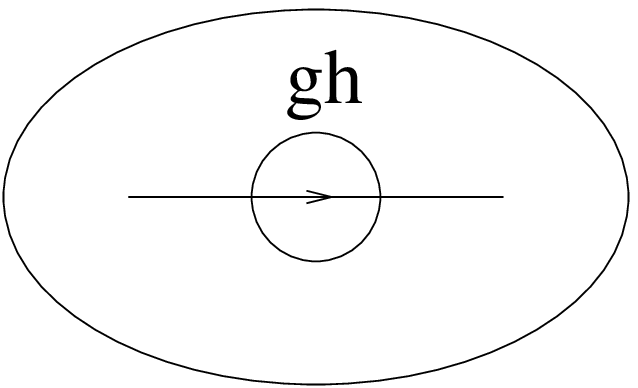}}
\end{center}
\end{minipage}
.
\]
Apply relation~\eqref{e:gr} if necessary,
so that both vertices incident to the edge $e$ point in the same direction.
Contract it and make the new vertex also point in the same direction.
And label it with the product of the two labels
(in the order specified by the direction).

We will call the homology of this bivalent complex
as the graph homology of the group $K$
and denote it $H_*(K,^{*})$.
In the notation of Section~\ref{s:gh},
it coincides with both $H_*(\mathcal C)$ and $H_*(\mathcal B)$ .

\subsection{Graph homology computations}  \label{ss:cal}

The problem of computing graph homology 
seems to be difficult in general.
The only instance where computations 
have been made is the commutative case (\ref{sss:c});
see the thesis of Ferenc Gerlits~\cite{\ferenc}.
This is the case where we are dealing with usual graphs.

We give two instances where graph homology $H_*(\mathcal C)$
can be completely computed.
They are the first two cases of the example
under discussion.
It is easy to generalise them but we do not do it here.

\subsubsection{The trivial group} \label{sss:ctg}

If $K$ is the trivial group then $P=u$ and $Q=uu$,
see (\ref{ss:uu}).
The computation for this case,
which coincides with $H_*(\mathcal P)$, see~(\ref{ss:cc}),
is discussed in~\cite{\god}.
The graphs are polygons with only fake vertices.
These graphs have two kinds of automorphisms: rotations and reflections.
It is easy to see that a simple rotation is orientation reversing if and
only if $k$, the number of vertices, is even, and a reflection reverses
the orientation iff $k \equiv 1$ or $2$ modulo $4$.  So there is one
non-zero chain in each degree of the form $4i+3$, and hence all the
boundary maps are zero, therefore these are the homology groups also.
Hence 
$$
  H_k(\{1\},\text{identity}) \ = \left\{ \begin{array}{c l}
    \Q & \mbox{ if $k \equiv 3 \ (\text{mod} \ 4)$}, \\
    0  & \mbox{ otherwise.} 
                    \end{array} \right. $$
By Corollary to Theorem~\ref{t:con} (Section~\ref{s:con}),
this is same as the primitive homology
of the Lie algebra ${{{\mathfrak{sp}}(2 \infty)}}$.

\subsubsection{The group $\Z_2$}  

For this case, the graphs are polygons with two types of vertices
corresponding to the two group elements.
Computing the homology directly from the definition is not easy.
So we rerun the spectral sequence argument
(Proposition~\ref{p:ss})
with some modification.
The stable page is $E^2$ with non-zero terms
only on the positive X and Y axis.
The conclusion is
$$
  H_k(\Z_2,\text{identity}) \ = \left\{ \begin{array}{c l}
    \Q \oplus \Q & \mbox{ $k \equiv 3 \ (\text{mod} \ 4)$}, \\
    0           & \mbox{ otherwise.} 
                    \end{array} \right. $$
This example, we hope, gives an idea 
of the complexity of computing graph homology.

\section{Graph cohomology}  \label{s:mgh}

The ideas in this section are based on a letter of Kontsevich.
It leads to the birth of many interesting operations
on graphs.
I thank Jim Conant for helping me understand its contents
and also for providing his notes related to this.

We begin by defining graph cohomology.
The homology and cohomology are related 
by an interesting and highly non-trivial pairing on graphs.
This is explained in (\ref{ss:pg}-\ref{ss:ad}).
Later in (\ref{ss:dmg}),
we use it to define a deformation map on graphs.
These ideas will be used in Section~\ref{s:sta} 
in the proof of Theorem~\ref{t:sta}
that deals with stability.
Throughout this section,
we assume that $Q$ is a mated species based on sets. 

\subsection{The blowup coboundary operator $\delta_E$}  \label{ss:stacb}

Let $\Gamma$ be an oriented $Q$-graph and
$I(\Gamma)$ 
be the set of its ideal edges.
These are the edges that are ``present''
in the internal structure of the vertices of $\Gamma$,
see (\ref{ss:mf}).

\begin{definition}  
The coboundary map $\delta_E : \mathcal G_k\to \mathcal G_{k+1}$
is defined using ideal edge expansions, i.e. breakups.
It is given by the formula

$$\delta_E (\Gamma,\sigma) = \sum_{e \in I(\Gamma)} (\Gamma \backslash e,
\sigma \backslash e),$$
where $\Gamma \backslash e$ is the graph $\Gamma$ 
with the ideal edge $e$ expanded
(see figure), and 
$\sigma \backslash e$ is obtained the following way: 
choose a representative of $\sigma$ 
where the vertex of the ideal edge $e$ has label 1,
give the two new vertices arising from the breakup
the labels 1 and 2,
add 1 to the label of each of the other vertices, and
direct the newly formed edge $e$ from vertex 1 to 2;
finally, keep the
orientations on the other edges unchanged.
\[
\begin{minipage}{\linewidth}
\begin{center}
\begin{minipage}{1.2 in}
\begin{center}
\psfrag{1}{\Huge $1$}
\psfrag{2}{\Huge $2$}
\psfrag{e}{\Huge $e$}
\resizebox{1 in}{.6 in}{\includegraphics{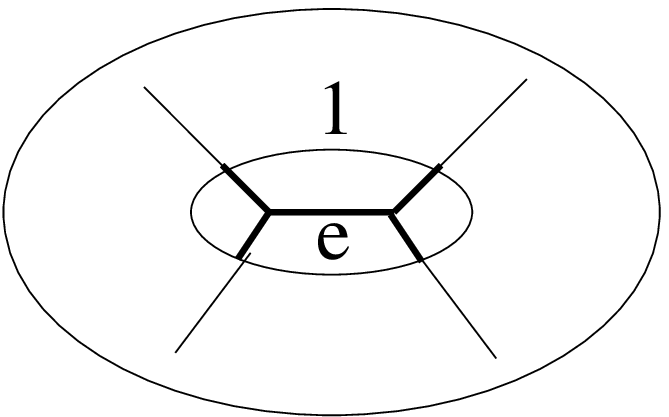}}
\end{center}
\end{minipage}
$\qquad \buildrel {\text{expand} \ e}\over \longmapsto \qquad$
\begin{minipage}{1.2 in}
\begin{center}
\psfrag{1}{\Huge $1$}
\psfrag{2}{\Huge $2$}
\psfrag{e}{\Huge $e$}
\resizebox{1 in}{.6 in}{\includegraphics{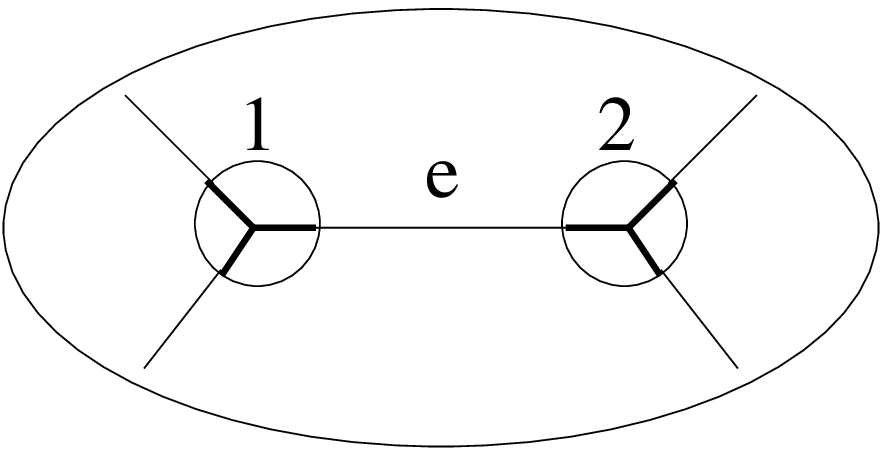}}
\end{center}
\end{minipage}\\
\vspace{.1in}
A local picture of a breakup along an ideal edge $e$.
\end{center}
\end{minipage}
\]
\end{definition}
Similarly, one can define a coboundary operator
$\delta_H$ as the co-analogue of $\partial_H$,
see (\ref{ss:finb}).
We will not deal with $\delta_H$ in this paper.

\subsection{A family of pairings $M(n)$}  \label{ss:pg}

For every $n$,
we will define a pairing
$M(n) : \mathcal G \otimes \mathcal G \to \Q$.
The definition, though explicit, will be somewhat complicated.
A way to understand it is given in
(\ref{ss:pl}-\ref{ss:pr}), where we derive it 
as a restriction of a simpler pairing
defined on a larger space.
The reader, who is more interested in Lie algebras
or the main theorem rather than just graphs,
may read that part first,
referring back as necessary.

Let $(\Gamma_1,\sigma_1),(\Gamma_2,\sigma_2) \in \mathcal G$
be two oriented graphs.
For simplicity,
we suppress orientations from the notation.

\begin{definition}  \label{d:gm}
A matching $m:\Gamma_1 \to \Gamma_2$ is a bijection
between the vertex sets
$V(\Gamma_1)$ and $V(\Gamma_2)$
that preserves the internal structure of the vertices.
We imagine this as an overlaying of the vertex sets of
$\Gamma_1$ and $\Gamma_2$.
\[
\begin{minipage}{1.6 in}
\begin{center}
\begin{minipage}{1.2 in}
\begin{center}
\resizebox{1.2 in}{1.2 in}{\includegraphics{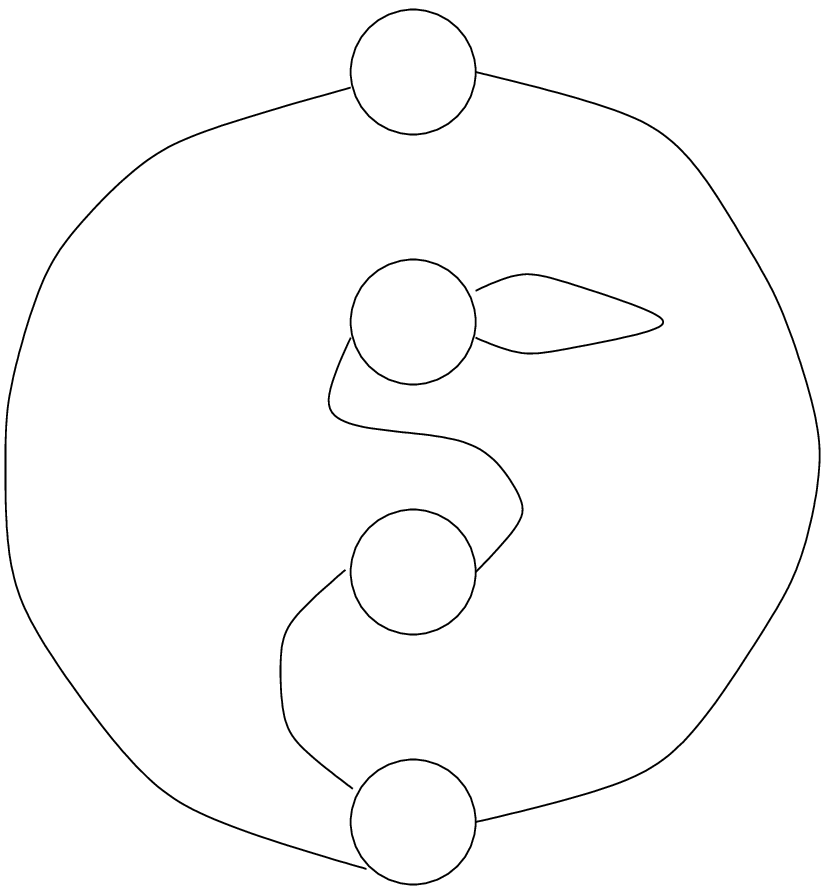}}
\end{center}
\end{minipage}\\
\vspace{.1in}
$\Gamma_1$
\end{center}
\end{minipage}
\begin{minipage}{1.6 in}
\begin{center}
\begin{minipage}{1.2 in}
\begin{center}
\resizebox{.7 in}{1.2 in}{\includegraphics{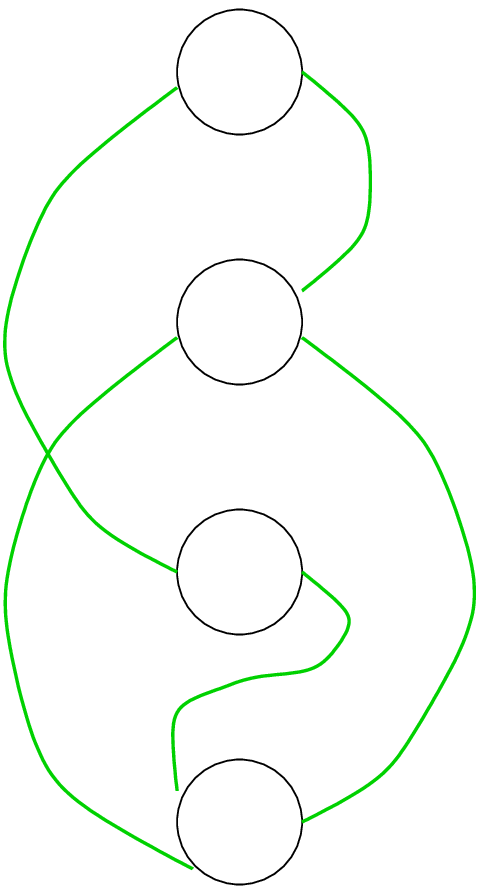}}
\end{center}
\end{minipage}\\
\vspace{.1in}
$\Gamma_2$
\end{center}
\end{minipage}
\begin{minipage}{1.6 in}
\begin{center}
\begin{minipage}{1.2 in}
\begin{center}
\psfrag{d}{\Huge $d$}
\psfrag{c}{\Huge $c$}
\psfrag{b}{\Huge $b$}
\psfrag{a}{\Huge $a$}
\resizebox{1.2 in}{1.2 in}{\includegraphics{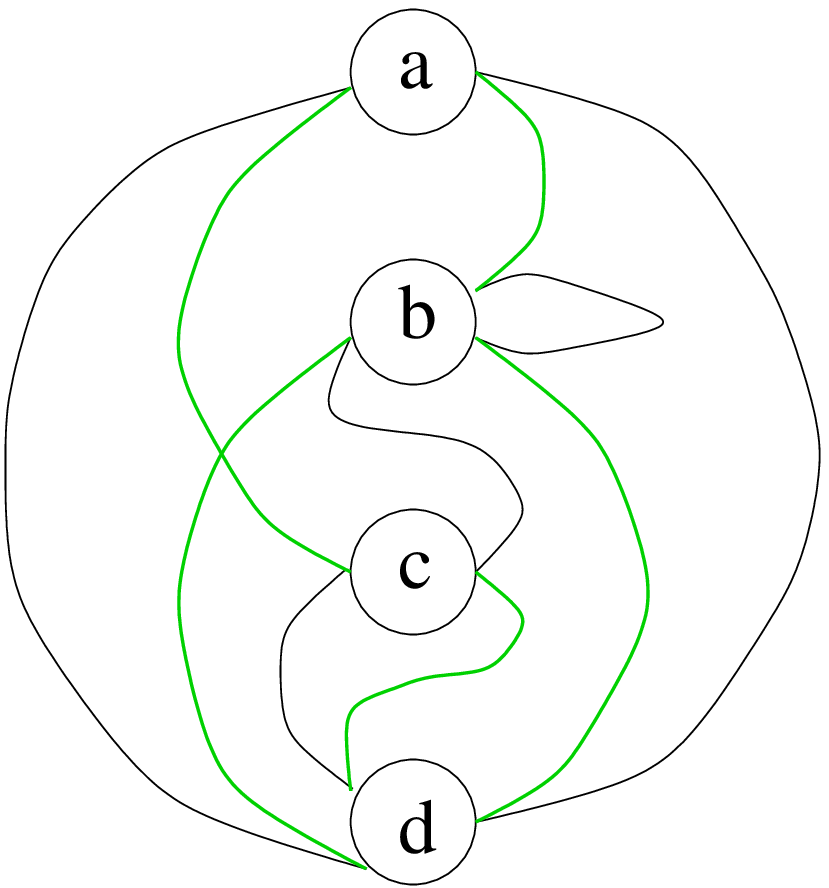}}
\end{center}
\end{minipage}\\
\vspace{.1in}
Matching of $\Gamma_1$ and $\Gamma_2$.
\end{center}
\end{minipage}
\]
An example of a matching is shown above.
For clarity, orientations 
and the internal structure of the vertices
have been omitted.
It is clear that for a matching to exist,
$\Gamma_1$ and $\Gamma_2$
must have the same number of vertices and edges.
\end{definition}

\begin{definition} \label{d:csm}
We now define the number of components 
$c(m)$ and $\sign(m)$
of a matching
$m:\Gamma_1 \to \Gamma_2$.
Delete the vertices from the overlaying of
$\Gamma_1$ and $\Gamma_2$ specified by $m$.
We are left with a disjoint union of even sided polygons
whose edges alternate between those of
$\Gamma_1$ and $\Gamma_2$.
The quantity $c(m)$ counts the number of these polygons.
In the example above, we get a hexagon and a square;
so $c(m)=2$.
\[
\begin{minipage}{.7 in}
\begin{center}
\psfrag{a}{\Huge $d$}
\psfrag{b}{\Huge $a$}
\psfrag{c}{\Huge $c$}
\psfrag{d}{\Huge $d$}
\psfrag{e}{\Huge $c$}
\psfrag{f}{\Huge $b$}
\resizebox{.7 in}{.7 in}{\includegraphics{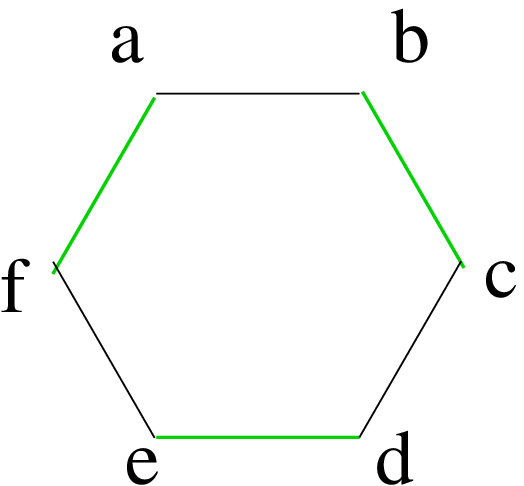}}
\end{center}
\end{minipage}
\hspace{1 in}
\begin{minipage}{.6 in}
\begin{center}
\psfrag{a}{\Huge $b$}
\psfrag{b}{\Huge $b$}
\psfrag{c}{\Huge $a$}
\psfrag{d}{\Huge $d$}
\resizebox{.5 in}{.5 in}{\includegraphics{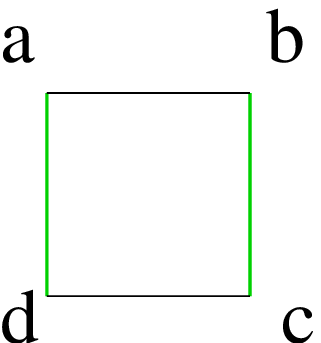}}
\end{center}
\end{minipage}
\]
The purpose of $\sign(m)$ is to take 
orientations into account.
Choose representatives for $\sigma_1$ and $\sigma_2$.
The matching $m$ gives a bijection of the vertex sets.
First note the sign of this permutation.
Then for each polygon, fix say the clockwise direction
for $\Gamma_1$ and anticlockwise direction for $\Gamma_2$. 
And write a minus sign for all edges that are out of order.
The product of all signs gives $\sign(m)$.
\end{definition}

\begin{definition}  
Define the pairing 
$M(n) : \mathcal G \otimes \mathcal G \to \Q$ by the formula 
\begin{equation}  \label{e:pg}
M(n)(\Gamma_1,\Gamma_2) = \underset{m:\Gamma_1 \to \Gamma_2}
{\sum} \sign(m) (2n)^{c(m)},
\end{equation}
where the sum is over all matchings
$m:\Gamma_1 \to \Gamma_2$
and $\sign(m)$ and $c(m)$ 
are the sign and the number of components of the matching $m$.
\end{definition}

\begin{remark}
For $\Gamma_1$ and $\Gamma_2$ fixed,
the maximum value of $c(m)$ is
$e=$ number of edges of $\Gamma_1$ (and $\Gamma_2$).
This happens iff every polygon in the matching $m$
has exactly two sides, i.e.
iff $m:\Gamma_1 \to \Gamma_2$ is an isomorphism.
We will call this a perfect matching.
And in this case, 
the coefficient of $(2n)^e$, upto sign, is
$|\Aut(\Gamma_1)| = |\Aut(\Gamma_2)|$.
%$|\text{Aut}(\Gamma_1)| = |\text{Aut}(\Gamma_2)|$.
\end{remark}

\subsection{The adjoint property} \label{ss:ad}

We now relate the finite boundary operator $\partial_n$
in (\ref{ss:finb})
to the (stable) coboundary operator $\delta_E$
in (\ref{ss:stacb}).
There is a simpler relation between $\partial_E$
and $\delta_E$ which can be derived from this one,
by looking at the leading coefficient.
This is explained in (\ref{ss:dmg}).
Recall that $\partial_n$ is defined by contracting quasi-edges
and $\delta_E$ is defined by expanding ideal edges.

\begin{proposition}  \label{p:ad}
The maps $\partial_n$ and $\delta_E$ are adjoints
with respect to the pairing $M(n)$. In other words,
$M(n)(\partial_n \Gamma_1,\Gamma_2)=M(n)(\Gamma_1,\delta_E \Gamma_2).$
\end{proposition}
\begin{proof}
To prove the above identity, 
we express both sides as weighted state sums,
and then give a bijection of the state space in the LHS
with that in the RHS such that it respects weights.
Note that for either side to be nonzero,
$\Gamma_1$ must have one vertex and edge more than $\Gamma_2$.

Define a state $S_L$ in the LHS to be a pair $(e,m)$
of a quasi-edge $e \in E_q(\Gamma_1)$
and a matching $m:\Gamma_1 /e \to \Gamma_2$.
Then by the definition of the boundary map $\partial_n$ 
and the pairing $M(n)$, the
LHS = $M(n)(\partial_n \Gamma_1,\Gamma_2)=\sum_{S_L} w(S_L), \ 
\text{where}$
$$
  w(S_L) \ = \left\{ \begin{array}{c l}
    \sign(m) (2n)^{c(m)+1} &
\mbox{ if $e \in E(\Gamma_1)$}, \\
    \sign(m) (2n)^{c(m)} &
\mbox{ if $e \in E_q(\Gamma_1) \setminus E(\Gamma_1).$}
                    \end{array} \right. $$
An example of a state $S_L$ for
$e \in E_q(\Gamma) \setminus E(\Gamma)$
is shown in the figure.
The quasi-edge $f$ is the partner of $e$.
The graph $\Gamma_2$ is not shown separately,
since it is visible from the overlaying.
Also, in the graph $\Gamma_1 /e$,
the edge $e$ has been shown as an ideal edge.
\[
\begin{minipage}[b]{2.4 in}
\begin{center}
\psfrag{e}{\Huge $e$}
\psfrag{f}{\Huge $f$}
\resizebox{1 in}{1 in}{\includegraphics{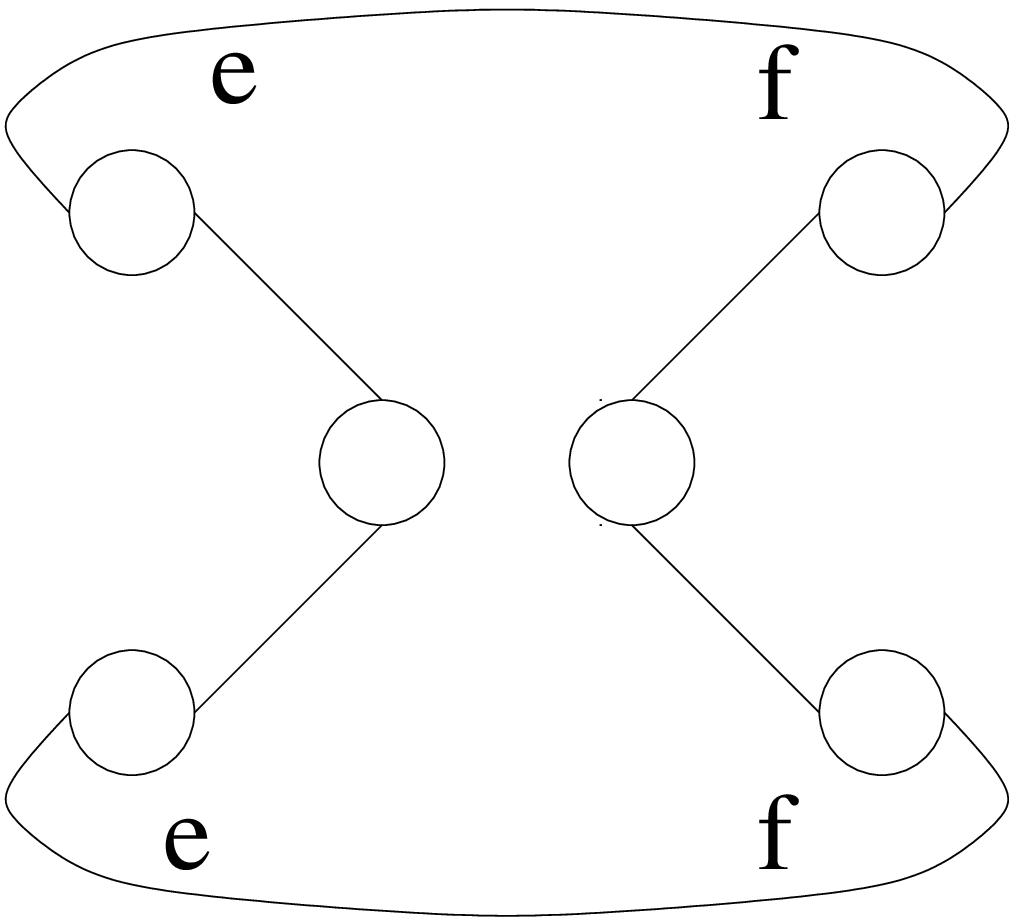}}\\
$\Gamma_1$ with a quasi-edge $e$.
\end{center}
\end{minipage}
\begin{minipage}[b]{2.4 in}
\begin{center}
\psfrag{e}{\Huge $e$}
\psfrag{f}{\Huge $f$}
\resizebox{1.4 in}{.8 in}{\includegraphics{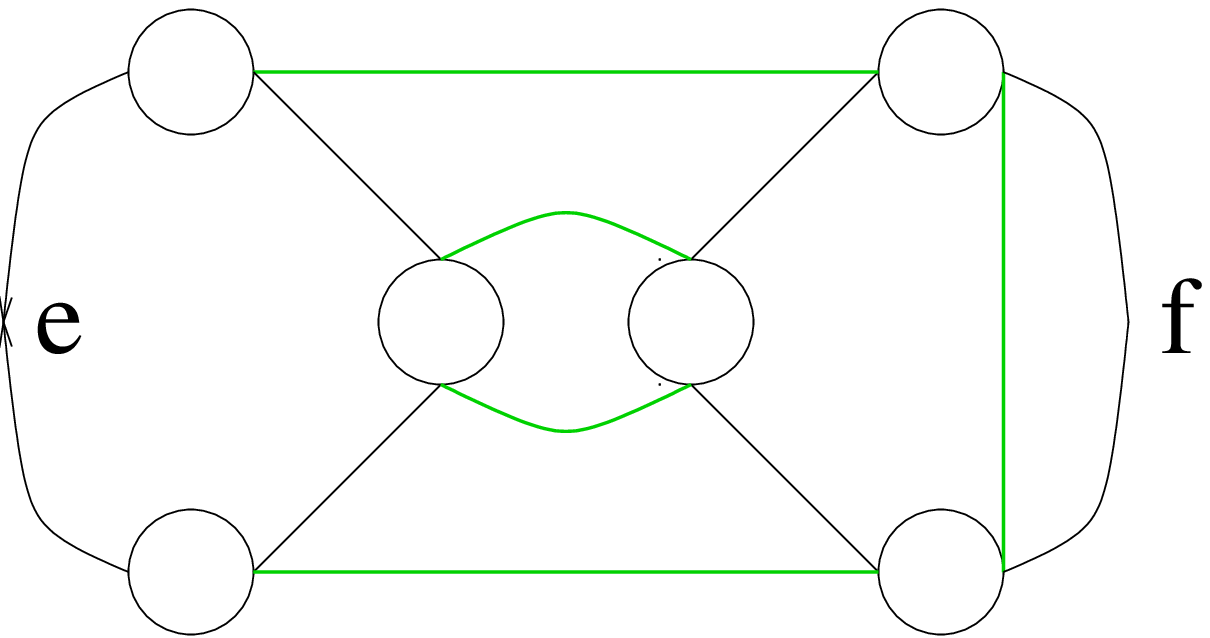}}\\
$S_L$ = A matching of $\Gamma_1 /e$ and $\Gamma_2$.
\end{center}
\end{minipage}
\]
Similarly,
a state $S_R$ in the RHS is defined to be a pair $(e,m)$
of an ideal edge $e \in I(\Gamma_2)$
and a matching $m:\Gamma_1 \to \Gamma_2 \backslash e$.
Then by definition the
RHS = $M(n)(\Gamma_1,\delta_E(\Gamma_2))=\sum_{S_R} w(S_R), \ 
\text{where}$
$w(S_R) = \sign(m) (2n)^{c(m)}$.

We now indicate the bijection between the two state spaces
by continuing our example.
We draw the state $S_R$ 
that corresponds to the state $S_L$ shown in the figure above.
\[
\begin{minipage}[b]{2.4 in}
\begin{center}
\psfrag{e}{\Huge $e$}
\resizebox{1 in}{.8 in}{\includegraphics{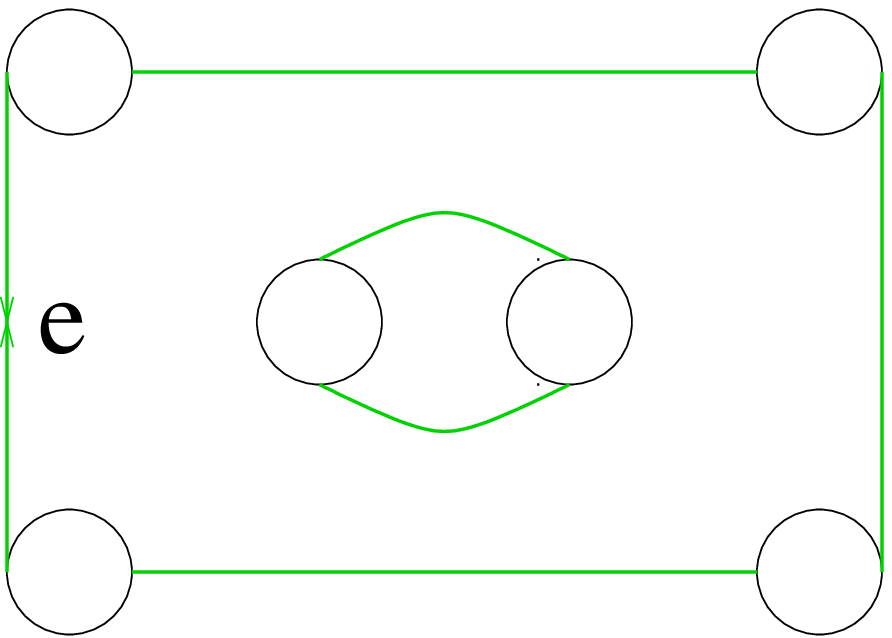}}\\
$\Gamma_2$ with an ideal edge $e$.
\end{center}
\end{minipage}
\begin{minipage}[b]{2.4 in}
\begin{center}
\psfrag{e}{\Huge $e$}
\resizebox{1 in}{1 in}{\includegraphics{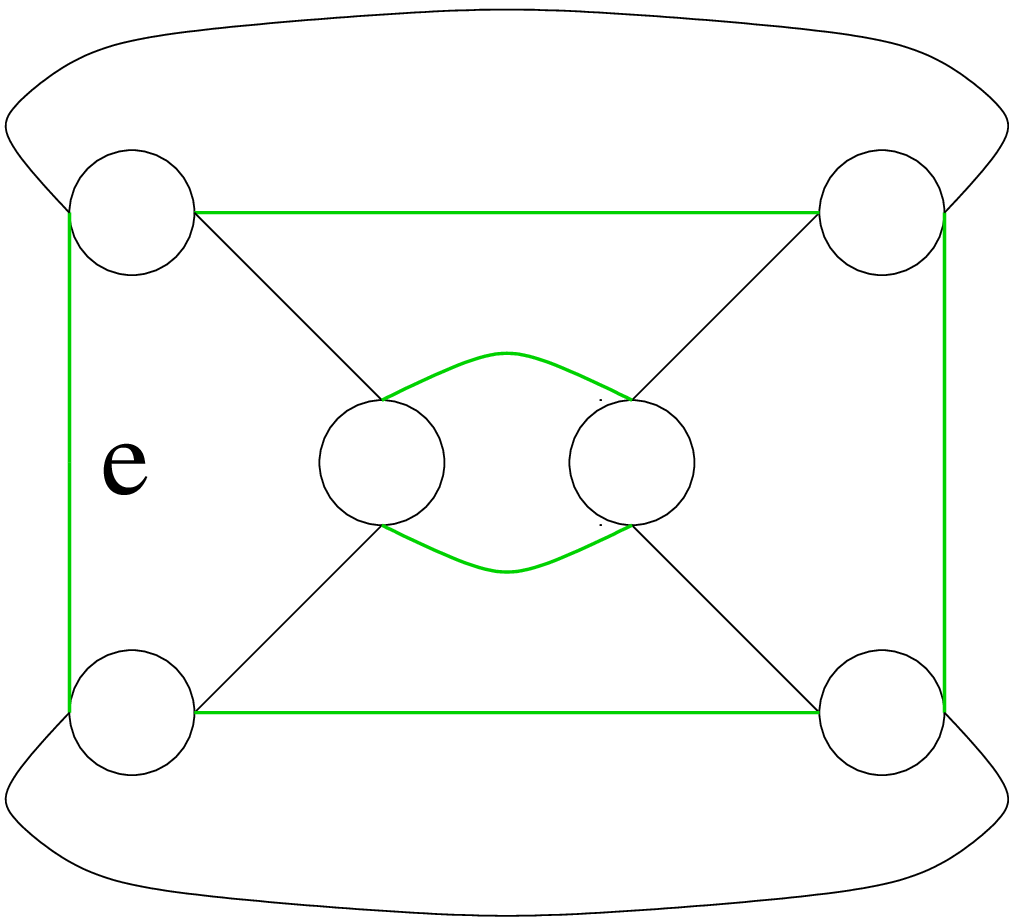}}\\
$S_R$ = A matching of $\Gamma_1$ and $\Gamma_2 \backslash e$.
\end{center}
\end{minipage}
\]
The two pictures that represent 
the states $S_L$ and $S_R$
are identical except for the local behaviour at the edge $e$.

In general, to see that $w(S_L) = w(S_R)$,
note that a state $S_R$ is of two kinds 
depending on whether the ideal edge $e \in I(\Gamma_2)$
is overlayed on a quasi-edge
$e \in E(\Gamma_1)$ or $e \in E_q(\Gamma_1) \setminus E(\Gamma_1).$
In the first case, $S_R$ has one more component than $S_L$
while in the second they are the same.
We illustrated the second case in our pictures.
This fits in with the two cases for $w(S_L)$.
It can also be checked that the signs work out correctly.
\end{proof}

\subsection{Non-degeneracy of the pairing}  \label{ss:pnd}

Let $\mathcal G_k^e$ be the span of all oriented graphs
with $k$ vertices and $e$ edges.
Observe that the subspaces $\mathcal G_k^e$ of $\mathcal G$,
as $k$ and $e$ vary,
are mutually orthogonal with respect to the pairing $M(n)$.
By our earlier notation (\ref{ss:is}),
we have
$\mathcal G_k^e=\mathcal G_k^{(r)}$,
with $r=1- \chi=e-k+1$.
The notation $\mathcal G_k^e$ is local to this section
and is introduced to improve clarity.

For any $k$ and $e$, 
the pairing $M(n)$ defines a map $\mathcal G_k^e \to (\mathcal G_k^e)^*$,
which we write $M_k^e(n)$.
Regard $M_k^e(n)$
as a matrix whose rows and columns are indexed by 
(isomorphism classes of) oriented graphs $\Gamma$
with $k$ vertices and $e$ edges and
where the $(\Gamma_i,\Gamma_j)$ entry is given by
$<M_k^e(n)(\Gamma_i),\Gamma_j>$.
By the definition of the pairing $M(n)$, 
the entries of this matrix 
are polynomials in $2n$ of degree $\leq e$
and moreover the degree $e$ entries
are precisely the diagonal entries.
Furthermore, the coefficient of $(2n)^e$ of
$<M_k^e(n)(\Gamma),\Gamma>$ is just $|\Aut(\Gamma)|$.
This is clear from the remark at the end of
(\ref{ss:pg}).

With $e$ and $k$ fixed and $n$ large enough,
the diagonal entries dominate,
so the matrix $M_k^e(n)$ is invertible, and
$\mathcal G_k^e \isoto (\mathcal G_k^e)^*$.
In other words,
the pairing $M(n)$ is non-degenerate
for $e$ and $k$ fixed and $n$ large enough.

\subsection{The deformation map $D(n)$} \label{ss:dmg}

The adjoint property of $\partial_n$ and $\delta_E$
established in Proposition~\ref{p:ad} says that
\begin{equation}  \label{e:ad}
M_{k-1}^{e-1}(n)
\circ \partial_n = 
\delta_E^* \circ M_k^e(n).
\end{equation}
Let $A_k^e$ denote the diagonal matrix with $(\Gamma,\Gamma)$
entry equal to $|\Aut(\Gamma)|$.
Then the degree $e$ part of the above equation is just
\begin{equation}  \label{e:sad}
A_{k-1}^{e-1}
\circ \partial_E = 
\delta_E^* \circ A_k^e.
\end{equation}
Equation~\eqref{e:ad} (resp.~\eqref{e:sad}) shows the precise sense in which
$\partial_n$ (resp $\partial_E$) 
is dual to $\delta_E$.
The results for $n$ large enough can be summarised 
in the following commutative diagram.
\begin{equation*}
\begin{CD}
\mathcal G_k^e @>M_k^e(n)>\iso> (\mathcal G_k^e)^* @<A_k^e<\iso< \mathcal G_k^e \\ 
@V{2n \partial_E + \partial_H}VV @V{\delta_E^*}VV @V{\partial_E}VV\\
\mathcal G_{k-1}^{e-1} @>\iso>M_{k-1}^{e-1}(n)> (\mathcal G_{k-1}^{e-1})^* @<\iso<A_{k-1}^{e-1}< \mathcal G_{k-1}^{e-1}.\\ 
\end{CD}
\end{equation*}
Now switching back to our earlier notation, 
the above commutative diagram can be rewritten as
\begin{equation}  \label{e:cdg}
\begin{CD}
\mathcal G_k^{(r)} @>M(n)>\iso> (\mathcal G_k^{(r)})^* @<A<\iso< \mathcal G_k^{(r)} \\ 
@V{2n \partial_E + \partial_H}VV @V{\delta_E^*}VV @V{\partial_E}VV\\
\mathcal G_{k-1}^{(r)} @>\iso>M(n)> (\mathcal G_{k-1}^{(r)})^* @<\iso<A< \mathcal G_{k-1}^{(r)}.\\ 
\end{CD}
\end{equation}
We point out that $M(n)$ is used to denote 
both the pairing and the map it induces.

\begin{corollary}
Let $k$ and $r$ be fixed. Then for $n$ large enough
$$H_k(\mathcal G^{(r)},\partial_n) = H_k(\mathcal G^{(r)},\partial_E).$$
\end{corollary}
\noindent
We will call the map $G^{(r)} \to G^{(r)}$,
which induces the above isomorphism,
the deformation map $D(n)$.
It is the map $M(n)$ followed by $A^{-1}$.
We note that for a graph $\Gamma$,
its image $D(n) \Gamma$
is a polynomial in $2n$ of degree $e$
with coefficients in $\mathcal G$.
Furthermore, the coefficient of $(2n)^e$ is exactly $\Gamma$.
Here $e=|E(\Gamma)|$, the number of edges in $\Gamma$.
The deformation map is a key step in the proof of 
the main theorem.
We will also study it briefly in Appendix~\ref{ss:dg}.

\section{The main theorem} \label{s:mt} 

In this section,
we give a precise statement of the main theorem.
We will also discuss the classical case briefly.
First recall some facts from Section~\ref{s:cm}.

Let $P$ be a reversible operad and 
$Q$ be its mated species.
There are two algebraic objects $PA$ and $QA$
associated to $P$ and $Q$ respectively.
One thinks of $QA$ as 
``Hamiltonian functions on the symplectic $P$-manifold''.
It has the structure of a Lie algebra.
The Lie algebra $QA$
depends on the dimension of the symplectic $P$-manifold.
So we write $QA_n$ when the dimension is $2n$.
The dimension is necessarily even.
We then have a sequence of Lie algebra inclusions
$$QA_1 \subset \ldots \subset QA_n \subset QA_{n+1} \subset \ldots$$
We denote the direct limit by $QA_{\infty}$.
And $QA_n$ always contains the symplectic Lie algebra
$\fsp$ as an anti-subalgebra.
We write
${{{\mathfrak{sp}}(2 \infty)}}$
for the corresponding direct limit.
The rational homology of $QA_{\infty}$, 
which we denote $H_*(QA_{\infty})$,
has the structure of a Hopf algebra.
We write $PH_*(QA_{\infty})$
for the subspace of primitive elements.

\subsection{Statement of the main theorem}
Kontsevich's result can now be stated as

\begin{theorem}  \label{t:mt}
$PH_*(QA_{\infty}) = 
H_*(Q\mathcal G,\partial_E) \oplus
PH_*({{{\mathfrak{sp}}(2 \infty)}}).$
\end{theorem}

\noindent
The term $H_*(Q\mathcal G,\partial_E)$
is the graph homology of the mated species $Q$.
It is the homology of the chain complex
$(Q\mathcal G,\partial_E)$ of graphs
that was defined in Section~\ref{s:gh}.

There are two conditions that we require in the theorem.
We assume that the operad $P$ is based on sets
rather than vector spaces.
This is because the proof involves a pairing on graphs
(\ref{ss:pg}),
which we know how to define only in the former case.
We hope that eventually this restriction
would not be necessary.
Hence, as of now,
the theorem cannot be applied to the Lie operad,
which is one of the cases claimed in~\cite{\god}.

Secondly, we assume that $P[1]$ 
is a singleton consisting of the unit element $u$,
see (\ref{ss:o}).
This is done mainly for simplicity.
If we drop this assumption 
then the summand $PH_*({{{\mathfrak{sp}}(2 \infty)}})$
has to be replaced by the homology of the bivalent graph complex
$(\mathcal B,\partial_E)$ defined in Section~\ref{s:gh}.

\begin{remark}
The primitive homology of ${{{\mathfrak{sp}}(2 \infty)}}$
is known (\ref{sss:ctg}).
It will be computed in the course of proving the theorem.
\end{remark}

\subsection{The classical case} \label{ss:ecc}

This is the commutative case, $P=c$ and $Q=cc$,
see (\ref{ss:tcc}).
In this case, the ``symplectic operad manifold''
actually exists and is simply $(\R^{2n},\omega_0)$.
The Lie algebra $QA_n$,
which we write as $ccA_n$, consists of 
polynomial functions in $2n$ variables 
with no constant or linear terms. 
The Lie structure is given by
the usual Poisson bracket 
(equation~\eqref{e:pb} in Section~\ref{s:sog}).
There is an evident subalgebra $uuA_n$
consisting of all homogeneous polynomials
of degree 2.
The reason for this notation is that 
$uuA_n$ can also be seen as an example of the theory
with $P=u$ and $Q=uu$, see (\ref{ss:sl}).
The Lie algebra $uuA_n$ is anti-isomorphic to 
the symplectic Lie algebra $\fsp$.

It is also clear that we have a sequence of Lie algebra inclusions
$$ccA_1 \subset \ldots \subset ccA_n \subset ccA_{n+1} \subset \ldots$$
The direct limit $ccA_{\infty}$ consists of finite polynomials
in infinitely many variables 
$p_1,p_2,\ldots,q_1,q_2,\ldots$.

\section{Proof of the main theorem-Part I}  \label{s:fin}

We now begin the proof of the main theorem.
It will be done in three steps.
In this section, we take the first step of relating 
the homology of a Lie algebra to graph homology.
We will prove the following theorem.

\begin{theorem}  \label{t:fin}
$H_*(QA_n) = H_*(\mathcal G,\partial_n).$
\end{theorem}

\noindent
For definitions of the above terms, see
(\ref{ss:fa}-\ref{ss:mate}) and (\ref{ss:gc}-\ref{ss:finb}).
For the commutative case, 
the Lie algebra $QA_n$ is easy to define, see (\ref{ss:ecc}).
The main ideas of the proof are already present 
in the commutative case.
Hence the reader may specialise to this case
on a first reading.

\begin{corollary}
$H_*(\fsp) = H_*(\mathcal F,\partial_n).$
\end{corollary}
\begin{proof}
We apply the theorem to the unit species, i.e. $Q=uu$.
In this case, $QA_n=uuA_n$, which is anti-isomorphic to $\fsp$,
see (\ref{ss:sl}).
Due to the trivial nature of the species,
the graphs have only fake (bivalent) vertices.
Hence the chain complex $\mathcal G$,
in this case, 
is simply the fake bivalent complex $\mathcal F$,
see (\ref{ss:cc}). 
\end{proof}

We now start the proof of Theorem~\ref{t:fin}.
It is best summarised in the following Kontsevich sentence.
\begin{quotation}
The spirit of the (quite simple) computations
is somewhere between Gelfand-Fuks computations
(see~\cite{\gelfand} and~\cite{\fuks}) and cyclic homology.
\end{quotation}

\subsection{Lie algebra homology}  \label{ss:lah}

A good introduction to Lie algebra homology
can be found in Weibel's book~\cite[Chapter 7]{\weibel}.
Recall that the homology of the Lie algebra $QA_n$ 
can be computed using
the Chevalley-Eilenberg or standard complex
$$\ldots\longrightarrow \mathcal C_{k+1}\buildrel \partial_{k+1}\over \longrightarrow 
\mathcal C_{k}\buildrel 
\partial_{k}\over \longrightarrow \mathcal C_{k-1}\buildrel \partial_{k-1}\over
\longrightarrow \ldots $$ with
$\mathcal C_k=\Lambda^k(QA_n)$ 
and 
\begin{equation}  \label{e:lb}
\partial_k(F_1 \wedge \ldots \wedge F_k)
=\sum_{1\leq s<t\leq k}(-1)^{s+t-1}\{F_s,F_t\}\wedge F_1\ldots \hat 
F_s\ldots \hat
F_t\ldots
\wedge F_k.
\end {equation}
The boundary operator commutes with the action
of the Lie algebra $QA_n$ 
on the exterior powers $\Lambda^k(QA_n)$. 

We have learnt to think of the Lie bracket $\{\ ,\ \}$ on $QA_n$
as a mating.
For the commutative case, this was discussed in
(\ref{ss:tcc}).
In the formula~\eqref{e:lb},
we take $k$ elements of $QA_n$,
say $F_1,\dots,F_k$, and to apply the boundary map
we do pairwise matings.
Now let us think for a moment about graphs.
The boundary operators on a graph are all defined
using quasi-edge contractions,
which are again matings.
Hence to relate the two notions,
all one needs to do is imagine the elements
$F_1,\dots,F_k$ as being the vertices of a graph.
We now work towards making this precise.

\subsection{Passing to the subcomplex of $\fsp$ invariants}  \label{ss:in}

Recall from (\ref{ss:fa}) that
\begin{equation}  \label{e:fsa}
QA_n = \bigoplus_{j\geq 2} (Q[j] \otimes V^{\otimes j})_{\Sigma_j}
= \bigoplus_{j\geq 2} QA_n^j.
\end{equation}
In the commutative case, $QA_n^j = S^j(V)$
is the $j$th symmetric tensor power of $V$.
This is just the space of commuting polynomials of degree $j$
in a basis of $V$.
 
The space $QA_n^j$ is a left $\fsp$ module with the action 
induced by the usual action on $V$ and
the trivial action on $Q[j]$, see (\ref{ss:sl}).
Thus, the Lie algebra $QA_n$
is a direct sum of finite dimensional $\fsp$ modules.
Since the Lie algebra $\fsp$ is simple,
it follows that $QA_n$ and $\Lambda^k(QA_n) = \mathcal C_k$ 
are semisimple $\fsp$ modules.
This allows us to write
$\mathcal C_{k} = (\mathcal C_{k})^{\fsp} \oplus \fsp \cdot \mathcal C_k$.
A standard argument now shows that
the subcomplex $\fsp \cdot \mathcal C$ is exact.
Hence
the homology of the standard complex is the same as the 
homology of the subcomplex of $\fsp$-invariants:
$$\ldots\longrightarrow (\mathcal C_{k+1})^{\fsp}\buildrel \partial_{k+1}\over \longrightarrow 
(\mathcal C_{k})^{\fsp}{ \buildrel \partial_{k}\over \longrightarrow} (\mathcal C_{k-1})^{\fsp}\buildrel
\partial_{k-1}\over
\longrightarrow \ldots\qquad $$

\begin{remark}
The standard argument is as under.
Let $\mathfrak g$ be a Lie algebra and
$\mathcal C$ be its standard complex.
Then for $\xi \in \mathfrak g$,
define Lie derivative $L_{\xi}$ (action of $\xi$)
and contraction operator $i_{\xi}$ (wedge with $\xi$)
as operators on $\mathcal C$.
The relations of Lemma~\ref{l:rel} (Section~\ref{s:st})
hold in this case also.
Cartan's formula implies that the Lie derivative
is zero on $H_*(\mathcal C)$.
Now let $\mathfrak h$ be a simple subalgebra of $\mathfrak g$
such that $\mathcal C$ is a semisimple $\mathfrak h$ module.
Then the subcomplex
$\mathfrak h \cdot \mathcal C$ is exact.

For the semisimplicity of $\Lambda^k(QA_n)$,
we used two facts.
Any finite dimensional module of a simple Lie algebra is semisimple.
The tensor product of two finite dimensional semisimple modules 
over any Lie algebra
is again semisimple~\cite[pg 83]{\jacobson}.
\end{remark}

\subsection{Passing from invariants to oriented graphs}  \label{ss:inv}

We want to relate the above subcomplex of $\fsp$-invariants
to the graph complex $(\mathcal G,\partial_n)$.
To this end, we first try to get a better understanding
of the chain groups $(\wedge^k QA_n)^{\fsp}$.

The description of $QA_n$ given by equation~\eqref{e:fsa}
gives us the following formula:
\[
\wedge^k QA_n = \bigoplus_{\tiny \begin{array}{c}k_1+k_2+\ldots+k_r=k\\
                             2\leq j_1<j_2<\ldots <j_r\end{array}}
\left(\wedge^{k_1}QA^{j_1} \otimes \cdots \otimes
\wedge^{k_r}QA^{j_r}\right).
\]
We try to understand the $\fsp$ invariants in each summand.
This will be done in three stages.
As a concrete example, we will consider the summand
$\wedge^{1}QA^{2} \otimes \wedge^{2}QA^{3} \otimes 
\wedge^{1}QA^{4}$
and see how the analysis works on it at every stage.
We will illustrate it with the tree species $(Q=tt)$.

The term $\wedge^{k_1}QA^{j_1}$ is a quotient of the tensor power
$\left(Q[j_1] \otimes V^{\otimes j_1}\right)^{\otimes k_1}$.
Hence
each summand of $\wedge^kQA_n$ is a quotient of the tensor power
$\otimes_{t=1}^r (Q[j_t] \otimes V^{\otimes j_t})^{\otimes k_t}$. 
In the first two stages,
we will figure out the $\fsp$-invariants in this tensor power.
In the third stage,
we will mod out by appropriate actions of the symmetric groups.

\subsubsection{The first stage}
We find the invariants in
$V^{\otimes \sum_{t=1}^r j_tk_t}$. 
By the invariant theory of $\fsp$, we know that
a base for the invariants in $V^{\otimes \sum_{t=1}^r k_tj_t}$ 
is given by
oriented chord diagrams on $\sum_{t=1}^r k_tj_t$ vertices,
if $n$ is sufficiently large;
see~\cite[Appendix F]{\fh}. 

For the definition of a chord diagram, see~(\ref{sss:k}).
By an oriented chord diagram,
we mean that each chord is oriented 
and reversing the orientation of a single chord
incurs a minus sign.
\[
\begin{minipage}{\linewidth}
\begin{center}
\psfrag{1}{\Huge $1$}
\psfrag{2}{\Huge $2$}
\psfrag{3}{\Huge $3$}
\psfrag{4}{\Huge $4$}
\psfrag{5}{\Huge $5$}
\psfrag{6}{\Huge $6$}
\psfrag{7}{\Huge $7$}
\psfrag{8}{\Huge $8$}
\psfrag{9}{\Huge $9$}
\psfrag{10}{\Huge $10$}
\psfrag{11}{\Huge $11$}
\psfrag{12}{\Huge $12$}
\resizebox{1.5 in}{1.3 in}{\includegraphics{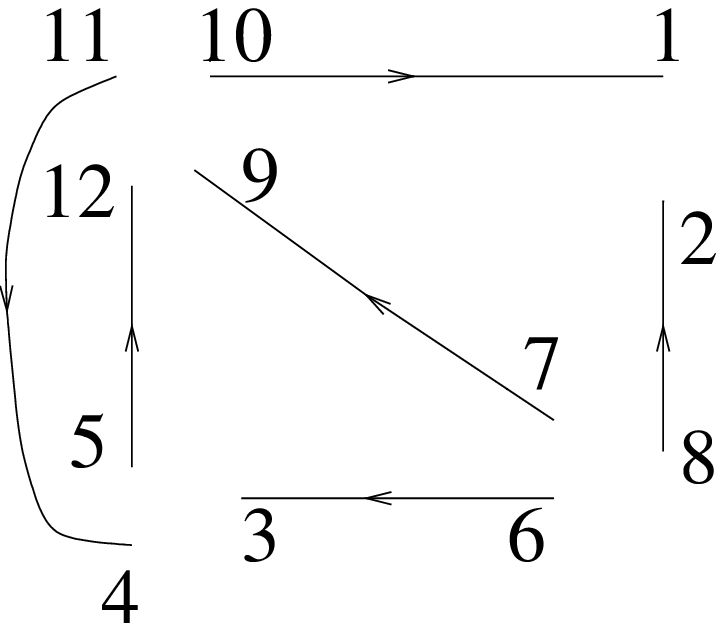}}\\
An oriented chord diagram on $12$ vertices.
It gives a $\fsp$ invariant in
$V^{\otimes 12} =
(V^{\otimes 2})^{\otimes 1} \otimes
(V^{\otimes 3})^{\otimes 2} \otimes
(V^{\otimes 4})^{\otimes 1}$.
\end{center}
\end{minipage}
\]
\subsubsection{On how an oriented chord diagram gives an invariant}
\label{sss:ci}
Each vertex of the diagram represents a tensor factor,
in the order given by the vertex labelling.
For each edge, we put a $p_i$ at the tail of the arrow
and a $q_i$ at the head or 
we put a $q_i$ at the tail of the arrow
and a $p_i$ at the head,
incurring a minus sign as a result.
We then sum over all possible choices to get the invariant.

The smallest $\fsp$ invariant lies in $V \otimes V$.
It is given by 
$\sum_{i=1}^n p_i \otimes q_i - q_i \otimes p_i$,
which we represent by the chord diagram
$1 \pointsat 2$.

\subsubsection{The second stage}
To get the invariants in 
$\otimes_{t=1}^r (Q[j_t] \otimes V^{\otimes j_t})^{\otimes k_t}$, 
we tensor the space of invariants obtained above by
$\otimes_{t=1}^r Q[j_t]^{\otimes k_t}$. 
This is alright because the $\fsp$ action is trivial 
on the $Q[j]$'s. 
This has the following effect on our picture.
\[
\begin{minipage}{1.8 in}
\begin{center}
\psfrag{1}{\Huge $1$}
\psfrag{2}{\Huge $2$}
\psfrag{3}{\Huge $3$}
\psfrag{4}{\Huge $4$}
\psfrag{5}{\Huge $5$}
\psfrag{6}{\Huge $6$}
\psfrag{7}{\Huge $7$}
\psfrag{8}{\Huge $8$}
\psfrag{9}{\Huge $9$}
\psfrag{10}{\Huge $10$}
\psfrag{11}{\Huge $11$}
\psfrag{12}{\Huge $12$}
\resizebox{1.5 in}{1.3 in}{\includegraphics{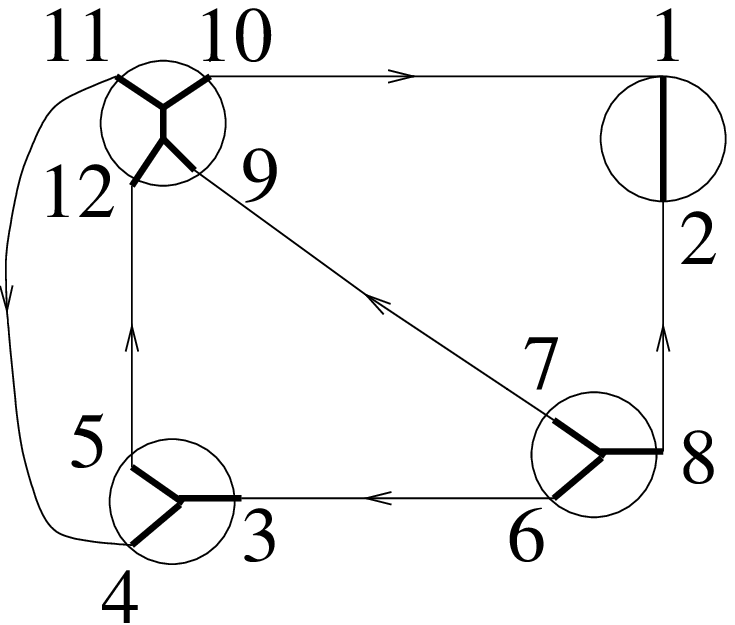}}
\end{center}
\end{minipage}
\]

\begin{remark}
We have drawn trees inside the circles because
we are illustrating with the tree species.
In the commutative case,
we will instead get ordinary graphs
with directed edges and an order on the set of half-edges.
And in the associative case, we get ribbon graphs
with similar data.
\end{remark}

\subsubsection{The third stage}  

Moding out the actions of the symmetric groups
has the following effect on our picture.
Again, we are illustrating with the tree case.
\[
\begin{minipage}{1.8 in}
\begin{center}
\psfrag{1}{\Huge $1$}
\psfrag{2}{\Huge $2$}
\psfrag{3}{\Huge $3$}
\psfrag{4}{\Huge $4$}
\resizebox{1.5 in}{1.3 in}{\includegraphics{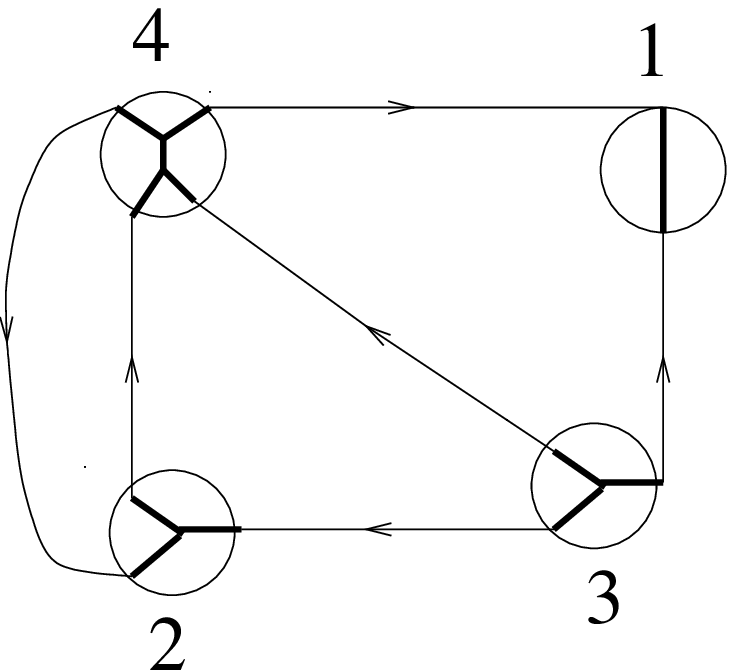}}
\end{center}
\end{minipage}
\]
\noindent
Firstly, we have removed all the labels on the half-edges
and instead given an ordering to the vertices.
This is because we are moding out the action of the $\Sigma_j$'s.
Secondly, we must now interpret the order on the vertices
of the graph in the sense of orientation.
That is, if we interchange the order of two consecutive vertices
then we pick a minus sign.
This is due to the presence of the wedges. 
What we are left with is precisely an oriented $Q$-graph,
see (\ref{ss:og}).
Note that the valence of the vertices is at least 2,
since the grading on $QA$ begins at 2.
Thus we get an isomorphism of chain groups
$(\mathcal C_k)^{\fsp} \cong \mathcal G_k$,
see (\ref{ss:gc}).

\subsection{Comparison of the boundary maps of the two complexes}  
\label{ss:cb}

To complete the proof of Theorem~\ref{t:fin},
we need to show that the following diagram commutes.
\begin{equation}  \label{e:comp}
\begin{CD}
\mathcal G_k @>\iso>> (\mathcal C_k)^{\fsp}\\
@V{2n \partial_E + \partial_H}VV @V{\partial}VV\\
\mathcal G_{k-1} @>\iso>> (\mathcal C_{k-1})^{\fsp}.\\
\end{CD}
\end{equation}
Define 
$\partial^{\prime} : \mathcal G_k \rightarrow \mathcal G_{k-1}$
as the composite of three maps in the above diagram
(the middle map is $\partial$).
We want to show that
$\partial^{\prime} = 2n \partial_E + \partial_H$.
For this, we first understand the map
$\mathcal G_k \rightarrow (\mathcal C_k)^{\fsp}$ better.

\subsubsection{On how an oriented graph gives an invariant}  \label{sss:gi}

Starting with an oriented graph with $k$ vertices,
we want to construct a $\fsp$ invariant in $\mathcal C_k=\Lambda^k(QA_n)$.
The description that we give follows
directly from the one that we gave 
for an oriented chord diagram (\ref{sss:ci}).

Let $(\Gamma,\sigma)$ be an oriented graph.
Choose a representative for $\sigma$.
This means that the vertices of $\Gamma$ are ordered 
and the edges are oriented.
Each vertex of the graph represents a tensor factor,
in the order given by the vertex labelling.
For each edge, we put a $p_i$ at the tail of the arrow
and a $q_i$ at the head or 
we put a $q_i$ at the tail of the arrow
and a $p_i$ at the head,
but picking a minus sign.
This is called a state of the edge.
And a \emph{state} of the graph is a choice of a state for every edge.
Summing over all states of $\Gamma$ and
passing to the wedge product
gives an $\fsp$ invariant in $\Lambda^*(QA_n)$.
\[
\begin{minipage}{4 in}
\begin{center}
\psfrag{p1}{\Huge $p_1$}
\psfrag{p2}{\Huge $p_2$}
\psfrag{p3}{\Huge $p_3$}
\psfrag{q1}{\Huge $q_1$}
\psfrag{q2}{\Huge $q_2$}
\psfrag{q3}{\Huge $q_3$}
\resizebox{1.1 in}{.7 in}{\includegraphics{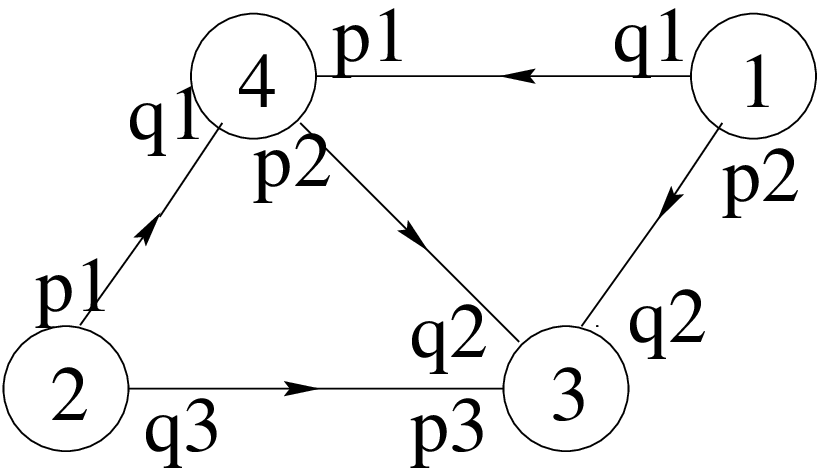}}
\end{center}
\end{minipage}
\]

\noindent
An example of a state is shown above.
It gives us the term
\[
\begin{minipage}{\linewidth}
\begin{center}
\begin{minipage}{1 in}
\begin{center}
\psfrag{p1}{\Huge $q_1$}
\psfrag{p2}{\Huge $p_2$}
\psfrag{q1}{\Huge $q_1$}
\psfrag{q2}{\Huge $p_2$}
\psfrag{1}{}
\psfrag{2}{}
\psfrag{uu}{}
\resizebox{.7 in}{.15 in}{\includegraphics{sf.eps}}
\end{center}
\end{minipage}
$\bigwedge$
\begin{minipage}{1 in}
\begin{center}
\psfrag{p1}{\Huge $p_1$}
\psfrag{p2}{\Huge $p_2$}
\psfrag{q1}{\Huge $q_1$}
\psfrag{q2}{\Huge $q_3$}
\psfrag{1}{}
\psfrag{2}{}
\psfrag{uu}{}
\resizebox{.7 in}{.15 in}{\includegraphics{sf.eps}}
\end{center}
\end{minipage}
$\bigwedge$
\begin{minipage}{1 in}
\begin{center}
\psfrag{p1}{\Huge $p_3$}
\psfrag{p2}{\Huge $q_2$}
\psfrag{q1}{\Huge $q_2$}
\resizebox{.8 in}{.6 in}{\includegraphics{qas.eps}}
\end{center}
\end{minipage}
$\bigwedge$
\begin{minipage}{1 in}
\begin{center}
\psfrag{p1}{\Huge $p_1$}
\psfrag{p2}{\Huge $p_2$}
\psfrag{q1}{\Huge $q_1$}
\resizebox{.8 in}{.6 in}{\includegraphics{qas.eps}}
\end{center}
\end{minipage}\\
\end{center}
\end{minipage}
\]
in $\Lambda^4(QA_n)$.
Since there are an even number of negative signs,
the net sign is positive.
To get the invariant, we sum over all states.

\subsubsection{Comparing the matings}  

Now we show that
$\partial^{\prime} = 2n \partial_E + \partial_H$.
Let $(\Gamma,\sigma) \in \mathcal G_k$. 
And let $I$ be the invariant in $(\mathcal C_k)^{\fsp}$
corresponding to $\Gamma$
obtained by the procedure described above.

In order to compute $\partial(I)$,
we apply the formula for $\partial$ given by equation \eqref{e:lb}
to each state of $I$.
So for each state,
we must compute the Lie bracket of the data
at each pair of vertices.
However, recall from (\ref{ss:mate})
that the Lie bracket is simply a sum of matings.
Hence $\partial(I)$ can be computed as a sum
over all states of $I$,
where for each state,
we mate a $p_i$ on a half-edge at a vertex
with a $q_i$ on another half-edge at a different vertex
in all possible ways.
Another way to say this is that
a mating occurs along a quasi-edge
(which is a pair of distinct half-edges)
that is not a quasi-loop (\ref{ss:finb}).

Now we rewrite $\partial(I)$ as a sum indexed by quasi-edges $e$
of $\Gamma$.
The summand for a quasi-edge $e$ 
is the sum of those terms for which
the mating occurs along $e$.
It is obtained in two steps.

1. Fix a state for $e$. Then sum over all possible states 
of the other edges.

2. Do this for each of the $2n$ states of $e$.

\noindent
We explain a little of how the orientation will work out.
We fix a representative $\sigma$ for $\Gamma$.
Let $v_1$ and $v_2$ with labels $s < t$ 
be the vertices of the quasi-edge $e$ 
that we want to contract.
In order to get a representative for $\sigma /e$,
we first reorder the vertices so that
$v_1$ and $v_2$ have labels $1$ and $2$ respectively.
In doing so, we pick up a factor of $(-1)^{s+t-1}$. 
This is precisely the factor that appears in formula \eqref{e:lb}
for $\partial(I)$.

Now for a quasi-edge $e$,
we perform the two steps above.
The analysis splits into two cases.
In the first case, assume that $e$ is an actual edge,
i.e. $e \in E(\Gamma)$.
Then the first step itself gives an invariant that comes from 
$(\Gamma /e,\sigma /e).$
Hence, after the second step, 
the net contribution is $2n (\Gamma /e,\sigma /e).$

In the second case, assume that $e$ is a quasi-edge
that is not an edge,
i.e. $e \in E_q(\Gamma) \setminus E(\Gamma)$. 
So the first step does not give an invariant.
This is because, if the state of $e$ is fixed
then the state of its partner, the quasi-edge $f$,
is also fixed. 
This is shown in the figure below.
\[
\begin{minipage}{1.6 in}
\begin{center}
\psfrag{v1}{\Huge $v_1$}
\psfrag{v2}{\Huge $v_2$}
\psfrag{e1}{\Huge $$}
\psfrag{e2}{\Huge $$}
\psfrag{e}{\Huge $e$}
\psfrag{f}{\Huge $f$}
\psfrag{pi}{\Huge $p_i$}
\psfrag{qi}{\Huge $q_i$}
\resizebox{1.4 in}{.9 in}{\includegraphics{bqi.eps}}
\end{center}
\end{minipage}
\]
However, after the second step, 
we do get the invariant $(\Gamma /e,\sigma /e).$

Hence putting the two cases together, we obtain
$$\partial^{\prime}_{} (\Gamma,\sigma) = 2n \sum_{e \in E(\Gamma)} (\Gamma /e,
\sigma /e) + \sum_{e \in E_q(\Gamma) \setminus E(\Gamma)} 
(\Gamma /e,\sigma /e).$$ 
The right hand side is precisely
$2n \partial_E + \partial_H$.
This shows the commutativity of diagram~\eqref{e:comp} and
completes the proof of Theorem~\ref{t:fin}.

\section{Proof of the main theorem-Part II}  \label{s:sta}

In the previous section, the Lie algebra homology
$H_*(QA_n)$ was related to graph homology.
In this section, we show that the homology 
of the Lie algebra $QA_n$, as $n$ varies, is ``stable''.
Recall that 
$QA_\infty = \underset{n \to \infty}{\lim} QA_n$,
see (\ref{ss:ha}).
We will prove the following theorem.

\begin{theorem}  \label{t:sta}
$H_*(QA_{\infty}) = H_*(\mathcal G,\partial_E).$
\end{theorem}
\noindent
As in the previous section,
we apply the theorem to the unit species $uu$
and obtain the following corollary.

\begin{corollary}
$H_*({{{\mathfrak{sp}}(2 \infty)}}) = H_*(\mathcal F,\partial_E).$
\end{corollary}

\noindent
Let us try to prove this theorem. We have
$$
\begin{array}{r c l}
H_k(QA_{\infty})
 & = & \underset{n \to \infty}{\lim} H_k(QA_n). \\
 & = & \underset{n \to \infty}{\lim} H_k(\mathcal G,\partial_n). \\
 & = & \underset{r}{\bigoplus} \underset{n \to \infty}{\lim} H_k(\mathcal G^{(r)},\partial_n). 
\end{array}
$$
The first equality says that homology commutes with direct limits.
And the second equality is the content of Theorem~\ref{t:fin}.
The bonding maps
$H_k(\mathcal G,\partial_n) \to H_k(\mathcal G,\partial_{n+1})$
for the direct limit are defined using the isomorphism
$H_k(QA_n) \iso H_k(\mathcal G,\partial_n)$.
For the third equality, 
see the discussion in (\ref{ss:is}).
There is a subtle point here.
We need to know that the bonding maps restrict to
$H_k(\mathcal G^{(r)},\partial_n) \to H_k(\mathcal G^{(r)},\partial_{n+1})$.
Even assuming this, we are stuck.
To complete the argument, one needs the following proposition.

\begin{proposition}  \label{p:sta}
$\underset{n \to \infty}{\lim} H_k(\mathcal G^{(r)},\partial_n) = 
H_k(\mathcal G^{(r)},\partial_E)$.
\end{proposition}
\noindent
It will be proved using the ideas of Section \ref{s:mgh}. 
The proof will be completed by the end of this section.
The subtle point about the definition of the LHS
raised above will be dealt in the course of the proof.

\subsection{A pairing on $\Lambda^*(QA_n)$}  \label{ss:pl}

We have fixed a basis $p_1, \ldots,p_n,q_1, \ldots,q_n$
for $V$.
This gives us a basis for the Lie algebra $QA_n$
and its exterior algebra $\Lambda^*(QA_n)$.
To be explicit, a basis for $QA_n$ is given by distinct monomials.
And for $\Lambda^*(QA_n)$ is given by taking wedge products
of such monomials.
\[
\begin{minipage}{\linewidth}
\begin{center}
\begin{minipage}{1.2 in}
\begin{center}
\psfrag{x1}{\Huge $p_1$}
\psfrag{x2}{\Huge $p_2$}
\psfrag{x4}{\Huge $p_4$}
\resizebox{.8 in}{.6 in}{\includegraphics{qa.eps}}
\end{center}
\end{minipage}
$\bigwedge$
\begin{minipage}{1.2 in}
\begin{center}
\psfrag{p1}{\Huge $p_1$}
\psfrag{p2}{\Huge $p_2$}
\psfrag{q1}{\Huge $q_1$}
\psfrag{q2}{\Huge $q_2$}
\psfrag{1}{}
\psfrag{2}{}
\psfrag{uu}{}
\resizebox{.7 in}{.15 in}{\includegraphics{sf.eps}}
\end{center}
\end{minipage}
$\bigwedge$
\begin{minipage}{1.2 in}
\begin{center}
\psfrag{p1}{\Huge $p_1$}
\psfrag{p2}{\Huge $p_2$}
\psfrag{q1}{\Huge $q_1$}
\psfrag{q2}{\Huge $q_2$}
\resizebox{.8 in}{.6 in}{\includegraphics{qas.eps}}
\end{center}
\end{minipage}\\
\vspace{.1 in} A basis element of $\Lambda^3(QA_4)$.
\vspace{.1 in}
\end{center}
\end{minipage}
\]
Define a pairing
$M^{\prime}(n) : \Lambda^*(QA_n)\otimes \Lambda^*(QA_n)\to \Q$ by
$M^{\prime}(n)(\psi_i,\psi_j) = \delta_{ij}$,
where $\{\psi_i\}$ is the basis of the exterior algebra
that was chosen above.
The subspaces $\Lambda^k(QA_n)$, as $k$ varies,
are clearly orthogonal with respect to 
$M^{\prime}(n)$.

To be more explicit, let
$\varphi_1 = F_1 \wedge \ldots \wedge F_k$ and
$\varphi_2 = H_1 \wedge \ldots \wedge H_k$ be wedges of monomials.
Then
$M^{\prime}(n)(\varphi_1,\varphi_2)$ is nonzero 
only if $\varphi_1=\varphi_2$ or $\varphi_1=-\varphi_2$.
This happens if
$F_i = H_{\pi(i)}$ for a permutation $\pi$,
i.e. there is a matching of the $F_i$'s with the $H_j$'s.
This way of saying it brings us closer to the pairing $M(n)$,
that was defined on graphs (\ref{ss:pg}). We say this more precisely.

The isomorphism
$\mathcal G_k \isoto (\Lambda^k(QA_n))^{\fsp}$,
explained in (\ref{sss:gi}),
gives us an inclusion map 
$\mathcal G \into \Lambda^*(QA_n)$.
Under this inclusion,
the pairing $M^{\prime}(n)$ has two important restrictions, namely,
$M^{\prime}(n) : \Lambda^*(QA_n)\otimes \mathcal G \to \Q$ and 
$M^{\prime}(n) : \mathcal G \otimes \mathcal G \to \Q$.
We will see in Proposition \ref{p:pr}
that the second pairing is just the pairing
$M(n)$,
that was defined in (\ref{ss:pg}).

\subsection{Connection between the pairings $M^{\prime}(n)$ and $M(n)$} 
\label{ss:pr}

We will now see how formula \eqref{e:pg} for the pairing $M(n)$
emerges naturally by analysing $M^{\prime}(n)$.

\begin{proposition}  \label{p:pr}
The pairing $M^{\prime}(n) : \Lambda^*(QA_n)\otimes \Lambda^*(QA_n)\to \Q$ restricts to the pairing
$M(n) : \mathcal G \otimes \mathcal G \to \Q$ 
under the inclusion map $\mathcal G \into \Lambda^*(QA_n)$.
\end{proposition}
\begin{proof}
Let $(\Gamma_1,\sigma_1),(\Gamma_2,\sigma_2) \in \mathcal G$.
We will suppress the inclusion map and
orientations from the notation.
We want to show
$M^{\prime}(n)(\Gamma_1,\Gamma_2) = M(n)(\Gamma_1,\Gamma_2)$.
We compute the LHS by summing over all pairs of states of
$\Gamma_1$ and $\Gamma_2$.
Recall that a state of a graph is a choice
of $p_i \lhor q_i$ for every edge (\ref{sss:gi}).
We write, 
\[
M^{\prime}(n)(\Gamma_1,\Gamma_2) = \sum_{S_1,S_2} M^{\prime}(n)(S_1,S_2),
\]
where $S_i$ is a state of $\Gamma_i$.
A term in this sum is nonzero only if
one can match $S_1$ and $S_2$ in the following sense.

There is a matching $m:\Gamma_1 \to \Gamma_2$ 
(Definition \ref{d:gm})
with the following additional piece of data.
For each component in the matching
(Definition \ref{d:csm}), 
there is an index $i$ such that
the vertices of the polygon are alternately
labelled $p_i$ and $q_i$.
\[
\begin{minipage}{.7 in}
\begin{center}
\psfrag{a}{\Huge $q_i$}
\psfrag{b}{\Huge $p_i$}
\psfrag{c}{\Huge $q_i$}
\psfrag{d}{\Huge $p_i$}
\psfrag{e}{\Huge $q_i$}
\psfrag{f}{\Huge $p_i$}
\resizebox{.7 in}{.7 in}{\includegraphics{ch.eps}}
\end{center}
\end{minipage}
\hspace{1 in}
\begin{minipage}{.6 in}
\begin{center}
\psfrag{a}{\Huge $q_j$}
\psfrag{b}{\Huge $p_j$}
\psfrag{c}{\Huge $q_j$}
\psfrag{d}{\Huge $p_j$}
\resizebox{.5 in}{.5 in}{\includegraphics{cs.eps}}
\end{center}
\end{minipage}
\]
Hence we can group terms and sum over all
matchings $m:\Gamma_1 \to \Gamma_2$.
For each matching $m$, there are $c(m)$ number of polygons.
And for each polygon, there are $2n$ choices of the index $i$.
This gives a factor of $(2n)^{c(m)}$ 
with a sign which can be checked to be 
$\sign(m)$. 
Hence
$M^{\prime}(n)(\Gamma_1,\Gamma_2) = \underset{m:\Gamma_1 \to \Gamma_2}{\sum}
\sign(m) (2n)^{c(m)}$,
which by definition is
$M(n)(\Gamma_1,\Gamma_2)$.
This completes the proof.

\end{proof}

\subsection{A commutative diagram}  \label{ss:cd}

Consider the restriction of the pairing
$M^{\prime}(n) : \Lambda^k(QA_n) \otimes \mathcal G_k^{(r)} \to \Q$.
This then defines a map
$\Lambda^k(QA_n) \to (\mathcal G_k^{(r)})^*$,
which we again denote by $M^{\prime}(n)$.
From Proposition \ref{p:pr},
it is clear that this map is an extension of the map
$\mathcal G_k^{(r)} \to (\mathcal G_k^{(r)})^*$
defined using the pairing $M(n)$.
We now claim that the following diagram commutes.
\begin{equation*}
\begin{CD}
\mathcal G_k^{(r)} @>>> \Lambda^k(QA_n)@>M^{\prime}(n)>> (\mathcal G_k^{(r)})^* @<A<\iso< \mathcal G_k^{(r)} \\ 
@V{2n \partial_E + \partial_H}VV @V{\partial}VV @V{\delta_E^*}VV @V{\partial_E}VV\\
\mathcal G_{k-1}^{(r)} @>>> \Lambda^{k-1}(QA_n)@>>M^{\prime}(n)> (\mathcal G_{k-1}^{(r)})^* @<\iso<A< \mathcal G_{k-1}^{(r)}.\\ 
\end{CD}
\end{equation*}
This is a fattening of the diagram \eqref{e:cdg} 
at the end of Section \ref{s:mgh},
with the Chevalley-Eilenberg complex added in.

The commutativity of the first square follows from 
diagram \eqref{e:comp} in (\ref{ss:cb}).
Saying that the second square commutes 
is equivalent to saying that
the boundary map $\partial:\Lambda^*(QA_n)\to \Lambda^*(QA_n)$
and the coboundary map $\delta_E:\mathcal G \to \mathcal G$ are adjoints
with respect to the pairing $M^{\prime}(n)$.
This follows by generalising the proof of Proposition \ref{p:ad}.
This proves the claim.

\subsection{Invariance of the pairing}  

The pairing 
$M^{\prime}(n) : \Lambda^*(QA_n)\otimes \Lambda^*(QA_n)\to \Q$
is $\fsp$ invariant 
in a precise sense.
It is more convenient to express it using 
the Lie algebra $uuA_n$, see (\ref{ss:sl}), 
which is anti-isomorphic to $\fsp$.

\begin{proposition}  \label{p:in}
We have
$M^{\prime}(n)(\varphi_1,\varphi_2 \cdot H) + 
M^{\prime}(n)(\varphi_1 \cdot JH, \varphi_2) = 0$,
where $H \in uuA_n$ and $J$ 
is the element of the symplectic group defined by
$J(p_i)=q_i$ and $J(q_i)= - p_i$.
\end{proposition}

\begin{proof}
We may assume that $\varphi_1,\varphi_2$ 
are wedges of monomials, 
i.e. they are two basis elements (upto sign) of
$\Lambda^k(QA_n)$.
We may also assume that $H$ is a monomial in $uuA_n$.
Say, for definiteness, that $H=q_1 q_2$ and so
$JH = p_1 p_2$.
Note that $\varphi_2 \cdot H$ involves replacing an occurrence of
$p_2$ (resp. $p_1$) in $\varphi_2$ by $q_1$ (resp. $q_2$),
see (\ref{ss:sl}).
And $\varphi_1 \cdot JH$ has exactly the opposite effect.
It involves replacing an occurrence of
$q_1$ (resp. $q_2$) in $\varphi_1$ by
$p_2$ (resp. $p_1$) and picking a minus sign.
It is now fairly clear how terms 
from the two pairing expressions would cancel.
\end{proof}

\begin{corollary}
The subspaces $(\Lambda^k(QA_n))^{\fsp}$ and
$\fsp \cdot (\Lambda^k(QA_n))$ are orthogonal 
with respect to the pairing $M^{\prime}(n)$.
\end{corollary}

\subsection{Stability of the pairings}  

Consider the Lie algebra $QA_{n+1}$.
The underlying vector space $V_{n+1}$ has dimension $2(n+1)$.
The basis of $V_{n+1}$ is an extension of the one on $V$ to
$p_1, \ldots,p_{n+1},q_1, \ldots,q_{n+1}$.
Hence the basis for $QA_n$ and $\Lambda^*(QA_n)$
includes into the basis for $QA_{n+1}$ and $\Lambda^*(QA_{n+1})$
respectively.

We now show that the restricted pairings 
$M^{\prime}(n) : \Lambda^*(QA_n)\otimes \mathcal G \to \Q$ 
are stable, that is,
for $\varphi \in \Lambda^*(QA_n)$ and $\Gamma \in \mathcal G$,
\begin{equation}  \label{e:ps}
M^{\prime}(n)(\varphi,\Gamma) = M^{\prime}(n+1)(\varphi,\Gamma). 
\end{equation}
The $\varphi$ in the RHS is obtained by using the inclusion of
$\Lambda^*(QA_n)$ in $\Lambda^*(QA_{n+1})$.

Express both sides as sums over the states of $\Gamma$,
see (\ref{sss:gi}).
Though the graph $\Gamma$ is the same on both sides,
it has more states for ``$n+1$'' than ``$n$''
because there are more variables available.
The extra states of $\Gamma$ in the RHS are the ones that
involve at least one $p_{n+1}$ or $q_{n+1}$.
However, they do not contribute anything to the pairing
$M^{\prime}(n+1)$, since by assumption, 
$\varphi$ does not involve either
$p_{n+1}$ or $q_{n+1}$. 
This proves that equation \eqref{e:ps} holds.

\subsection{The stability commutative diagram}  

The chain group 
$\mathcal C_{k} = \Lambda^k(QA_n)$ can be written as
$\mathcal C_{k} = (\mathcal C_{k})^{\fsp} \oplus \fsp \cdot \mathcal C_k$,
see (\ref{ss:in}).
It depends on $n$ 
but we will suppress that in our notation.
The stability of the pairings $M^{\prime}(n)$ (equation \eqref{e:ps})
says that the following diagram commutes.
\begin{equation*}
\begin{CD}
(\mathcal C_{k})^{\fsp} \oplus \fsp \cdot \mathcal C_k @>M^{\prime}(n)>> (\mathcal G_k^{(r)})^* @<A<\iso< \mathcal G_k^{(r)} \\ 
@VVV @V{\id}VV @V{\id}VV\\
(\mathcal C_{k})^{{\mathfrak{sp}}(2n+2)} \oplus {{\mathfrak{sp}}(2n+2)} \cdot \mathcal C_k 
@>>M^{\prime}(n+1)> (\mathcal G_{k}^{(r)})^* @<\iso<A< \mathcal G_{k}^{(r)}.\\ 
\end{CD}
\end{equation*}
The first vertical map does not restrict to a map
on the invariants.
However, by projecting on the first factor,
we obtain maps
$(\mathcal C_{k})^{\fsp} \to (\mathcal C_{k})^{{{\mathfrak{sp}}(2n+2)}}$,
one for each $n$.
And these are the maps that one can use 
for computation of stable homology
$\underset{n \to \infty}{\lim} H_k(QA_n)$.
This gives us the diagram
\begin{equation*}
\begin{CD}
\mathcal G_k \iso (\mathcal C_{k})^{\fsp} @>M^{\prime}(n)>> (\mathcal G_k^{(r)})^* @<A<\iso< \mathcal G_k^{(r)} \\ 
@VVV @V{\id}VV @V{\id}VV\\
\mathcal G_k \iso (\mathcal C_{k})^{{\mathfrak{sp}}(2n+2)}
@>>M^{\prime}(n+1)> (\mathcal G_{k}^{(r)})^* @<\iso<A< \mathcal G_{k}^{(r)}.\\ 
\end{CD}
\end{equation*}
By the Corollary to Proposition \ref{p:in},
the map
$\fsp \cdot \mathcal C_k \to (\mathcal G_k)^*$ is zero.
Hence the above diagram still commutes.
Furthermore, if one thinks of the map
$(\mathcal C_{k})^{\fsp} \to (\mathcal C_{k})^{{{\mathfrak{sp}}(2n+2)}}$
in terms of graphs,
then it is clear that it preserves the Euler characteristic.
This gives a map
$\mathcal G_k^{(r)} \to \mathcal G_k^{(r)}$ for each $n$.
The direct limit
$\underset{n \to \infty}{\lim} H_k(\mathcal G^{(r)},\partial_n)$
is taken with respect to these maps.
Consider the stability commutative diagram
\begin{equation*}
\begin{CD}
\mathcal G_k^{(r)} @>>> (\mathcal C_{k})^{\fsp} @>M^{\prime}(n)>> (\mathcal G_k^{(r)})^* @<A<\iso< \mathcal G_k^{(r)} \\ 
@VVV @VVV @V{\id}VV @V{\id}VV\\
\mathcal G_{k}^{(r)} @>>> (\mathcal C_{k})^{\fsp+2} @>>M^{\prime}(n)> (\mathcal G_{k}^{(r)})^* @<\iso<A< \mathcal G_{k}^{(r)}.\\ 
\end{CD}
\end{equation*}
We know that the map
$M(n): \mathcal G_k^{(r)} \to (\mathcal G_k^{(r)})^*$
is an isomorphism for all $n$ large enough (\ref{ss:pnd}).
In the above diagram,
it is the composite of the first two horizontal maps.
Proposition \ref{p:sta} now follows
by the isomorphism between the leftmost and rightmost columns.

\begin{remark}
In \cite{\god}, Kontsevich conjectured that
the stable homology groups
$H_k(QA_\infty)$ are finite dimensional.
In terms of graphs, this says that
$H_k(\mathcal G^{(r)},\partial_E) = 0$ 
for fixed $k$ and sufficiently large $r$.
\end{remark}

\section{Proof of the main theorem-Part III}  \label{s:con}

The limit Lie algebra,
$QA_\infty = \underset{n \to \infty}{\lim} QA_n$,
has the structure of a Hopf algebra on its homology
(\ref{ss:ha}).
Let
$PH_*(QA_{\infty})$
denote the subspace of primitive elements of
$H_*(QA_{\infty})$.
Also let $(\mathcal C,\partial_E)$ be the subcomplex of
$(\mathcal G,\partial_E)$ spanned by oriented connected graphs.
We will prove the following theorem.

\begin{theorem}  \label{t:con}
$PH_*(QA_{\infty}) = H_*(\mathcal C,\partial_E).$
\end{theorem}
\begin{proof}
The Hopf algebra $H_*(QA_{\infty})$ is commutative
and cocommutative in the graded sense.
Hence by the structure theorem of Milnor and Moore
\cite{\mm},
$H_*(QA_{\infty})$ is a graded polynomial algebra
generated by the primitive elements.
By Theorem \ref{t:sta},
we know $H_*(QA_{\infty}) = H_*(\mathcal G,\partial_E).$
This along with the identification given by Proposition \ref{p:du}
(which we are going to prove)
shows that the subspace of primitive elements 
is precisely $H_*(\mathcal C,\partial_E).$
\end{proof}
\noindent
As always,
we apply the above theorem to the unit species $uu$
and obtain the following corollary.

\begin{corollary}
$PH_*({{{\mathfrak{sp}}(2 \infty)}}) = H_*(\mathcal P,\partial_E).$
\end{corollary}

\begin{proposition}  \label{p:du}
The product in the stable homology
$H_*(QA_{\infty})$ can be defined on the chain complex
$(\mathcal G,\partial_E)$ as disjoint union of graphs.
\end{proposition}
\begin{proof}
We make use of the deformation map
$D(n) : \Lambda^k(QA_n) \to \mathcal G_k^{(r)}$.
This is a composition of the maps
$\Lambda^k(QA_n) \to (\mathcal G_k^{(r)})^* \to \mathcal G_k^{(r)}$
that we had in the previous section (\ref{ss:cd}).
It is an extension of the deformation map
$D(n) : \mathcal G_k^{(r)} \to \mathcal G_k^{(r)}$ 
defined at the end of (\ref{ss:dmg}).
For simplicity,
we will suppress the dependence on $n$
and use the notation
$D(\varphi) = \sum_{\Gamma} D(\varphi,\Gamma) \ \Gamma$
for $\varphi \in \Lambda^k(QA_{n})$ and
$\Gamma \in \mathcal G_k^{(r)}$.
The usage of $D$ for $D(n)$ 
is local to this proposition.
Later in Appendix~\ref{ss:dg},
the letter $D$ will be used for
something slightly different.

The map $QA_n \oplus QA_m \into QA_{n+m}$ induces a map
$\Lambda^i(QA_{n}) \otimes \Lambda^{k-i}(QA_{m}) \to
\Lambda^k(QA_{n+m})$.
More explicitly,
$\varphi_1 \otimes \varphi_2 \mapsto \varphi_1 \wedge \varphi_2$, 
with the indices of $\varphi_2$ shifted up by $n$.
In order to complete the proof,
it is enough to show that the following diagram commutes.

\begin{equation*}
\begin{CD}
\Lambda^i(QA_{n}) \otimes \Lambda^{k-i}(QA_{m}) @>D \otimes D>> \mathcal G_i \otimes \mathcal G_{k-i} \\
 @VVV @V \mu VV \\
\Lambda^k(QA_{n+m}) @>D>> \mathcal G_k. \\
\end{CD}
\end{equation*}
Here the map $\mu$ stands for disjoint union of graphs.
The commutativity is equivalent to showing the identity
$$D(\varphi_1 \wedge \varphi_2,\Gamma) = \underset{\Gamma_1 \sqcup \Gamma_2 = \Gamma}{\sum}
D(\varphi_1,\Gamma_1) D(\varphi_2,\Gamma_2).$$
To calculate the LHS,
we match the $i$ pieces of $\varphi_1$ 
and $k-i$ pieces of $\varphi_2$ 
with the $k$ vertices of $\Gamma$ and then contract indices.
And we sum over all matchings.

The key observation is that $\varphi_1$ and $\varphi_2$
must land on disconnected parts of $\Gamma$,
say $\Gamma_1$ and $\Gamma_2$,
in order to get a non-zero contraction.
This is because the indices that occur in $\varphi_1$
are disjoint from those in $\varphi_2$.
The rest of the argument is now fairly clear and we omit it.
\end{proof}

\begin{remark}
In fact, the chain complex $(\mathcal G,\partial_E)$ 
is a differential graded Hopf algebra
with product $\mu$ given by disjoint union and
the coproduct $\Delta$ defined for a connected graph $\Gamma$
by $\Delta (\Gamma) = 1 \otimes \Gamma + \Gamma \otimes 1$
and extended to $\mathcal G$ as a morphism of algebras.
Here 1 stands for the unit in $\mathcal G$ and
may be thought of as the empty graph.
The boundary map $\partial_E$ is a
derivation with respect to the product $\mu$
and a coderivation with respect to the coproduct $\Delta$.
And this induces the Hopf algebra structure on
$H_*(\mathcal G,\partial_E)$. 
\end{remark}

\subsection*{Proof of the main theorem concluded}  

Our main theorem (Theorem \ref{t:mt}) now follows by putting together
Theorem \ref{t:con} and its corollary
and the corollary to Proposition \ref{p:ss}.

\appendix

\section{Deformation quantisation} \label{ss:dq}

For a good review of the problem of
deformation quantisation, 
see the notes by Voronov~\cite{\voronov}.
The main object of interest is a Poisson algebra.
The algebra of functions on a symplectic,
or more generally, a Poisson manifold 
form a Poisson algebra.
This is the classical case.
In (\ref {ss:tp}-\ref {ss:se}),
we provide some background on the classical case.
This would also be useful for some of the ideas 
in Appendix~\ref{ss:dg}.
The material in~\ref{ss:tp} is taken from~\cite{\voronov}.

In (\ref {ss:qoa}),
we speculate on the form of this problem for operads.
The main question is:
What is a ``Poisson operad manifold''?
Unlike the classical case,
a symplectic operad manifold is not automatically
a Poisson operad manifold.
The relation between our viewpoint and 
the standard deformation theory of operads
considered by Balavoine in~\cite{\lsv}
is not clear.

\subsection{The problem}  \label{ss:tp}

Let $A$ be a commutative algebra over a field $k$
of characteristic zero.
A formal deformation of $A$ is a $k[[t]]$ bilinear product
(which we denote $\star$) on the space $A[[t]]$
of formal power series in a variable $t$ satisfying:
$$F \star H = F \cdot H + \mu_1(F,H) t + \mu_2(F,H) t^2 +
\ldots \ \ \text{for} \ \ F,H \in A,$$
where $F \cdot H$ is the original product on $A$
and the star product $\star$ is associative.
The deformation is called \emph{trivial} if
all the higher products $\mu_1,\mu_2,\ldots$ are zero.

Suppose that $\star$ is a deformation of $\cdot$,
the original product.
Then define $\{F,H\} = \frac{1}{2} (\mu_1(F,H) - \mu_1(H,F))$.
One can check that $\{\ ,\ \}$ 
is a Lie bracket on $A$ 
and further that the triple $(A,\cdot,\{\ ,\ \})$ is a Poisson algebra.
Recall that a Poisson algebra is a space
with compatible commutative and Lie structures.
The compatibility relation is 
$\{F \cdot G,H\} = F \cdot \{G,H\} + G \cdot \{F,H\}$.

In physical terms, one regards
the Poisson algebra $A$ as the quasi-classical limit 
of the associative algebra $A[[t]]$, 
and the algebra $A[[t]]$ as a deformation quantisation
of the Poisson algebra $A$.
The deformation quantisation problem is the inverse problem:
given a Poisson algebra $(A,\cdot,\{\ ,\ \})$,
find a formal deformation $\star$ returning the original
Poisson algebra structure on $A$ in the quasi-classical limit.
 
\subsection{Gauge equivalence}  \label{ss:ge}

There is a natural gauge group acting on star products.
This group consists of automorphisms $D$ of $A[[t]]$
which are $k[[t]]$ linear.
They have the form
$D = D_0 + t D_1 + t^2 D_2 + \ldots$,
where $D_i : A \to A$ are operators 
with $D_0$ being the identity.
This means that
$$F \mapsto F + t D_1(F) + t^2 D_2(F) + \ldots \ \ \text{for}
\qquad F \in A,$$
and for a general element in $A[[t]]$,
one uses $k[[t]]$ linearity.

Two star products $\star$ and $\star^{\prime}$
are gauge equivalent if there is an automorphism $D$
as above so that the following diagram commutes.
\begin{equation*}
\begin{CD}
A[[t]] \otimes A[[t]] @>D \otimes D>> A[[t]] \otimes A[[t]] \\
 @V \star VV @V \star^{\prime} VV \\
A[[t]] @>D>> A[[t]]. \\
\end{CD}
\end{equation*}

\subsection{The simplest example}  \label{ss:se}

Let $A$ be the space of all polynomial functions on $\R^{2n}$,
i.e. polynomials in the $2n$ variables 
$p_1, \ldots,p_{n},q_1, \ldots,q_{n}$.
Then the triple
$(A,\cdot,\{\ ,\ \})$ is a Poisson algebra,
where $\cdot$ is the usual product in $A$
and $\{\ ,\ \}$ is the Poisson bracket
defined by equation~\eqref{e:pb} in (\ref{s:sog}).

In this case, the solution to the
deformation quantisation problem is 
given by the Moyal $\star$ product.
In fact, we will have
$F \star H = F \cdot H + \{F,H\} t + \ \text{higher terms}$;
that is, $\mu_1(F,H) = \{F,H\}$.

\subsubsection{The Moyal product}  \label{ss:mp}

Let $B$ be any operator on $A \otimes A$.
For $F,H \in A$, define
$B(F,H)$ to be the element of $A$ obtained
by applying the product in $A$ to 
$B(F \otimes H)$.
Namely, if
$B = \sum B_{(1)} \otimes B_{(2)}$ then
$B(F,H) = \sum B_{(1)} F \cdot B_{(2)} H$.

Now let 
$B = \sum_{i=1}^n \frac{\partial}{\partial p_i} \otimes \frac{\partial}{\partial q_i} - 
\frac{\partial}{\partial q_i} \otimes \frac{\partial}{\partial p_i}$
be a bi-differential operator on $A \otimes A$.
By our notation,
$$B(F,H)=\sum_{i=1}^n \frac{\partial F}{\partial p_i} \frac{\partial H}{\partial q_i} - 
\frac{\partial F}{\partial q_i} \frac{\partial H}{\partial p_i}.$$  
In other words, $B(F,H) = \{F,H\}$ is just the Poisson bracket.
Now define the Moyal $\star$ product by
\begin{equation}  \label{e:mp}
F \star H = e^{tB}(F,H) = \sum_{n \geq 0} t^n
\frac{B^n}{n!}(F,H).
\end{equation}
Observe that
$F \star H = F \cdot H + \{F,H\} t + \ \text{higher terms}$,
as claimed earlier.
To show that the Moyal product 
solves the deformation quantisation problem,
one must prove that it is associative.
We will do this by interpreting equation~\eqref{e:mp} using pictures.
  
\subsubsection{Associativity of the Moyal product}  \label{ss:amp}

If $F \in A$ is a monomial, say
$F = p_1^2 p_2 q_2$,
then we represent it as 
$F = 
\fbox{
\begin{minipage}{.4 in}
\begin{center}
\psfrag{a}{\Huge $p_1$}
\psfrag{b}{\Huge $p_1$}
\psfrag{c}{\Huge $p_2$}
\psfrag{d}{\Huge $q_2$}
\resizebox{.4 in}{.4 in}{\includegraphics{ccd.eps}}
\end{center}
\end{minipage} 
}
$.
And if $F$ is a polynomial rather than a monomial
then we represent it as a formal sum of pictures.
It is clear that the product of two monomials in $A$ is defined
as the disjoint union of pictures.

Our goal is to understand the operator $\frac{B^n}{n!}$
via pictures.
As a start, we understand the term
$\frac{\partial F}{\partial p_i} \frac{\partial H}{\partial q_i}$.
This is a sum of pictures of the form
\begin{minipage}{1 in}
\begin{center}
\psfrag{pi}{\Huge $p_i$}
\psfrag{qi}{\Huge $q_i$}
\psfrag{F}{\Huge $F$}
\psfrag{H}{\Huge $H$}
\resizebox{1 in}{.4 in}{\includegraphics{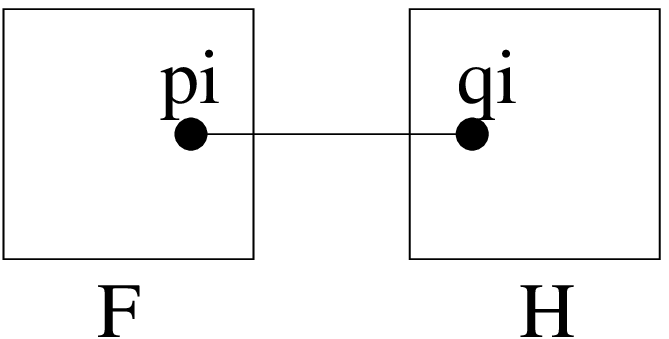}}
\end{center}
\end{minipage} 
,
where the edge goes from a point labelled $p_i$
to a point labelled $q_i$.
The other points in $F$ and $H$ are not shown in the picture.
We will have the following interpretation.
A picture of the form
\begin{minipage}{1 in}
\begin{center}
\psfrag{pi}{\Huge $$}
\psfrag{qi}{\Huge $$}
\psfrag{F}{\Huge $$}
\psfrag{H}{\Huge $$}
\resizebox{1 in}{.4 in}{\includegraphics{cm.eps}}
\end{center}
\end{minipage} 
means that we delete the edge and its endpoints
and then take disjoint union of the two parts.
Earlier in the paper,
such an operation was called a mating.

Now we claim that
\begin{equation}  \label{e:mpp}
\frac{B^n}{n!}(F,H) = 
\sum_P (-1)^{\sign(P)}
\left(
\begin{minipage}{1.2 in}
\begin{center}
\psfrag{F}{\Huge $F$}
\psfrag{H}{\Huge $H$}
\psfrag{n}{\Huge $n$}
\resizebox{1 in}{.4 in}{\includegraphics{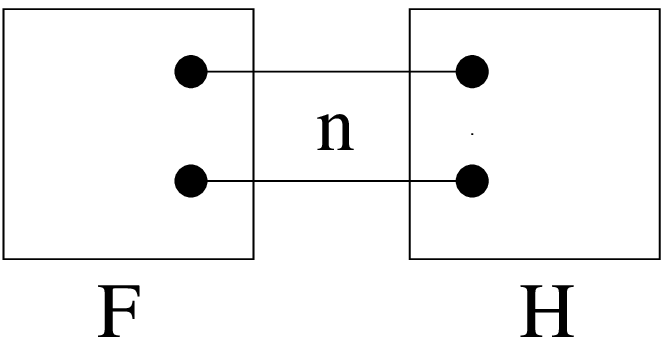}}
\end{center}
\end{minipage} 
\right)
,
\end{equation}
where $P$ is a picture as shown.
There are $n$ edges between $F$ and $H$ 
each connecting a $p_i$ in $F$ with a $q_i$ in $H$ or vice-versa.
And $\sign(P)$ is the number of edges 
that are out of order;
i.e. they connect
a $q_i$ in $F$ to a $p_i$ in $H$.  
The interpretation for $P$ is that 
we delete the $n$ edges and their endpoints
and take disjoint union of what remains.

Note that the factor $n!$ is no longer necessary
in the pictorial description.
Each of the $n$ $B$'s in $B^n$ contributes
to one of the $n$ edges in the picture.
So there are $n!$ ways to obtain the same picture.
After dividing by $n!$,
we get every picture exactly once,
giving us equation~\eqref{e:mpp}.

The associativity of the Moyal product can be seen from
the formula
\begin{equation}  \label{e:amp}
F \star G \star H = 
\sum_{\substack {n \geq 0 \\ n=k+l+m}} t^n
\sum_P (-1)^{\sign(P)}
\left(
\begin{minipage}{1.2 in}
\begin{center}
\psfrag{F}{\Huge $F$}
\psfrag{G}{\Huge $G$}
\psfrag{H}{\Huge $H$}
\psfrag{k}{\Huge $k$}
\psfrag{l}{\Huge $l$}
\psfrag{m}{\Huge $m$}
\psfrag{n}{\Huge $n$}
\resizebox{1 in}{.8 in}{\includegraphics{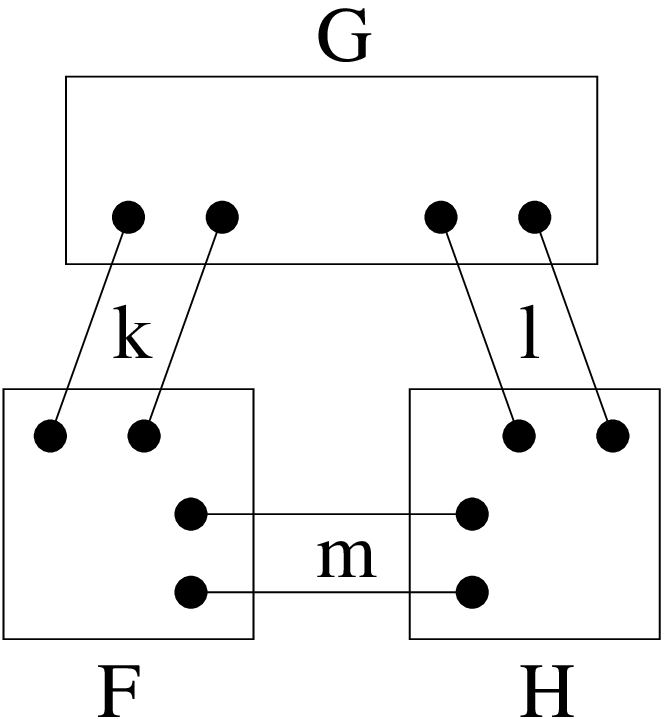}}
\end{center}
\end{minipage} 
\right)
,
\end{equation}
where the picture $P$ means that
delete all edges and their endpoints and
then take disjoint union of what remains.
And  $\sign(P)$ is the number of edges that connect
a $q_i$ to a $p_i$. 
(The edges go from $F$ to $G$, $G$ to $H$ 
and $F$ to $H$.)
The letters $k$, $l$ and $m$ refer to the number of connecting edges.

\subsection{An example based on the surface species (\ref{ss:ss})}  \label{ss:ae}

Before going to the general case,
let us do another example.
Consider a picture of the form
$ 
\fbox{
\begin{minipage}{.6 in}
\begin{center}
\psfrag{a}{\Huge $p_1$}
\psfrag{b}{\Huge $p_2$}
\resizebox{.6 in}{.4 in}{\includegraphics{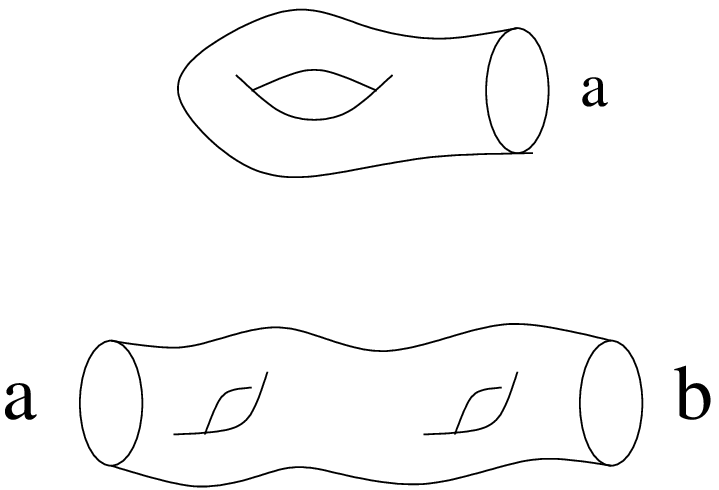}}
\end{center}
\end{minipage} 
}
$,
namely, we have a compact orientable surface
(not necessarily connected)
whose boundary circles are labelled by
$p_1,\ldots,p_n,q_1,\ldots,q_n$.
Let $A$ be the vector space spanned by such pictures,
each picture being a monomial.
The space $A$ is a Poisson algebra 
with the commutative product $\cdot$ given by disjoint union of pictures
and the Lie bracket $\{\ ,\ \}$ given by
\begin{equation}  \label{e:mps}
\{F,H\} = 
\sum_P (-1)^{\sign(P)}
\left(
\begin{minipage}{1 in}
\begin{center}
\psfrag{pi}{\Huge $$}
\psfrag{qi}{\Huge $$}
\psfrag{F}{\Huge $F$}
\psfrag{H}{\Huge $H$}
\resizebox{1 in}{.4 in}{\includegraphics{cm.eps}}
\end{center}
\end{minipage} 
\right)
,
\end{equation}
where the edge connects a boundary circle 
labelled $p_i$ in $F$ to a boundary circle 
labelled $q_i$ in $H$ or vice-versa.
The picture means that the two boundary circles
joined by the edge are glued together
to form a single surface.
As an example,
\[
\left\{
\begin{minipage}{.8 in}
\begin{center}
\psfrag{a}{\Huge $p_1$}
\resizebox{.6 in}{.2 in}{\includegraphics{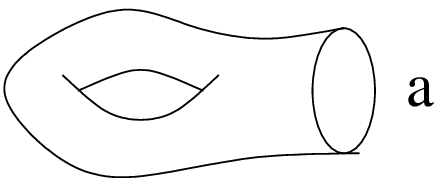}}
\end{center}
\end{minipage}
\ \ , \ \
\begin{minipage}{.8 in}
\begin{center}
\psfrag{a}{\Huge $q_1$}
\psfrag{b}{\Huge $p_2$}
\resizebox{.7 in}{.17 in}{\includegraphics{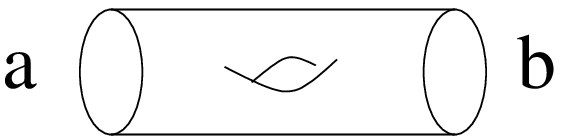}}
\end{center}
\end{minipage}
\right\}
\ \ = \ \
\begin{minipage}{.8 in}
\begin{center}
\psfrag{b}{\Huge $p_2$}
\resizebox{.7 in}{.2 in}{\includegraphics{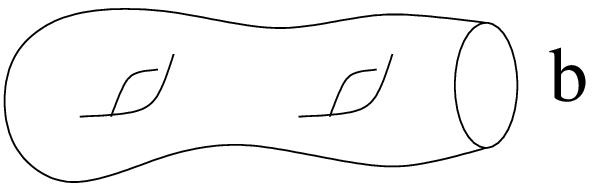}}
\end{center}
\end{minipage}
.
\]
The deformation quantisation for the Poisson algebra $A$
is again given by equations~\eqref{e:mp} and~\eqref{e:mpp}.
The picture in equation~\eqref{e:mpp} means that
we glue together the boundary circles that are
connected by edges.
There will be a total of $n$ gluings.
The associativity of the star product follows
from equation~\eqref{e:amp}.

\subsection{Quantization of operad algebras}  \label{ss:qoa}

We now initiate a general discussion of deformation quantization
for operads.
Precise definitions are not given;
so this section should be read mainly for the philosophy.

As a candidate for the Poisson algebra $A$,
take the free algebra over a mated species $Q$,
which we have previously denoted $QA$, see (\ref{ss:fa}).
On the geometric side,
we are now dealing with
a symplectic operad manifold.
The examples in (\ref{ss:se}) and (\ref{ss:ae})
correspond to the commutative and surface species
respectively.

\begin{remark}
By our construction, the vector space $QA$ is graded
beginning with degree $2$.
For example,
$ccA$ are polynomials with no constant or linear terms.
Also the monomials in $ssA$ have at least two boundary circles.
However, for the deformation quantization problem,
one needs to enlarge $QA$ suitably 
by adding stuff in degrees 0 and 1,
as was done for the previous two examples.
\end{remark}

We first want to give a Poisson structure to $QA$
and then deform it by the Moyal product.
We know that $QA$ has a Lie structure $\{\ ,\ \}$.
This was discussed in (\ref{ss:mate}).
In terms of pictures, it is given by equation~\eqref{e:mps}.
The interpretation of the picture is that
a mating occurs along the edge.
But what is the commutative structure on $QA$ ?
We have emphasised that an element of $QA$ 
can be represented by a picture.
So take $F \cdot H$ to be the
``disjoint union of pictures''
that represent $F$ and $H$.
The implicit assumption is that
the species $Q$ can be written as the exponential
of some other species.
This is our proposal for the Poisson structure on $QA$.
From the viewpoint of pictures,
the compatibility of the commutative and Lie products is clear.

Next we want to define the Moyal product
and show that it solves the deformation quantisation problem.
In other words, we want to understand the meaning of
$\frac{B^n}{n!}(F,H)$ for $n \geq 2$,
where $B$ is the bi-differential operator
defined in (\ref{ss:mp}).
For $n=0,1$, we know that 
$B^0(F,H)=F \cdot H$ and $B(F,H) = \{F,H\}$.
Writing
$\frac{B^n}{n!} = \sum B_{(1)} \otimes B_{(2)}$,
where $B_{(1)}$ and $B_{(2)}$ are $n$th order differential operators,
we obtain
$\frac{B^n}{n!}(F,H) = \sum B_{(1)} F \otimes B_{(2)} H$.
One needs to make sense of this.
We first explain how a higher order derivative works
by showing it on a schematic example.
\begin{equation}  \label{e:hd}
\begin{minipage}{.8 in}
\begin{center}
\psfrag{a}{\Huge $x_2$}
\psfrag{b}{\Huge $x_1$}
\psfrag{c}{\Huge $x_1$}
\psfrag{d}{\Huge $x_4$}
\resizebox{.6 in}{.6 in}{\includegraphics{q.eps}}
\end{center}
\end{minipage}
\qquad \buildrel \frac{\partial^2}{\partial x_1^2}\over \longmapsto \qquad
2 \left(
\begin{minipage}{.8 in}
\begin{center}
\psfrag{x1}{\Huge $x_1$}
\psfrag{x2}{\Huge $x_2$}
\psfrag{x4}{\Huge $x_4$}
\resizebox{.7 in}{.6 in}{\includegraphics{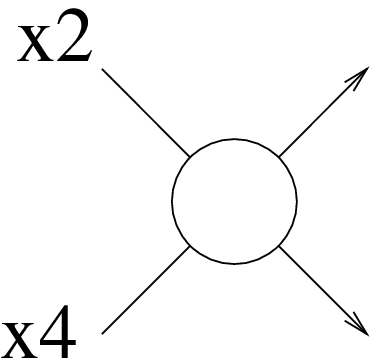}}
\end{center}
\end{minipage}
\right).
\end{equation}
Thus a $n$th order differential operator acting on $F$
produces an algebraic object with $n$ outputs
(\ref{ss:prop}).
We now interpret the tensor sign as a $n$th order mating.
We write
$$
\frac{B^n}{n!}(F,H) = 
\sum_P (-1)^{\sign(P)}
\left(
\begin{minipage}{1.2 in}
\begin{center}
\psfrag{F}{\Huge $F$}
\psfrag{H}{\Huge $H$}
\psfrag{n}{\Huge $n$}
\resizebox{1 in}{.4 in}{\includegraphics{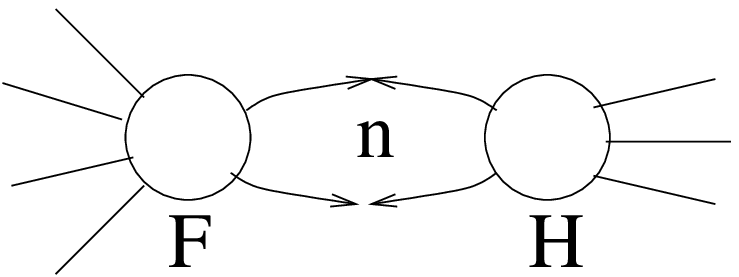}}
\end{center}
\end{minipage} 
\right)
,$$
where the picture $P$ shows a higher order mating.
There are $n$ ideal edges between $F$ and $H$ 
each connecting a $p_i$ with a $q_i$ or vice-versa.
And $\sign(P)$ is the number of edges that connect
a $q_i$ on $F$ to a $p_i$ on $H$.
For the definition of an ideal edge, see (\ref{ss:mf}).

The examples discussed in Section~\ref{s:me}
have associated PROPs.
In those cases, 
we understand the meaning of higher order derivatives,
i.e. the right hand side of equation~\eqref{e:hd}
makes sense.
Similarly, higher order matings 
have a natural meaning in those examples.
The example of the surface species was illustrated in
(\ref{ss:ae}).
The deformation quantisation for the Poisson algebra $A$
is again given by the formula in equation~\eqref{e:mp}.
The Moyal $\star$ product is associative
in these cases for the same reason 
as in the commutative case $(Q=cc)$.
The associativity can be seen from the formula
$$
F \star G \star H = 
\sum_{\substack {n \geq 0 \\ n=k+l+m}} t^n
\sum_P (-1)^{\sign(P)}
\left(
\begin{minipage}{1.2 in}
\begin{center}
\psfrag{F}{\Huge $F$}
\psfrag{G}{\Huge $G$}
\psfrag{H}{\Huge $H$}
\psfrag{k}{\Huge $k$}
\psfrag{l}{\Huge $l$}
\psfrag{m}{\Huge $m$}
\psfrag{n}{\Huge $n$}
\resizebox{1 in}{.8 in}{\includegraphics{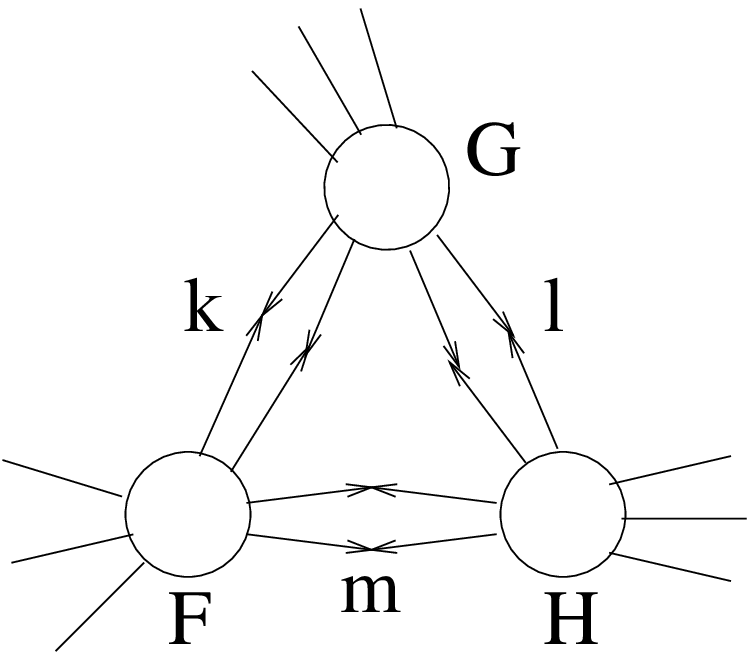}}
\end{center}
\end{minipage} 
\right)
,$$
\noindent
where $\sign(P)$ is the number of ideal edges that connect
a $q_i$ to a $p_i$. 
(The edges go from $F$ to $G$, $G$ to $H$ 
and $F$ to $H$.)
The letters $k,l,m$ stand for the number of edges.

The interested reader can work out the details 
for the examples
based on graphs (\ref{ss:g}-\ref{ss:gt}).
Another nice example is obtained 
by exponentiating the associative species.

\section{The deformation map on graphs}  \label{ss:dg}

Recall that the chain complex $(\mathcal G,\partial_E)$ 
is a differential graded Hopf algebra
with product $\mu$ given by disjoint union.
For a definition of $(\mathcal G,\partial_E)$,
 see Section~\ref{s:gh}
and for the Hopf algebra structure,
see the remark at the end of Section~\ref{s:con}.
In addition, recall that we defined 
another boundary operator $\partial_H$ on $\mathcal G$
that satisfied
$\partial_E \partial_H = - \partial_H \partial_E$,
see (\ref{ss:finb}).
So $\partial_H$ induces a map on homology
$H_*(\mathcal G,\partial_E) \to H_*(\mathcal G,\partial_E)$. 

The boundary map $\partial_E$ is a
derivation with respect to the product $\mu$.
This can be written as $\partial_E \mu = \mu \partial_E$.
However $\partial_H$ is not a derivation with respect to $\mu$.
Following~\cite{\cv},
define a bracket
$[\ ,\ ] : \mathcal G \otimes \mathcal G \to \mathcal G$ by
\[
[\ ,\ ] = \partial_H \mu - \mu \partial_H.
\]
In other words,
$[\ ,\ ]$ measures the failure of $\partial_H$
to be a derivation with respect to $\mu$.
As a formal consequence of the previous three identities,
it follows that
$\partial_E [\ ,\ ] = - [\ ,\ ] \partial_E$.
This gives us an induced map on homology
$[\ ,\ ] : H_*(\mathcal G,\partial_E) \otimes H_*(\mathcal G,\partial_E)
\to H_*(\mathcal G,\partial_E)$.

In this section, we show the following.

\begin{proposition}  \label{p:triv}
The maps $\partial_H$ and $[\ ,\ ]$ are trivial on homology.
\end{proposition}

This answers some of the questions raised in~\cite{\cv}.
The proof is largely a matter of putting together
things we already know.
The main step is to use 
the pairing $M(n)$ on $\mathcal G$ 
and its adjoint property 
proved in Proposition~\ref{p:ad}. 

The coproduct on $\mathcal G$ is irrelevant to Proposition~\ref{p:triv}.
So we write the differential algebra $\mathcal G$ as a triple
$(\mathcal G,\partial_E,\mu)$.
Extending $\partial_E$ and $\mu$ by $\Q[[t]]$ linearity,
we get a differential algebra
$(\mathcal G[[t]],\partial_E,\mu)$. 
In other words, this is just the
trivial deformation of the algebra $\mathcal G$,
see (\ref{ss:tp}).
We now consider another deformation of $\mathcal G$ that is
gauge equivalent to this one (\ref{ss:ge}).

\subsection{The deformation map}

Let $D(n) : \mathcal G \to \mathcal G$ be the deformation map
defined at the end of (\ref{ss:pnd}).
It was defined by $D(n) = A^{-1} \circ M(n)$,
where the map $M(n)$ was induced by a pairing on $\mathcal G$
(equation \eqref{e:pg}),
while the map $A$ scaled a graph $\Gamma$ by 
the factor $|\Aut(\Gamma)|$.

Recall that $D(n) \Gamma$
is a polynomial in $2n$ of degree $e$
with coefficients in $\mathcal G$.
Furthermore, the coefficient of $(2n)^e$ is exactly $\Gamma$.
Here $e=|E(\Gamma)|$, the number of edges in $\Gamma$.
Write
$D(n) \Gamma = (2n)^e (D_0 + \frac{1}{2n} D_1 + \ldots)(\Gamma)$.
This defines maps $D_i : \mathcal G \to \mathcal G$ 
which do not depend on $n$. Also $D_0$ is the identity.

Consider the $\Q[[t]]$ linear automorphism of $\mathcal G[[t]]$
given by $D = D_0 + t D_1 + t^2 D_2 + \ldots$,
where $D_i : \mathcal G \to \mathcal G$ are as above.
The map $D$ defines an element of the gauge group
(\ref{ss:ge}).

\subsection{Comparing the two gauge equivalent situations}

Let $(\mathcal G[[t]],\partial_t,\mu_t)$
be the differential algebra got by applying $D$
to the trivial deformation $(\mathcal G[[t]],\partial_E,\mu)$. 
Stated differently, one has two commutative diagrams. 
\begin{equation}  \label{e:dcdiag}
\begin{CD}
\mathcal G[[t]] @>D>> \mathcal G[[t]]\\
 @V \partial_E VV @V \partial_t VV \\
\mathcal G[[t]] @>D>> \mathcal G[[t]]. \\
\end{CD}
\hspace{.6 in}
\begin{CD}
\mathcal G[[t]] \otimes \mathcal G[[t]] @>D \otimes D>> \mathcal G[[t]] \otimes \mathcal G[[t]] \\
 @V \mu VV @V \mu_t VV \\
\mathcal G[[t]] @>D>> \mathcal G[[t]]. \\
\end{CD}
\end{equation}
Since $\partial_E$ is a
derivation with respect to the product $\mu$,
we have $\partial_E \mu = \mu \partial_E$.
Hence
we know $\partial_t \mu_t = \mu_t \partial_t$.
Write
$\partial_t = \partial_0 + t \partial_1 + t^2 \partial_2 + \ldots$
and
$\mu_t = \mu_0 + t \mu_1 + t^2 \mu_2 + \ldots$.
Then the commutativity of the diagrams in~\eqref{e:dcdiag}
implies
$\partial_0 = \partial_E$ and $\mu_0 = \mu$.

\begin{lemma} \label{l:drel}
The following relations hold.
\begin{itemize}
\item[(1)]
$- \partial_1 = - D_1 \partial_0 + \partial_0 D_1$.

\item[(2)]
$\mu_1 = D_1 \mu_0 - \mu_0 D_1$.

\item[(3)]
$- \partial_1 \mu_0 + \mu_0 \partial_1 = \partial_0 \mu_1 - \mu_1 \partial_0$.
\end{itemize}
\end{lemma}
\begin{proof}
The first two items follow by looking at the coefficient of $t$
in the diagrams in~\eqref{e:dcdiag}
and the third item follows by looking at the coefficient of $t$
in $\partial_t \mu_t = \mu_t \partial_t$.
\end{proof}

\begin{lemma}
$\partial_1 = - \partial_H$.
\end{lemma}
\begin{proof}
The commutative diagram~\eqref{e:cdg} in Section~\ref{s:mgh}
can be redrawn as under
\begin{equation*} 
\begin{CD}
\mathcal G[[t]] @>D>> \mathcal G[[t]]\\
 @V \partial_E + t \partial_H VV @V \partial_E VV \\
\mathcal G[[t]] @>D>> \mathcal G[[t]]. \\
\end{CD}
\end{equation*}
Comparing the coefficient of $t$ gives
$\partial_H = - D_1 \partial_0 + \partial_0 D_1$.
This together with Lemma~\ref{l:drel} (item (1))
shows that $\partial_1 = - \partial_H$.
\end{proof}

\subsection*{Proof of Proposition~\ref{p:triv}}
Items $(1)$ and $(3)$ in Lemma~\ref{l:drel}
can now be rewritten as 
$\partial_H = - D_1 \partial_E + \partial_E D_1$ and
$[\ ,\ ] = \partial_E \mu_1 - \mu_1 \partial_E$.
Thus $\partial_H$ and $[\ ,\ ]$ induce the zero map on homology,
with $D_1$ and $\mu_1$ providing the respective chain homotopies.
This completes the proof of Proposition~\ref{p:triv}.

\begin{remark}
It is possible to just give an explicit definition of $\mu_1$
and check the relation $[\ ,\ ] = \partial_E \mu_1 - \mu_1 \partial_E$
directly.
This will give us that $[\ ,\ ]$ is zero on homology.
And similarly for $D_1$.
But we prefer the more conceptual approach 
via deformation theory.
It would be interesting to also describe explicitly
the higher products
$\mu_2, \mu_3, \ldots$ and
$\partial_2, \partial_3, \ldots$, etc.

Recall that $\mathcal G$ also has a coproduct $\Delta$ and
$\partial_E$ is a coderivation with respect to the coproduct $\Delta$.
However, $\partial_H$ is not and
this failure can be measured by a cobracket
$\theta : \mathcal G \to \mathcal G \otimes \mathcal G$ given by
$\theta = \partial_H \Delta - \Delta \partial_H$.
Similar arguments show that 
$\theta$ also induces the zero map on homology.
\end{remark}

\subsection*{Acknowledgements}
I thank Ken Brown for being a constant source of
knowledge and encouragement.
Among the participants of the Bernstein seminar
at Cornell University (Fall 2000),
Dan Ciubotaru, Ferenc Gerlits and Jim Conant
deserve a special thank you. 
I also thank M. Aguiar, Y. Berest and M. Yakimov
for helpful comments.
Finally, I would like to thank Kontsevich
without whom this paper would never have been written.
I hope that the reader will go back 
and read his original papers,
for which there can be no substitute.
It is fitting to conclude with a Kontsevich sentence.
We choose the opening sentence of~\cite{\kont}.
\begin{quotation}
We shall describe a program here relating Feynman diagrams,
topology of manifolds, homotopical algebra,
non-commutative geometry and
several kinds of ``topological physics.''
\end{quotation}
Enjoy!

%\input e

%\input app

%\bibliographystyle{hamsplain}
%\bibliography{quantum}

%\input a.bbl
\def\cprime{$'$} \def\cprime{$'$} \def\cprime{$'$}
\providecommand{\bysame}{\leavevmode\hbox to3em{\hrulefill}\thinspace}

\end{document}